\documentclass[10pt]{amsart}
\usepackage{amsmath}
\usepackage{amsxtra}
\usepackage{amscd}
\usepackage{amsthm}
\usepackage{amsfonts}
\usepackage{amssymb}
\usepackage{eucal}

\newtheorem{cor}[subsection]{Corollary}
\newtheorem{lem}[subsection]{Lemma}
\newtheorem{prop}[subsection]{Proposition}

\newtheorem{conj}[subsection]{Conjecture}
\newtheorem{thm}[subsection]{Theorem}

\theoremstyle{definition}

\theoremstyle{remark}

\newcommand{\thmref}[1]{Theorem~\ref{#1}}
\newcommand{\secref}[1]{Sect.~\ref{#1}}
\newcommand{\lemref}[1]{Lemma~\ref{#1}}
\newcommand{\propref}[1]{Proposition~\ref{#1}}
\newcommand{\corref}[1]{Corollary~\ref{#1}}

\numberwithin{equation}{section}

\newcommand{\nc}{\newcommand}
\nc{\renc}{\renewcommand}
\nc{\ssec}{\subsection}
\nc{\sssec}{\subsubsection}
\nc{\on}{\operatorname}

\nc\ol{\overline}
\nc\wt{\widetilde}
\nc\tboxtimes{\wt{\boxtimes}}
\nc{\alp}{\alpha}

\nc{\ZZ}{{\mathbb Z}}
\nc{\NN}{{\mathbb N}}
\nc{\OO}{{\mathbb O}}
\renc{\SS}{{\mathbb S}}
\nc{\DD}{{\mathbb D}}
\nc{\GG}{{\mathbb G}}

\nc{\Fq}{{\mathbb F}_q}
\nc{\Fqb}{\ol{{\mathbb F}_q}}
\nc{\Ql}{\ol{{\mathbb Q}_\ell}}
\nc{\id}{\text{id}}
\nc\X{\mathcal X}

\nc{\Hom}{\on{Hom}}
\nc{\Lie}{\on{Lie}}
\nc{\Loc}{\on{Loc}}
\nc{\Pic}{\on{Pic}}
\nc{\Bun}{\on{Bun}}
\nc{\IC}{\on{IC}}
\nc{\Aut}{\on{Aut}}
\nc{\rk}{\on{rk}}
\nc{\Sh}{\on{Sh}}
\nc{\Perv}{\on{Perv}}
\nc{\pos}{{\on{pos}}}
\nc{\Conv}{\on{Conv}}
\nc{\Sph}{\on{Sph}}
\nc{\Sym}{\on{Sym}}
\nc{\BunBb}{\overline{\Bun}_B}
\nc{\Buno}{\overset{o}{\Bun}}
\nc{\BunPb}{{\overline{\Bun}_P}}
\nc{\BunBM}{\Bun_{B(M)}}
\nc{\BunBMb}{\overline{\Bun}_{B(M)}}
\nc{\BunPbw}{{\widetilde{\Bun}_P}}
\nc{\BunBP}{\widetilde{\Bun}_{B,P}}
\nc{\GUb}{\overline{G/U}}
\nc{\GUPb}{\overline{G/U(P)}}

\nc{\Hhom}{\underline{\on{Hom}}}
\nc\syminfty{\on{Sym}^{\infty}}
\nc\lal{\ol{\lambda}}
\nc\xl{\ol{x}}
\nc\thl{\ol{\theta}}
\nc\nul{\ol{\nu}}
\nc\mul{\ol{\mu}}
\nc\Sum\Sigma
\nc{\oX}{\overset{o}{X}{}}
\nc{\hl}{\overset{\leftarrow}h{}}
\nc{\hr}{\overset{\rightarrow}h{}}
\nc{\M}{{\mathcal M}}
\nc{\N}{{\mathcal N}}
\nc{\F}{{\mathcal F}}
\nc{\D}{{\mathcal D}}
\nc{\Q}{{\mathcal Q}}
\nc{\Y}{{\mathcal Y}}
\nc{\G}{{\mathcal G}}
\nc{\E}{{\mathcal E}}
\nc{\CalC}{{\mathcal C}}
\nc\Dh{\widehat{\D}}

\nc{\C}{{\mathcal C}}
\nc{\K}{{\mathcal K}}
\renewcommand{\H}{{\mathcal H}}

\nc{\T}{{\mathcal T}}
\nc{\V}{{\mathcal V}}
\renc{\P}{{\mathcal P}}
\nc{\A}{{\mathcal A}}
\nc{\B}{{\mathcal B}}
\nc{\U}{{\mathcal U}}

\nc{\Gr}{{\on{Gr}}}

\nc{\frn}{{\check{\mathfrak u}(P)}}
\nc{\p}{\mathfrak p}
\nc{\q}{\mathfrak q}
\nc\f{{\mathfrak f}}

\nc{\qo}{{\mathfrak q}}
\nc{\po}{{\mathfrak p}}
\nc{\s}{{\mathfrak s}}
\nc\w{\text{w}}

\renewcommand{\mod}{{\text{-mod}}}

\nc\mathi\iota
\nc\Spec{\on{Spec}}
\nc\Mod{\on{Mod}}
\nc{\tw}{\widetilde{\mathfrak t}}
\nc{\pw}{\widetilde{\mathfrak p}}
\nc{\qw}{\widetilde{\mathfrak q}}
\nc{\jw}{\widetilde j}

\nc{\grb}{\overline{\Gr}}
\nc{\I}{\mathcal I}

\nc{\lambdach}{{\check\lambda}}
\nc{\Lambdach}{{\check\Lambda}{}}
\nc{\much}{{\check\mu}}
\nc{\omegach}{{\check\omega}}
\nc{\nuch}{{\check\nu}}
\nc{\etach}{{\check\eta}}
\nc{\alphach}{{\check\alpha}}
\nc{\betach}{{\check\beta}}
\nc{\rhoch}{{\check\rho}}
\nc{\ch}{{\check h}}

\nc{\Hb}{\overline{\H}}

\nc{\BA}{{\mathbb{A}}}
\nc{\BC}{{\mathbb{C}}}
\nc{\BG}{{\mathbb{G}}}
\nc{\BM}{{\mathbb{M}}}
\nc{\BD}{{\mathbb{D}}}
\nc{\BN}{{\mathbb{N}}}
\nc{\BP}{{\mathbb{P}}}
\nc{\BR}{{\mathbb{R}}}
\nc{\BZ}{{\mathbb{Z}}}
\nc{\BS}{{\mathbb{S}}}

\nc{\CA}{{\mathcal{A}}}
\nc{\CB}{{\mathcal{B}}}

\nc{\CE}{{\mathcal{E}}}
\nc{\CF}{{\mathcal{F}}}

\nc{\CL}{{\mathcal{L}}}
\nc{\CC}{{\mathcal{C}}}
\nc{\CM}{{\mathcal{M}}}
\nc{\CN}{{\mathcal{N}}}
\nc{\CK}{{\mathcal{K}}}
\nc{\CO}{{\mathcal{O}}}
\nc{\CP}{{\mathcal{P}}}
\nc{\CQ}{{\mathcal{Q}}}
\nc{\CR}{{\mathcal{R}}}
\nc{\CS}{{\mathcal{S}}}
\nc{\CT}{{\mathcal{T}}}
\nc{\CU}{{\mathcal{U}}}
\nc{\CV}{{\mathcal{V}}}
\nc{\CW}{{\mathcal{W}}}
\nc{\CZ}{{\mathcal{Z}}}
\nc{\CX}{{\mathcal{X}}}
\nc{\CI}{{\mathcal{I}}}

\nc{\csM}{{\check{\mathcal A}}{}}
\nc{\oM}{{\overset{\circ}{\mathcal M}}{}}
\nc{\obM}{{\overset{\circ}{\mathbf M}}{}}
\nc{\oCA}{{\overset{\circ}{\mathcal A}}{}}
\nc{\obA}{{\overset{\circ}{\mathbf A}}{}}
\nc{\ooM}{{\overset{\circ}{M}}{}}
\nc{\osM}{{\overset{\circ}{\mathsf M}}{}}
\nc{\vM}{{\overset{\bullet}{\mathcal M}}{}}
\nc{\nM}{{\underset{\bullet}{\mathcal M}}{}}
\nc{\oD}{{\overset{\circ}{\mathcal D}}{}}
\nc{\obD}{{\overset{\circ}{\mathbf D}}{}}
\nc{\oA}{{\overset{\circ}{\mathbb A}}{}}
\nc{\op}{{\overset{\bullet}{\mathbf p}}{}}
\nc{\cp}{{\overset{\circ}{\mathbf p}}{}}
\nc{\oU}{{\overset{\bullet}{\mathcal U}}{}}
\nc{\oZ}{{\overset{\circ}{\mathcal Z}}{}}
\nc{\ofZ}{{\overset{\circ}{\mathfrak Z}}{}}
\nc{\oF}{{\overset{\circ}{\fF}}}

\nc{\fa}{{\mathfrak{a}}}
\nc{\fb}{{\mathfrak{b}}}
\nc{\fd}{{\mathfrak{d}}}
\nc{\fg}{{\mathfrak{g}}}
\nc{\fgl}{{\mathfrak{gl}}}
\nc{\fh}{{\mathfrak{h}}}
\nc{\fj}{{\mathfrak{j}}}
\nc{\fl}{{\mathfrak{l}}}
\nc{\fm}{{\mathfrak{m}}}
\nc{\fn}{{\mathfrak{n}}}
\nc{\fu}{{\mathfrak{u}}}
\nc{\fp}{{\mathfrak{p}}}
\nc{\fr}{{\mathfrak{r}}}
\nc{\fs}{{\mathfrak{s}}}
\nc{\fsl}{{\mathfrak{sl}}}
\nc{\hsl}{{\widehat{\mathfrak{sl}}}}
\nc{\hgl}{{\widehat{\mathfrak{gl}}}}
\nc{\hg}{{\widehat{\mathfrak{g}}}}
\nc{\chg}{{\widehat{\mathfrak{g}}}{}^\vee}
\nc{\hn}{{\widehat{\mathfrak{n}}}}
\nc{\chn}{{\widehat{\mathfrak{n}}}{}^\vee}

\nc{\fA}{{\mathfrak{A}}}
\nc{\fB}{{\mathfrak{B}}}
\nc{\fD}{{\mathfrak{D}}}
\nc{\fE}{{\mathfrak{E}}}
\nc{\fF}{{\mathfrak{F}}}
\nc{\fG}{{\mathfrak{G}}}
\nc{\fK}{{\mathfrak{K}}}
\nc{\fL}{{\mathfrak{L}}}
\nc{\fM}{{\mathfrak{M}}}
\nc{\fN}{{\mathfrak{N}}}
\nc{\fP}{{\mathfrak{P}}}
\nc{\fU}{{\mathfrak{U}}}
\nc{\fV}{{\mathfrak{V}}}
\nc{\fZ}{{\mathfrak{Z}}}

\nc{\bb}{{\mathbf{b}}}
\nc{\bc}{{\mathbf{c}}}
\nc{\bd}{{\mathbf{d}}}
\nc{\be}{{\mathbf{e}}}
\nc{\bj}{{\mathbf{j}}}
\nc{\bn}{{\mathbf{n}}}
\nc{\bp}{{\mathbf{p}}}
\nc{\bq}{{\mathbf{q}}}
\nc{\bu}{{\mathbf{u}}}
\nc{\bv}{{\mathbf{v}}}
\nc{\bx}{{\mathbf{x}}}
\nc{\bs}{{\mathbf{s}}}
\nc{\by}{{\mathbf{y}}}
\nc{\bw}{{\mathbf{w}}}
\nc{\bA}{{\mathbf{A}}}
\nc{\bK}{{\mathbf{K}}}
\nc{\bB}{{\mathbf{B}}}
\nc{\bC}{{\mathbf{C}}}
\nc{\bG}{{\mathbf{G}}}
\nc{\bD}{{\mathbf{D}}}
\nc{\bH}{{\mathbf{H}}}
\nc{\bM}{{\mathbf{M}}}
\nc{\bN}{{\mathbf{N}}}
\nc{\bV}{{\mathbf{V}}}
\nc{\bW}{{\mathbf{W}}}
\nc{\bX}{{\mathbf{X}}}
\nc{\bZ}{{\mathbf{Z}}}
\nc{\bS}{{\mathbf{S}}}

\nc{\sA}{{\mathsf{A}}}
\nc{\sB}{{\mathsf{B}}}
\nc{\sC}{{\mathsf{C}}}
\nc{\sD}{{\mathsf{D}}}
\nc{\sF}{{\mathsf{F}}}
\nc{\sK}{{\mathsf{K}}}
\nc{\sM}{{\mathsf{M}}}
\nc{\sO}{{\mathsf{O}}}
\nc{\sW}{{\mathsf{W}}}
\nc{\sQ}{{\mathsf{Q}}}
\nc{\sP}{{\mathsf{P}}}
\nc{\sZ}{{\mathsf{Z}}}
\nc{\sfp}{{\mathsf{p}}}
\nc{\sr}{{\mathsf{r}}}
\nc{\sk}{{\mathsf{k}}}
\nc{\sg}{{\mathsf{g}}}
\nc{\sff}{{\mathsf{f}}}
\nc{\sfb}{{\mathsf{b}}}
\nc{\sfc}{{\mathsf{c}}}
\nc{\sd}{{\mathsf{d}}}

\nc{\BK}{{\bar{K}}}

\nc{\tA}{{\widetilde{\mathbf{A}}}}
\nc{\tB}{{\widetilde{\mathcal{B}}}}
\nc{\tg}{{\widetilde{\mathfrak{g}}}}
\nc{\tG}{{\widetilde{G}}}
\nc{\TM}{{\widetilde{\mathbb{M}}}{}}
\nc{\tO}{{\widetilde{\mathsf{O}}}{}}
\nc{\tU}{{\widetilde{\mathfrak{U}}}{}}
\nc{\TZ}{{\tilde{Z}}}
\nc{\tx}{{\tilde{x}}}
\nc{\tbv}{{\tilde{\bv}}}
\nc{\tfP}{{\widetilde{\mathfrak{P}}}{}}
\nc{\tz}{{\tilde{\zeta}}}
\nc{\tmu}{{\tilde{\mu}}}

\nc{\urho}{\underline{\rho}}
\nc{\uB}{\underline{B}}
\nc{\uC}{{\underline{\mathbb{C}}}}
\nc{\ui}{\underline{i}}
\nc{\uj}{\underline{j}}
\nc{\ofP}{{\overline{\mathfrak{P}}}}
\nc{\oB}{{\overline{\mathcal{B}}}}
\nc{\og}{{\overline{\mathfrak{g}}}}
\nc{\oI}{{\overline{I}}}

\nc{\eps}{\varepsilon}
\nc{\hrho}{{\hat{\rho}}}

\nc{\one}{{\mathbf{1}}}
\nc{\two}{{\mathbf{t}}}

\nc{\Rep}{{\mathop{\operatorname{\rm Rep}}}}
\nc{\Tot}{{\mathop{\operatorname{\rm Tot}}}}
\nc{\Ker}{{\mathop{\operatorname{\rm Ker}}}}
\nc{\Hilb}{{\mathop{\operatorname{\rm Hilb}}}}
\nc{\End}{{\mathop{\operatorname{\rm End}}}}
\nc{\Ext}{{\mathop{\operatorname{\rm Ext}}}}
\nc{\CHom}{{\mathop{\operatorname{{\mathcal{H}}\it om}}}}
\nc{\GL}{{\mathop{\operatorname{\rm GL}}}}
\nc{\gr}{{\mathop{\operatorname{\rm gr}}}}
\nc{\Id}{{\mathop{\operatorname{\rm Id}}}}
\nc{\de}{{\mathop{\operatorname{\rm def}}}}
\nc{\length}{{\mathop{\operatorname{\rm length}}}}
\nc{\supp}{{\mathop{\operatorname{\rm supp}}}}

\nc{\Cliff}{{\mathsf{Cliff}}}
\nc{\Fl}{\on{Fl}}
\nc{\Fib}{{\mathsf{Fib}}}
\nc{\Coh}{{\mathsf{Coh}}}
\nc{\FCoh}{{\mathsf{FCoh}}}

\nc{\reg}{{\text{\rm reg}}}

\nc{\cplus}{{\mathbf{C}_+}}
\nc{\cminus}{{\mathbf{C}_-}}
\nc{\cthree}{{\mathbf{C}_*}}
\nc{\Qbar}{{\bar{Q}}}
\nc\Eis{\on{Eis}}
\nc\Eisb{\ol\Eis{}}
\nc\wh{\widehat}
\nc{\Def}{\on{Def_{\check{\fb}}(E)}}

\nc{\fq}{\mathfrak q}

\nc{\fqb}{\ol{\fq}{}}
\nc{\fpb}{\ol{\fp}{}}

\nc{\hattimes}{\wh\otimes}

\nc{\bh}{{\bar{h}}}
\nc{\bOmega}{{\overline{\Omega(\check \fn)}}}

\nc{\seq}[1]{\stackrel{#1}{\sim}}

\nc{\cT}{{\check{T}}}
\nc{\cG}{{\check{G}}}
\nc{\cM}{{\check{M}}}
\nc{\cB}{{\check{B}}}

\nc{\ct}{{\check{\mathfrak t}}}
\nc{\cg}{{\check{\fg}}}
\nc{\cb}{{\check{\fb}}}
\nc{\cn}{{\check{\fn}}}

\nc{\cLambda}{{\check\Lambda}}

\nc{\cla}{{\check\lambda}}
\nc{\cmu}{{\check\mu}}
\nc{\cnu}{{\check\nu}}
\nc{\ceta}{{\check\eta}}

\nc{\DefbE}{{\on{Def}_{\cB}(E_\cT)}}

\nc{\imathb}{{\ol{\imath}}}

\begin{document}

\title{Deformations of local systems and Eisenstein series}

\author{Alexander Braverman and Dennis Gaitsgory}

\dedicatory{To our teacher J.~Bernstein}

\address
{\newline
A.B.: Dept. of Math., Brown Univ., Providence, RI 02912, USA;\newline
D.G.: Dept. of Math., Harvard Univ., Cambridge, MA 02138, USA}

\email{braval@math.brown.edu; gaitsgde@math.harvard.edu}

\date{May 2006, revised January 2007}

\maketitle

\section*{Introduction}

\ssec{}

The goal of this paper is to realize a suggestion made by V.~Drinfeld. To explain
it let us recall the general framework of the geometric Langlands correspondence.

\medskip

Let $G$ be a reductive group and $X$ a (smooth and complete) curve.
Let $\Bun_G$ denote the moduli stack of $G$-torsors on $X$; let
$D^b(\Bun_G)$ be the appropriately defined derived category
of constructible sheaves. Let $\cG$ denote the Langlands dual 
group of $G$. 

The nature of $\cG$ depends on the sheaf-theoretic context we work in. The 
following are the main options. If we live over a ground field $k$ of 
characteristic $0$, we can work with holonomic D-modules on schemes over $k$, 
and $\cG$ will be an algebraic group over $k$.
If the ground field is $\BC$, we can work with sheaves 
of $k'$-vector spaces (here $k'$ is another
field of characteristic $0$) in the analytic topology; in this case $\cG$ is
an algebraic group over $k'$. Finally, over any ground field, we can 
work with $\ell$-adic sheaves (where $\ell$ is different from the 
ground field); in this case $\cG$ is an algebraic group over $\Ql$.

For the duration of this introduction we can work in either of the
above sheaf-theoretic contexts.

\medskip

Let $E_\cG$ be a $\cG$-local system on $X$, thought of as
a tensor functor $V\mapsto V_{E_\cG}$ from the category $\Rep(\cG)$ 
of finite-dimensional $\cG$-representations to that of local 
systems (=lisse sheaves) on $X$.

In this case one introduces the notion of {\it Hecke eigensheaf},
which is an object $\CS(E_\cG)\in D^b(\Bun_G)$, satisfying
\begin{equation} \label{Hecke eigen}
\on{H}^V(\CS(E_\cG))\simeq \CS(E_\cG)\boxtimes V_{E_\cG},
\end{equation}
where $\on{H}^V:D^b(\Bun_G)\to D^b(\Bun_G\times X)$
are the Hecke functors, defined for each $V\in \Rep(\cG)$;
the isomorphisms \eqref{Hecke eigen} are required to satisfy certain
compatibility conditions, that we will not list here.

\medskip

A basic (but in general unconfirmed, and perhaps even imprecise) expectation
is that for every $E_\cG$ there corresponds a non-zero Hecke eigensheaf
$\CS(E_\cG)$. This is a weak form of the geometric Langlands conjecture.
A stronger form of the conjecture, which only makes sense in the context of 
D-modules, says that the assignment $E_\cG\mapsto \CS(E_\cG)$ 
should work in families. In other words, if $E_{\cG,Y}$ is a 
$Y$-family of $\cG$-local systems, where $Y$ is a scheme 
over $k$, then to it there should correspond a $Y$-family
$\CS(E_{\cG,Y})$. The necessity to use D-modules here, as opposed
to any other sheaf-theoretic context, is that it is only in this case that
we have a reasonable notion of $Y$-families of objects of
$D^b(\CX)$ on a scheme (or stack) $\CX$. 

The strongest (and boldest)
form of the geometric Langlands conjecture says that the above
assignment should give rise to an equivalence between the category
$D^b(\Bun_G)$ and the appropriately defined derived category of
quasi-coherent sheaves on the stack $\on{LocSys}_{\cG}$,
classifying $\cG$-local systems on $X$.

Let us, however, consider the following intermediate case. Let $E_\cG$
be a fixed local system, and let $E_{\cG,Y}$ be its {\it formal} deformation.
We note that the notion of a $Y$-family of sheaves makes sense in
any of the above sheaf-theoretic contexts, when $Y$ is a {\it formal} 
scheme. Suppose we have found $\CS(E_\cG)$ which is a Hecke
eigensheaf with respect to $E_\cG$. 

Can we extend $\CS(E_\cG)$
to a $Y$-family $\CS(E_{\cG,Y})$ of eigensheaves? 
This is the question that Drinfeld asked on several occasions.

\ssec{}

In addition to posing this question, Drinfeld emphasized the following
general principle that should lead to a solution.

Let $\CM$ be an object of an abelian category $\CC$. (We will take
$\CM$ to be $\CS(E_\cG)$ as an object of the category of perverse
sheaves on $\Bun_G$, when $\CS(E_\cG)$ is known to exist and
be perverse.) 

Let $Y$ be a formal scheme, or more generally a formal DG-scheme,
which means by definition that $Y=\Spec(A^\bullet)$, where $A^\bullet$ is 
a (super)-commutative formal DG-algebra. Assume that $A^\bullet$ is 
isomorphic, or rather quasi-isomorphic, to the standard (=Chevalley) 
complex of a DG Lie algebra $L^\bullet$. 

Recall that the standard complex equals
the symmetric algebra on the topological graded vector space 
$(L^\bullet)^*[-1]$ with the differential induced by the differential
on $L^\bullet$ and the Lie bracket. (In our main example $L^\bullet$
will only have cohomology in degree $1$, implying that $A^\bullet$
is acyclic off cohomological degree $0$, so that $Y$ is an
honest (and not a DG) formal scheme.)

\medskip

The following principle is known as the Quillen (or Koszul) duality:

\medskip

\hskip0.5cm{\it A data of deformation of $\CM$ over the base $Y$
is equivalent, up to quasi-isomorphism, to a data of action of
$L^\bullet$ on $\CM$.}

\medskip

Returning to our problem, let us take $Y$ to be $\on{Def}(E_\cG)$--
the base of the universal deformation of $E_\cG$ as a $\cG$-local system. 
In this case, $Y$ is indeed quasi-isomorphic to the standard 
complex of a DG algebra canonically attached to $E_\cG$.
Namely, let $\cg_{X,E_\cG}$ be the sheaf of Lie algebras on $X$,
associated with the adjoint representation of $\cG$. A general
result of deformation theory says:

\medskip

\hskip0.5cm{\it The DG formal scheme $\on{Def}(E_\cG)$ corresponds 
to the Lie 
algebra $R\Gamma(X,\cg_{X,E_\cG})$}.

\medskip

Our understanding is that results of this type were first discussed in a letter
by V.~Drinfeld to V. Schechtman, and worked out in a series of
papers by V.~Hinich and V.~Schechtman \cite{HS} and V.~Hinich \cite{Hi,Hi1}. 
In any case, the above assertion is a theorem
in the context of local systems of D-modules and sheaves in
the classical topology over $\BC$. We are not sure of its status
for $\ell$-adic sheaves, and for that reason we will
avoid evoking it in this context.

\medskip

Summarizing, we obtain that the existence of the the family 
$\CS(E_{\cG,\on{Def}(E_\cG)})$ of Hecke eigen-sheaves
parametrized by $\on{Def}(E_\cG)$ is equivalent to the
existence of the action of the DG Lie algebra $R\Gamma(X,\cg_{X,E_\cG})$
on $\CS(E_\cG)$.

\medskip

Let us note now that the existence of such an action is heuristically
very natural: we expect the assignment $E_\cG\to \CS(E_\cG)$ to
be functorial; in particular we want automorphisms of $E_\cG$
to act on $\CS(E_\cG)$. Therefore, $R\Gamma(X,\cg_{X,E_\cG})$, which
can be thought of as the Lie algebra of derived endomorphisms of $E_\cG$,
should act by derived endomorphisms of $\CS(E_\cG)$, which is
what we are looking for.
 
\ssec{}

Unfortunately, even the assignment $E_\cG\to \CS(E_\cG)$ has been
constructed only in few cases.

\medskip

One such case is when $G=GL_n$ and $E_\cG=E_n$ is an $n$-dimensional
irreducible local system. In this case the existence of $\CS(E_n)$ is
known, and the action of
$R\Gamma(X,\cg_{X,E_\cG})\simeq R\on{End}(E_n)$ on $\CS(E_\cG)$,
has been fully investigated by S.~Lysenko in \cite{Lys}. In particular,
in {\it loc.cit.} it was shown that the family $\CS(E_{n,\on{Def}(E_n)})$
has the properties expected from the most general form of the
geometric Langlands conjecture, mentioned above.

\medskip

The case that we will study in this paper is, in some sense, the opposite
one. We will take $\cG$ to be arbitrary, but $E_\cG$ will be assumed
"maximally reducible", i.e., $E_\cG$ is induced from a local system
$E_\cT$ with respect to the Cartan group $\cT\subset \cG$. 
In this case the corresponding Hecke eigensheaf was constructed
in \cite{BG}, under the name "geometric Eisenstein series".
In this paper we denote it by $\Eisb(E_\cT)$. We will review the
construction of $\Eisb(E_\cT)$ later on. 

Remarkably, the action of $R\Gamma(X,\cg_{X,E_\cG})$ on
$\Eisb(E_\cT)$ has been essentially constructed in \cite{FFKM}, for
independent reasons. Thus, we do have the object
$\CS(E_{\cG,\on{Def}(E_\cG)})$, and our current goal is
to describe it more explicitly.

\medskip

At the moment, however, an explicit description of $\CS(E_{\cG,\on{Def}(E_\cG)})$
is beyond what we know how to do. We will be able
to address a more modest question, though:

\medskip

Namely, let $\DefbE$ be the base of the universal deformation of $E_\cT$,
thought of as a $\cB$-local system (here $\cB$ is the Borel subgroup of
$\cG$), such that the induced $\cT$-local system under the canonical
projection $\cB\twoheadrightarrow \cT$ is fixed to be $E_\cT$.

For example, in the case $G=GL_2$, in which case $\cG$ is also isomorphic
to $GL_2$, a $\cT$-local system can be thought of as a pair of 
$1$-dimensional local systems $(E_1,E_2)$, and we will be looking for 
$2$-dimensional local systems of the form
$$0\to E_1\to E\to E_2\to 0.$$
The (formal) scheme of such local systems is isomorphic to (the completion
at the origin of) the vector space $\on{Ext}^1(E_2,E_1)$, if we ignore
the DG complications.

The DG formal scheme $\DefbE$ maps naturally to $\on{Def}(E_\cG)$,
so we can restrict and obtain a $\DefbE$-family of Hecke eigensheaves
$\CS(E_{\cG,\DefbE})$. In fact, $\DefbE$ corresponds to a DG Lie
subalgebra in $R\Gamma(X,\cg_{X,E_\cG})$, namely, $R\Gamma(X,\cn_{X,E_\cT})$,
where $\cn$ is the nilpotent radical of $\cb$, and $\cn_{X,E_\cT}$ is the
corresponding local system of Lie algebras on $X$, twisted by 
$E_\cT$ using the adjoint action of $\cT$ on $\cn$.

\ssec{}

A concrete question posed by Drinfeld was that of explicit description of
$\CS(E_{\cG,\DefbE})$. However, more recently (in the fall of 2003) he himself 
suggested an answer:

\medskip

Along with the geometric Eisenstein series $\Eisb(E_\cT)$ there
exists a more naive object, that we can call "classical" Eisenstein series;
in this paper we denote it by $\Eis_!(E_\cT)$. When we work over the finite
ground field and $\ell$-adic sheaves, $\Eis_!(E_\cT)$ goes over under
the {\it faisceaux-fonctions} correspondence to the usual Eisenstein
series as defined in the theory of automorphic forms.

\medskip

Drinfeld's conjecture was that the family $\CS(E_{\cG,\DefbE})$ is
nothing but (a certain completion of) $\Eis_!(E_\cT)$. 
This statement has an ideological significance also for the classical
(i.e., function theoretic vs. sheaf-theoretic) Langlands correspondence:

\smallskip

\hskip0.3cm{\it The classical Eisenstein series correspond not to homomorphisms
$\on{Galois}\to \cG$ that factor through $\cT$, but rather to the universal
family of homomorphisms $\on{Galois}\to \cB$ with a fixed composition
$\on{Galois}\to \cB\twoheadrightarrow \cT$.}

\medskip

The present paper is devoted to the proof of Drinfeld's conjecture, under a certain 
simplifying hypothesis on $E_\cT$.

\medskip

Namely, we will assume that $E_\cT$ is {\it regular}, i.e., that it is as 
non-degenerate as possible. This means, by definition, that 
for every co-root $\check\alpha$ of $G$, which is the
same as a root of $\cG$, the induced $1$-dimensional
local system $\check\alpha(E_\cT)$ is non-trivial. For example,
in the case of $GL_2$ this means that the two local systems $E_1$
and $E_2$ are non-isomorphic.

This regularity assumption is equivalent to requiring that the DG
formal scheme $\DefbE$  be an "honest" scheme. In the $GL_2$
example this manifests itself in the absence of $\on{Hom}(E_2,E_1)$
and $\on{Ext}^2(E_2,E_1)$.

Moreover, we show that in this case both versions
of  Eisenstein series, namely, $\Eisb(E_\cT)$ and $\Eis_!(E_\cT)$ are 
perverse sheaves. This fact and the simplified nature of $\DefbE$
makes life significantly easier, since we can avoid a lot of complications 
of homotopy-theoretic nature. However, we are sure that, once
properly formulated, Drinfeld's conjecture that $\CS(E_{\cG,\DefbE})\simeq
\Eis_!(E_\cT)$, is true for any $E_\cT$.

\medskip

Proving the above conjecture amounts to the following:

\begin{itemize}

\item Exhibiting the action of the commutative
algebra $\CO_{\DefbE}$ of functions on $\DefbE$ on 
$\Eis_!(E_\cT)$.

\item Establishing an isomorphism
$\BC\overset{L}{\underset{\CO_{\DefbE}}\otimes}\Eis_!(E_\cT)\simeq
\Eisb(E_\cT)$, compatible with the $R\Gamma(X,\cn_{X,E_\cT})$ actions.

\item Verifying the Hecke property
$\on{H}^V(\Eis_!(E_\cT))\simeq \Eis_!(\cT)\underset{\CO_{\DefbE}}\boxtimes
V_{E_{\cG,\DefbE}}$,
where $V_{E_{\cG,\DefbE}}$ is the canonical family of $\cG$-local systems
over $\DefbE$.

\end{itemize}

As is to be expected, the verification of these properties is a nice simple
exercise when $G=GL_2$, that we will perform in Sections \ref{case of GL(2)}
and \ref{Hecke prop for GL(2)}, but not altogether trivial for other groups.

In the next section we will review the definitions of $\Eisb(E_\cT)$ and $\Eis_!(E_\cT)$,
and state each of the above properties precisely as a theorem.

\bigskip

\noindent{\bf Acknowledgments.} This introduction makes it clear
how much we owe V.~Drinfeld for the existence of this paper.
We would like to thank A.~Beilinson for developing and 
generously explaining to us the theory of chiral homology, 
which is crucial for a manageable description of the deformation 
base $\DefbE$. We would also like to thank M.~Finkelberg for 
numerous illuminating discussions about Drinfeld's compactifications
and, particularly, on the paper \cite{FFKM}.

\section{Background and overview}

\ssec{Notation and conventions}

The notation in this paper by and large follows that of \cite{BG}. We refer
the reader to {\it loc. cit.} for conventions regarding stacks,
derived categories, etc. 

We can work in any of the sheaf-theoretic contexts
mentioned in the introduction {\it excluding} the following one: 
$\ell$-adic sheaves on schemes over a ground field of positive 
characteristic {\it other than} ${\mathbb F}_q$. \footnote{The latter exclusion is due
to the fact that we will be using Kashiwara's conjecture, proved
by Drinfeld in \cite{Dr}. The only one place in this paper that
uses this result is \thmref{Eisenstein is perverse}. We have no
doubt, however, that \thmref{Eisenstein is perverse} remains
valid in the context of $\ell$-adic sheaves over fields of
positive characteristic, see \secref{eis is perv} for an additional
remark.} To simplify the notation, from now on we shall assume
that the field of coefficients of our sheaves is $\BC$, and 
for a scheme or a stack
$Y$ we shall denote by $\BC_Y$ {\it the constant sheaf} on $Y$.

\medskip

Let $G$ be a reductive group; as in \cite{BG} we will make a simplifying
assumption that the derived group of $G$ is simply connected. Let $B$
be a (fixed) Borel subgroup of $G$ and $T$ its Cartan quotient. We let
$\Lambda$ denote the lattice of characters of $T$ and $\cLambda$
the dual lattice. By $\Lambda^+$ (resp., $\Lambda^{pos}$) we will denote the
semi-group of dominant weights (resp, the semi-group generated by
positive linear combinations of simple roots); by $\cLambda^+$ and $\cLambda^{pos}$ 
we will denote the corresponding objects for the Langlands dual group.

We let $I$ denote the set of vertices of the Dynkin diagram of $G$. For
$i\in I$ we will denote by $\alpha_i\in \Lambda^{pos}$ and 
$\check\alpha_i\in \cLambda^{pos}$ the corresponding simple root and
co-root, respectively.

\ssec{Drinfeld's compactifications}

Let $\Bun_G$ (resp., $\Bun_B$, $\Bun_T$) be the moduli stacks of
$G$ (resp., $B$, $T$) torsors on $X$. We have the natural maps
$$\Bun_G\overset{\fp}\leftarrow \Bun_B\overset{\fq}\to \Bun_T.$$
The stacks $\Bun_B$ and $\Bun_T$ are both unions of connected components,
numbered by elements of $\cLambda$. We will sometimes use the
notation $\fp^\cmu$ and $\fq^\cmu$ for the restrictions of $\fp$ and
$\fq$ to the connected component $\Bun_B^\cmu$.

\medskip

Let us now recall the definition of the stack $\BunBb$. By definition,
it classifies triples $(\CP_G,\CP_T,\{\kappa^\lambda\})$, where $\CP_G$
is a $G$-torsor, $\CP_T$ is a $T$-torsor, and $\kappa^\lambda$ is an
injective map of coherent sheaves
$$\kappa^\lambda:\CL^\lambda_{\CP_T}\to V^\lambda_{\CP_G},$$
defined for each $\lambda\in \Lambda^+$, where $\CL^\lambda_{\CP_T}$
is the line bundle associated to $\CP_T$ and the character $T\overset{\lambda}\to \BG_m$,
$V^\lambda$ is the corresponding highest weight representation of $G$,
and $V^\lambda_{\CP_G}$ is the associated vector bundle. The collection
$\{\kappa^\lambda\}$ is required to satisfy the Pl\"ucker relations, see \cite{BG},
Sect. 1.2.1.

Let $\fpb$ (resp., $\fqb$) denote the natural morphism from $\BunBb$ to
$\Bun_G$ (resp., $\Bun_T$) that sends a triple $(\CP_G,\CP_T,\{\kappa^\lambda\})$
to $\CP_G$ (resp., $\CP_T$). It is a basic fact, established, e.g., in \cite{BG}, 
Proposition 1.2.2, that the map $\fpb$ is representable and proper.

The stack $\BunBb$ also splits into connected
components, denoted $\BunBb^\cmu$, $\cmu\in \cLambda$, according to the
degree of $\CP_T$. We shall denote by $\fpb^\cmu$ (resp., $\fqb^\cmu$) the
restriction of $\fpb$ (resp., $\fqb$) to the corresponding connected component.

\medskip

Consider the open substack of $\BunBb$ corresponding to the condition that
all the maps $\kappa^\lambda$ are (injective) bundle maps. This substack
identifies naturally with $\Bun_B$. We will denote by $\jmath$ the corresponding
open embedding, so that
$$\fp=\fpb\circ \jmath \text{ and } \fq=\fqb\circ \jmath.$$

Thus, $\BunBb$ can be regarded as a relative compactification of $\Bun_B$
over $\Bun_G$. 
It established in \cite{FGV}, Proposition 3.3.1, that the morphism 
$\jmath$ is affine. In fact, in \cite{BFG}, Theorem 11.6, a stronger assertion 
is proved: namely, that $\Bun_B$ is a complement to an effective 
Cartier divisor on $\BunBb$.

\ssec{}

The stack $\BunBb$ admits a natural stratification related to zeroes
of the maps $\kappa^\lambda$:

Let $\cla$ be an element of $\cLambda^{pos}$ equal to 
$\underset{i\in I}\Sigma\, n_i\cdot \check\alpha_i$,
let $|\cla|$ be the integer equal to $\Sigma\, n_i$, and let
$X^\cla$ denote the corresponding partially symmetrized power
of $X$, i.e., $$X^\cla=\underset{i\in I}\Pi\, X^{(n_i)}.$$
Points of $X^\cla$ can be thought of as effective $I$-coloured divisors
on $X$; each such point has the form $\Sigma\, \cla_k\cdot x_k$
with $x_k\neq x_{k'}$, $\cla_k\in \cLambda^{pos}$ and
$\Sigma\, \cla_k=\cla$.

\medskip

For every $\cla$ there corresponds a finite map
$$\imathb_\cla:X^\cla\times \BunBb^\cmu\to \BunBb^{\cmu-\cla}.$$
It is defined as follows. A triple $(\CP_G,\CP_T,\{\kappa^\lambda\})$
and $D^\cla\in X^\cla$ gets sent to $(\CP'_G,\CP'_T,\{\kappa'{}^\lambda\})$,
where $\CP'_G=\CP_G$, $\CP'_T=\CP_T(-D^\cla)$, i.e.,
$$\CL^\lambda_{\CP'_T}=\CL^\lambda_{\CP_T}(-\lambda(D^\cla)),$$
and $\kappa'{}^\cla$ is the composition
$$\CL^\lambda_{\CP'_T}\hookrightarrow \CL^\lambda_{\CP_T}
\overset{\kappa^\lambda}\longrightarrow V^\lambda_{\CP_G}.$$

Let $\imath_\cla$ denote the composition 
$$\imathb_\cla\circ \jmath:X^\cla\times \Bun_B\to \BunBb.$$

The maps $\imath_\cla$ are locally closed embeddings and their
images define a stratification of $\BunBb$. Since the map
$\jmath$ is affine and $\imathb_\cla$ finite, the map
$\imath_\cla$ is affine as well.

\ssec{Eisenstein series}

Let $E_\cT$ be a $\cT$-local system on $X$. Let $\CS(E_\cT)$ be the local
system on $\Bun_T\simeq \on{Pic}(X)\underset{\BZ}\otimes \cLambda$,
corresponding to $E_\cT$ via the geometric class field theory. (Our 
normalization of $\CS(E_\cT)$ is so that it is a {\it sheaf} and not a perverse
sheaf.) The defining property of $\CS(E_\cT)$ is that for $\cla\in \cLambda$,
its pull-back under the corresponding map $$\Bun_T\times X\to \Bun_T$$
is isomorphic to $\CS(E_\cT)\boxtimes \cla(E_\cT)$, where $\cla(E_\cT)$
is the induced $1$-dimensional local system on $X$
under $\cT\overset{\cla}\to \BG_m$.

\medskip

For $\cmu\in \cLambda$ the geometric (compactified) Eisenstein series 
$\Eisb^\cmu(E_\cT)$ is an object of $D^b(\Bun_G)$ defined as
$$\Eisb^\cmu(E_\cT):=\fpb^\cmu_!
\left(\IC_{\BunBb^\cmu}\otimes (\fqb^\cmu){}^*(\CS(E_\cT)\right).$$

The naive (non-compactified) Eisenstein series $\Eis_!^\cmu(E_\cT)$ is 
is defined as
$$\Eis_!^\cmu(E_\cT):=\fp^\cmu_!
\left(\IC_{\Bun_B^\cmu}\otimes (\fq^\cmu){}^*(\CS(E_\cT)\right),$$
or which is the same,
$$\fpb^\cmu_!
\left(\jmath_!(\IC_{\BunBb^\cmu})\otimes (\fqb^\cmu){}^*(\CS(E_\cT)\right).$$

Note, however, that since $\Bun_B$ (unlike $\BunBb$) is smooth,
$\IC_{\Bun_B^\cmu}$ is the constant sheaf $\BC_{\Bun_B^\cmu}$,
up to a cohomological shift.

We define $\Eisb(E_\cT)$ and $\Eis_!(E_\cT)$ as $\cLambda$-graded objects
of $D^b(\Bun_G)$, equal to $\underset{\cmu\in \cLambda}\oplus\, \Eisb^\cmu(E_\cT)$
and $\underset{\cmu\in \cLambda}\oplus\, \Eis_!^\cmu(E_\cT)$, respectively.

\medskip

Let us assume now that $E_\cT$ is regular, i.e., for every root $\check\alpha$ of $\cG$,
the $1$-dimensional local system $\check\alpha(E_\cT)$ is non-trivial.
We will prove the following:
\begin{thm}  \label{Eisenstein is perverse}  
The complexes $\Eis_!^\cmu(E_\cT)$ and $\Eisb^\cmu(E_\cT)$ are perverse sheaves.
\end{thm}

This result was conjectured in \cite{BG}; the idea of the proof 
belongs to V.~Drinfeld, and is based on the validity of 
Kashiwara's conjecture proved by him earlier.

Another fact, which we will use only marginally, is that for
$E_\cT$ regular, for every open substack $U\subset \Bun_G$ 
{\it of finite type}, the restrictions of $\Eisb^\cmu(E_\cT)|_U$
are non-zero for only finitely many $\cmu$'s. In particular,
the direct sum $\Eisb(E_\cT)$ makes sense as an object
of $\on{Perv}(\Bun_G)$.

\ssec{The space of deformations}

Before we state the main result of this paper, we need to
discuss the formal scheme of deformations $\DefbE$. From
now on we will assume that $E_\cT$ is regular. 

By definition, $\DefbE$ is a functor on the category of
local Artinian $\BC$-algebras that assigns to $R$ the set of
isomorphism classes of $R$-flat $\cB$-local systems $E_{\cB,R}$, 
such that the induced $\cT$-local system $E_{\cT,R}$, by means
of $\cB\twoheadrightarrow \cT$, is identified with $E_\cT\otimes R$, 
and such that the reduction modulo the maximal ideal, 
$E_{\cB,R/\fm_R}$ is identified with the
local system, induced from $E_\cT$ by means of 
$\cT\hookrightarrow \cB$, in a compatible way.
 
Since $H^0(X,\cn_{X,E_\cT})=0$, local systems as above
have no non-trivial automorphisms; so, by passing to
the set of isomorphism classes of objects, we do not
lose information. Moreover, since $H^2(X,\cn_{X,E_\cT})=0$,
the deformation theory is unobstructed. Hence, $\DefbE$
is representable by a smooth formal scheme.

The tangent space to $\DefbE$ at the origin is canonically
isomorphic to $H^1(X,\cn_{X,E_\cT})$. Hence, there exists
a {\it non-canonical} isomorphism between
$\DefbE$ and the completion of $H^1(X,\cn_{X,E_\cT})$ at
the origin. However, such an isomorphism is indeed
very non-canonical, and in order to proceed, we need
to describe $\DefbE$ explicitly in terms of $E_\cT$.
Such a description is provided by \thmref{deformation base}.

\medskip

Namely, in \secref{deformations} for $\cla\in \cLambda^{pos}$
we introduce a perverse sheaf $\Omega(\cn_{X,E_\cT})^{-\cla}$ on $X^\cla$.
Set $R_{E_\cT}^{-\cla}=H(X^\cla,\Omega(\cn_{X,E_\cT})^{-\cla})$. 
The regularity assumption on $E_\cT$ implies that the above
cohomology is concentrated in degree $0$.

The collection $\{\Omega(\cn_{X,E_\cT})^{-\cla}\}$ has a natural multiplicative
structure with respect to the addition operation on $\cLambda^{pos}$,
making $R_{E_\cT}:=\underset{\cla}\oplus\, R_{E_\cT}^{-\cla}$ into
a $-\cLambda^{pos}$-graded commutative algebra, with the $0$-graded
component isomorphic to $\BC$. The completion $\wh{R}_{E_\cT}$
of $R_{E_\cT}$ with respect to the augmentation ideal is isomorphic to
$\underset{\cla}\Pi\, R_{E_\cT}^{-\cla}$.

We will prove (see \thmref{proof of deformation base}) that $\wh{R}_{E_\cT}$
is canonically isomorphic to $\wh\CO_{\DefbE}$--the (topological) algebra
of functions on the formal scheme $\DefbE$. The $\cLambda$-grading
on $\wh{R}_{E_\cT}$ corresponds to the $\cT$-action on $\DefbE$,
which comes from the adjoint action of $\cT$ on $\cB$. We will
denote by $\CO_{\DefbE}$ the algebra equal to the sum of 
homogeneous components of $\wh\CO_{\DefbE}$; 
by the above, it is isomorphic to $R_{E_\cT}$.

\medskip

Let us comment on how one could guess the above description
of $\wh\CO_{\DefbE}$. As was mentioned above, $\wh\CO_{\DefbE}$ is
quasi-isomorphic to the standard complex 
$\bC^\bullet\left(R\Gamma(X,\cn_{X,E_\cT})\right)$
of the DG Lie algebra $R\Gamma(X,\cn_{X,E_\cT})$. 

The machinery of chiral algebras, developed in \cite{CHA}, implies that
$\bC^\bullet\left(R\Gamma(X,\cn_{X,E_\cT})\right)$ is quasi-isomorphic to
the chiral homology of
the (super)-commutative chiral algebra on $X$ equal to the standard
complex $\bC^\bullet(\cn_{X,E_\cT})$ of the sheaf of Lie algebras 
$\cn_{X,E_\cT}$. By definition, the above chiral homology is computed
as the homology of a sheaf associated to $\bC^\bullet(\cn_{X,E_\cT})$ 
on the {\it Ran space} of $X$. 

The latter sheaf splits into direct summands
corresponding to elements $\cla$, and each direct summand is isomorphic
to the direct image image of $\Omega(\cn_{X,E_\cT})^{-\cla}$ under the
natural map from $X^\cla$ to the Ran space. 

\ssec{The main result}

We are now ready to state one of the two main results of this paper:

\begin{thm} \label{main}  Let $E_\cT$ be regular.

\smallskip

\noindent{\em (1)}
There is a natural action of the $\cLambda$-graded
algebra $\CO_{\DefbE}$ on the $\cLambda$-graded object 
$\Eis_!(E_\cT)\in \on{Perv}(E_\cT)$.

\smallskip

\noindent{\em (2)}
The higher $\on{Tor}_i^{\CO_{\DefbE}}(\BC,\Eis_!(E_\cT))$
vanish, and we have a canonical isomorphism of $\cLambda$-graded
perverse sheaves:
$$\BC \underset{\CO_{\DefbE}}\otimes \Eis_!(E_\cT)\simeq \Eisb(E_\cT).$$
\end{thm}

In \secref{structure of open} we will give an explicit and intrinsic 
description of the pro-object 
$$\wh{\Eis}_!(E_\cT):=\Eis_!(E_\cT)\underset{\CO_{\DefbE}}
\otimes \wh\CO_{\DefbE},$$ which is our candidate for
$\CS(E_{\cG,\DefbE})$.

\medskip

Let us add a few comments on the strategy of the proof of this theorem.
Taking into account \corref{graded deformation base}, to define an action as
in point (1), we need to construct morphisms
\begin{equation} \label{action downstairs}
R^{-\cla}_{E_\cT}\otimes \Eis_!^{\cmu+\cla}(E_\cT)\to \Eis_!^\cmu(E_\cT)
\end{equation}
that are associative in the natural sense. The morphisms \eqref{action downstairs}
will be obtained by applying the functor $\fpb^\cmu_!$ to some canonical
morphism of perverse sheaves upstairs. 

Namely, in \thmref{extension by zero} we will show that there exists a
canonical map
$$\imath_\cla{}_!\left(\Omega(\cn_{X})^{-\cla}\boxtimes \IC_{\Bun_B^{\cmu+\cla}}\right)
\to \jmath_!(\IC_{\Bun_B^{\cmu}}),$$
where $\Omega(\cn_{X})^{-\cla}$ is the perverse sheaf corresponding to the trivial
twist. Tensoring both sides of the above expression by $(\fqb^\cmu){}^*(\CS(E_\cT))$
we obtain a map
$$\imath_\cla{}_!\biggl(\Omega(\cn_{X,E_\cT})^{-\cla}\boxtimes
\left(\IC_{\Bun_B^{\cmu+\cla}}\otimes (\fq^{\cmu+\cla}){}^*(\CS(E_\cT))\right)\biggr)\to
\jmath_!\left(\IC_{\Bun_B^{\cmu}}\otimes (\fq^{\cmu}){}^*(\CS(E_\cT))\right),$$
which gives rise to \eqref{action downstairs}.

\ssec{Koszul complex}

To prove point (2) of \thmref{main} we proceed as follows. In \secref{abstract Koszul}
for each $\cla\in \cLambda^{pos}$, we define a certain explicit complex of perverse 
sheaves on $X^\cla$, denoted $\fU(\cn_{X,E_\cT})^{\bullet,-\cla,*}$. 

One should think of $\Omega(\cn_{X,E_\cT})^{-\cla}$
and $\fU(\cn_{X,E_\cT})^{\bullet,-\cla,*}$ as Koszul-dual objects in the same way as
the universal enveloping algebra $U(\fh)$ of a Lie algebra $\fh$ 
is a Koszul-dual object to the co-standard complex $\bC_\bullet(\fh)$.

We show (see \thmref{exactness of Koszul upstairs}) that the graded perverse sheaf
on $\BunBb^\cmu$
\begin{equation}   \label{Koszul up intr}
\on{Kosz}_\cmu(\Eis_!(E_\cT))^\bullet:=
\underset{\cla\in \cLambda^{pos}}\bigoplus\,
\imath_\cla{}_!\left(\fU(\cn_{X,E_\cT})^{\bullet,-\cla,*}\boxtimes 
\jmath_!(\IC_{\Bun_B^{\cmu+\cla}})\right)
\end{equation}
acquires a natural differential, such that the resulting complex is
quasi-isomorphic to $\IC_{\BunBb^\cmu}$. This implies the assertion  
of point (2) of the theorem as follows:

\medskip

The regularity assumption on $E_\cT$ implies that 
$H^j(X^\cla,\fU(\cn_{X,E_\cT})^{i,-\cla,*})=0$ for $j\neq 0$. So, we obtain
well-defined complexes
$$\fu^{\bullet,-\cla,*}_{E_\cT}:=H^0(X^\cla,\fU(\cn_{X,E_\cT})^{\bullet,-\cla,*}).$$
Set $\fu^{\bullet,*}_{E_\cT}:=\underset{\cla}\oplus\, \fu^{\bullet,-\cla,*}_{E_\cT}$.
In \secref{abstract Koszul} we show that for any $R_{E_\cT}$-module $\CM$ the tensor product 
$$\sK(\CM):=\fu^{\bullet,*}_{E_\cT}\otimes \CM$$
acquires a natural differential, and the resulting complex is quasi-isomorphic 
to $\BC\overset{L}{\underset{R_{E_\cT}}\otimes}\CM$.

Taking the direct image of \eqref{Koszul up intr} with respect to
$\fpb^\cmu$, and summing up over $\cmu$, we obtain
$$\BC\overset{L}{\underset{R_{E_\cT}}\otimes} \Eis_!(E_\cT)
\simeq \sK(\Eis_!(E_\cT))\simeq 
\underset{\cmu}\oplus\, \fpb^\cmu_!(\on{Kosz}_\cmu(\Eis_!(E_\cT))^\bullet)
\simeq \Eisb(E_\cT),$$
as required.

\medskip

Let us add a comment on the nature of the complexes $\fu^{\bullet,-\cla,*}_{E_\cT}$.
In \secref{Koszul section} we show that the collection $\{\fU(\cn_{X,E_\cT})^{\bullet,-\cla,*}\}$ 
possesses a natural co-multiplicative structure, thereby endowing $\fu^{\bullet,*}_{E_\cT}$
with a structure of DG graded co-associative co-algebra (with a trivial differential (!)).

Let $\fu^{\bullet,\cla}_{E_\cT}$ be the dual of $\fu^{\bullet,-\cla,*}_{E_\cT}$, and set
$\fu^\bullet_{E_\cT}:=\underset{\cla}\oplus\, \fu^{\bullet,\cla}_{E_\cT}$. We obtain
that $\fu^\bullet_{E_\cT}$ is an associative DG algebra (also, with a trivial differential).
Although we do not state this explicitly, from \secref{abstract Koszul} one can
deduce that $\fu^\bullet_{E_\cT}$ is canonically quasi-isomorphic to the 
universal enveloping algebra of $R\Gamma(X,\cn_{X,E_\cT})$.

By construction, $\sK(\Eis_!(E_\cT))$ is a DG-comodule with respect to
$\fu^{\bullet,-\cla,*}_{E_\cT}$, and, hence, a DG-module over $\fu^\bullet_{E_\cT}$.
This structure can be viewed as an action of DG Lie algebra $R\Gamma(X,\cn_{X,E_\cT})$
on $\Eisb(E_\cT)$. By definition, this action equals the one arising via the
Koszul-Quillen duality on $\BC\overset{L}{\underset{R_{E_\cT}}\otimes} \Eis_!(E_\cT)$. 

A compatibility result proved in \secref{comparison} implies that the above
$R\Gamma(X,\cn_{X,E_\cT})$-action on $\Eisb(E_\cT)$, coincides with the one
given by the construction of \cite{FFKM}. This is equivalent to the fact that
$\Eis_!(E_\cT)$ {\it is} the $\DefbE$-family corresponding to $\Eisb(E_\cT)$
with the action of $R\Gamma(X,\cn_{X,E_\cT})$ constructed in \cite{FFKM}.

\ssec{The Hecke property}

The second main result of the present paper has to do with the verification
of the Hecke property of $\Eis_!(E_\cT)$. In order to simplify the exposition,
instead of the Hecke functor $\on{H}^V:D^b(\Bun_G)\to D^b(\Bun_G\times X)$,
we will consider the local Hecke functor
$$\on{H}_x^V:D^b(\Bun_G)\to D^b(\Bun_G),$$
corresponding to a fixed point $x\in X$. As will be clear from the proof,
the case of a moving point (or multiple points, as required for the verification
of the additional compatibility condition, which we did not even state
explicitly) is analogous.

Thus, to a point $x\in X$ there corresponds a $\cG$-torsor, equal to the
restriction to the fiber at $x$ of the universal $\cG$-local system over
$\DefbE$. For $V\in \Rep(\cG)$, let $V_{E_{\cB,\wh{univ},x}}$ be the 
corresponding locally free $\wh\CO_{\DefbE}$-module. Let 
$V_{E_{\cB,univ,x}}$ be the corresponding $\cLambda$-graded
version, which is a locally free module over $\CO_{\DefbE}$.

We will prove:
\begin{thm} \label{Hecke}
There exists a canonical isomorphism
$$\on{H}_x^V(\Eis_!(E_\cT))\simeq V_{E_{\cB,univ,x}}\underset{\CO_{\DefbE}}\otimes
\Eis_!(E_\cT).$$
\end{thm}

Let us explain the main ideas involved in the proof of this theorem. First, we need
to interpret $V_{E_{\cB,univ,x}}$ in terms of the isomorphism
$\CO_{\DefbE}\simeq R_{E_\cT}$. This is also done using chiral homology.

In \secref{Hecke section} for every $\cnu\in \cLambda$ we construct a perverse 
sheaf $\Omega(\cn_{X,E_\cT},V_{E_\cT,x})^\cnu$ that lives over an
appropriate ind-version $_{\infty\cdot x}X^\cnu$ of the space of coloured
divisors (essentially, we allow the divisor to be non-effective at $x$).
We set $$R(V_x)_{E_\cT}^\cnu:=H({}_{\infty\cdot x}X^\cnu,\Omega(\cn_{X,E_\cT},V_{E_\cT,x})^\cnu),
\text{ and } 
R(V_x)_{E_\cT}=\underset{\cnu}\oplus\, R(V_x)_{E_\cT}^\cnu.$$ We show that
$R(V_x)_{E_\cT}$ is a locally free $R_{E_\cT}$-module of rank $\dim(V)$, and
that $R(V_x)_{E_\cT}$ corresponds to
$V_{E_{\cB,univ,x}}$ under the isomorphism $\CO_{\DefbE}\simeq R_{E_\cT}$.

\medskip

Secondly, we need to reinterpret the LHS in \thmref{Hecke} locally in
terms of $\BunBb^\cmu$. To do this, as in \cite{BG}, we need to replace 
$\BunBb^\cmu$ by its ind-version $_{\infty\cdot x}\BunBb^\cmu$ that allows
the maps $\kappa^\lambda$ (see the definition of $\BunBb$) to have 
poles of an arbitrary order at $x$. Tautologically, the Hecke functors
$\on{H}_x^V$ that act on $D^b(\Bun_G)$ lift to functors, denoted
$\on{H}'{}_x^V$, that act on $D^b({}_{\infty\cdot x}\BunBb^\cmu)$,
in a way compatible with the push-forward functor
$\fpb^\cmu_!:D^b({}_{\infty\cdot x}\BunBb^\cmu)\to D^b(\Bun_G)$.

Thus, the LHS in \thmref{Hecke} is given by
$$\underset{\cnu}\bigoplus\,
\fpb^\cmu_!\biggl(\on{H}'{}_x^V\left(\jmath_!(\IC_{\Bun_B^\cmu})\otimes 
(\fqb^\cmu){}^*(\CS(E_\cT))\right)\biggr).$$

As in the case of $\BunBb^\cmu$, we
have a natural map
$$_{\infty\cdot x}\imath_\cnu:{}_{\infty\cdot x}X^\cnu\times
\Bun_B^{\cmu-\cnu}\to {}_{\infty\cdot x}\BunBb^\cmu.$$ 
In \thmref{Hecke of extension by zero}, we show that there exists
a canonical map
$$_{\infty\cdot x}\imath_\cnu{}_!
\left(\Omega(\cn_{X},V_{x})^\cnu\boxtimes
\IC_{\Bun_B^{\cmu-\cnu}}\right)\to 
\on{H}'{}_x^V\left(\jmath_!(\IC_{\Bun_B^\cmu})\right).$$

Tensoring with the pull-back of $\CS(E_\cT)$ under $\fqb^\cmu:
{}_{\infty\cdot x}\BunBb^\cmu\to \Bun_T$, we obtain a map
\begin{multline*}
_{\infty\cdot x}\imath_\cnu{}_!
\biggl(\Omega(\cn_{X,E_\cT},V_{E_\cT,x})^\cnu\boxtimes
\left(\IC_{\Bun_B^{\cmu-\cnu}}\otimes (\fq^{\cmu-\cnu})^*(\CS(E_\cT))\right)\biggr)\to \\
\to \on{H}'{}_x^V\biggl(\jmath_!\left((\IC_{\Bun_B^\cmu})\otimes 
(\fq^{\cmu})^*(\CS(E_\cT))\right)\biggr),
\end{multline*}
and applying the direct image under $\fpb^\cmu_!$, we obtain a map
in one direction ($\leftarrow$) in \thmref{Hecke}.

The proof that the resulting map is an isomorphism follows by
considering filtrations defined naturally on the two sides, and
showing that the map induced on the associated graded level
is an isomorphism.

\section{The case of $GL_2$}   \label{case of GL(2)}

In this section we will prove \thmref{main}
for $G=GL_2$ by an explicit calculation.
This will be a prototype of the argument in the general case.

\ssec{The base of deformation}

Let $E_1$ and $E_2$ be two non-isomorphic $1$-dimensional local
systems on $X$. We regard the pair $(E_1,E_2)$ as a local system $E_{\cT}$ 
with respect to the Cartan subgroup $\cT \simeq (\BG_m,\BG_m)$ of the group 
$\cG$, Langlands dual of $GL_2$, which is itself isomorphic to $GL_2$.

By definition, the formal scheme $\DefbE$ associates to a local Artinian $\BC$-algebra
$R$ the category of $R$-flat local systems $E_R$ that fit
into the short exact sequence
$$0\to E_1\otimes R\to E_R\to E_2\otimes R\to 0,$$
which modulo the maximal ideal of $R$ are identified with the
direct sum $E_1\oplus E_2$.

Note that since $E_1\neq E_2$, such short exact sequences admit no 
automorphisms, so that the above category is (equivalent to) a set.

Hence, we obtain that $\DefbE$ is naturally isomorphic to the completion
at $0$ of the $\BC$-vector space $\Ext^1(E_2,E_1)$. Let us denote the
dual vector space $H^1(X,E_2\otimes E_1^{-1})$ by $\sW$.
The corresponding complete commutative algebra $\wh\CO_{\DefbE}$ is isomorphic to 
$\wh{\Sym}(\sW)$. 

We will also consider the non-completed algebra $\Sym(\sW)$, endowed with a
grading by {\it negative} integers. We regard $\BZ$ as a subgroup of 
$\BZ\oplus \BZ\simeq \cLambda$
via $d\mapsto (d,-d)$.

\ssec{}

For each $(d_1,d_2)\in \BZ\oplus \BZ$, 
let $\Eis_!^{d_1,d_2}(E_\cT)\in D^b(\Bun_G)$ and $\Eisb^{d_1,d_2}(E_\cT)\in D^b(\Bun_G)$
be the corresponding component of the non-compactified and compactified 
Eisenstein series, respectively, attached to $E_\cT$. Both $\Eis_!^{d_1,d_2}(E_\cT)$ and
$\Eisb^{d_1,d_2}(E_\cT)$ are known (\cite{Ga}) to be perverse sheaves; the
generalization of this assertion for arbitrary $G$ (due to Drinfeld) will be established
in \secref{eis is perv}.

\medskip

Adapting the notation of \thmref{main}
to the present context, we obtain the
following:

\begin{thm} \label{main GL(2)}  \hfill

\smallskip

\noindent{\em (a)}
There exists a grading preserving action
$$\Sym(\sW)\otimes \Eis_!(E_\cT)\to \Eis_!(E_\cT).$$

\smallskip

\noindent{\em (b)} For each $(d_1,d_2)\in \BZ\oplus \BZ$,
the resulting Koszul complex
\begin{multline*}
\sK_{d_1,d_2}(\Eis_!(E_\cT))^\bullet:=
...\to \Lambda^d(\sW)\otimes \Eis_!^{d_1+d,d_2-d}(E_\cT)\overset{\partial_d}\to...\to \\
\to \Lambda^2(\sW)\otimes \Eis_!^{d_1+2,d_2-2}(E_\cT)\overset{\partial_2}\to
\sW\otimes \Eis_!^{d_1+1,d_2-1}(E_\cT)\overset{\partial_1}\to \Eis_!^{d_1,d_2}(E_\cT)
\end{multline*}
is quasi-isomorphic to $\Eisb^{d_1,d_2}(E_\cT)$.

\end{thm}

\ssec{Proof of \thmref{main GL(2)}}

Recall the stack $\BunBb^{d_1,d_2}$. For each non each non-negative integer
$d$ let $\imathb^{d_1,d_2}_d$ denote the finite map
$$X^{(d)}\times \BunBb^{d_1+d,d_2-d}\to \BunBb^{d_1,d_2}.$$
Let $\imath^{d_1,d_2}_d$ denote the composition of $\imathb^d$ and the open embedding
$$\imath_0^{d_1+d,d_2-d}:=\jmath^{d_1+d,d_2-d}:
\Bun^{d_1+d,d_2-d}_B\hookrightarrow \BunBb^{d_1+d,d_2-d}.$$

It is known (see \cite{BG}, Proposition 6.1.2) that $\imath_d$ is a locally closed embedding, and that these stacks define a stratification of $\BunBb^{d_1,d_2}$.
Moreover, the map $\jmath^{d_1+d,d_2-d}$ is known to be affine. 
Hence, the map $\imath^{d_1,d_2}_d$ is also affine.

\medskip

Recall that for $G=GL_2$, the stack $\BunBb^{d_1,d_2}$ is smooth; in 
particular, its intersection cohomology sheaf is constant. Thus, we
obtain an exact complex of perverse sheaves on $\BunBb^{d_1,d_2}$:

\begin{multline}  \label{Kosz GL(2) upstairs non-twisted}
...\to(\imath^{d_1,d_2}_d)_!(\IC_{X^{(d)}\times \Bun_B^{d_1+d,d_2-d}})\to...
...\to (\imath^{d_1,d_2}_2)_!(\IC_{X^{(2)}\times \Bun_B^{d_1+2,d_2-2}})\to \\
\to (\imath^{d_1,d_2}_1)_!(\IC_{X\times \Bun_B^{d_1+1,d_2-1}})\to
(\imath^{d_1,d_2}_0)_!(\IC_{\Bun_B^{d_1,d_2}})\to \IC_{\BunBb^{d_1,d_2}}.
\end{multline}

Let $\CS(E_\cT)$ be the local system on $\Bun_T\simeq \Pic(X)\times \Pic(X)$,
corresponding to $E_\cT$. We normalize it so that the pull-back of $\CS(E_\cT)$
under 
$$X^{(d'_1)}\times X^{(d'_2)}\overset{\on{AJ}\times \on{AJ}}\to \Pic(X)\times \Pic(X),$$
where $\on{AJ}$ is the Abel-Jacobi map $D\mapsto \CO_X(D)$,
is $E_1^{(d'_1)}\boxtimes E_2^{(d'_2)}$.

Recall that by definition,
$$\Eisb^{d_1,d_2}(E_\cT)=\fpb^{d_1,d_2}_!\left(\IC_{\BunBb^{d_1,d_2}}\otimes
\fqb^{d_1,d_2}{}^*(\CS(E_\cT))\right)$$ and 
$$\Eis_!^{d_1,d_2}(E_\cT)=\fp^{d_1,d_2}_!\left(\IC_{\Bun_B^{d_1,d_2}}\otimes
\fq^{d_1,d_2}{}^*(\CS(E_\cT))\right),$$
where $\fpb^{d_1,d_2}$ ($\fqb^{d_1,d_2}$) denote the natural projection
from $\BunBb^{d_1,d_2}$ to $\Bun_G$ (resp., $\Bun_T$), and 
$\fp^{d_1,d_2}=\fpb^{d_1,d_2}\circ \jmath^{d_1,d_2}$ 
(resp., $\fq^{d_1,d_2}=\fqb^{d_1,d_2}\circ \jmath^{d_1,d_2}$) is its
restriction to $\Bun_B^{d_1,d_2}$.

\medskip

Tensoring the complex \eqref{Kosz GL(2) upstairs non-twisted} with the local system
$\fqb^{d_1,d_2}{}^*(\CS(E_\cT))$ we obtain a complex
\begin{multline}  \label{Kosz GL(2) upstairs}
\on{Kosz}_{d_1,d_2}(\Eis_!(E_\cT))^\bullet:= \\
...\to(\imath^{d_1,d_2}_d)_!\left((E_2\otimes E_1^{-1})^{(d)}[d]\boxtimes 
\left(\IC_{\Bun_B^{d_1+d,d_2-d}}\otimes\fq^{d_1+d,d_2-d}{}^*(\CS(E_\cT))\right)\right)\to...\\
...\to(\imath^{d_1,d_2}_2)_!\left((E_2\otimes E_1^{-1})^{(2)}[2]\boxtimes 
\left(\IC_{\Bun_B^{d_1+2,d_2-2}}\otimes\fq^{d_1+2,d_2-2}{}^*(\CS(E_\cT))\right)\right)\to \\
\to(\imath^{d_1,d_2}_1)_!\left((E_2\otimes E_1^{-1})[1]\boxtimes 
\left(\IC_{\Bun_B^{d_1+1,d_2-1}}\otimes\fq^{d_1+1,d_2-1}{}^*(\CS(E_\cT))\right)\right)\to \\
\to(\imath^{d_1,d_2}_0)_!
\left(\IC_{\Bun_B^{d_1,d_2}}\otimes\fq^{d_1,d_2}{}^*(\CS(E_\cT))\right),
\end{multline}
which is quasi-isomorphic to 
$\IC_{\BunBb^{d_1,d_2}}\otimes \fqb^{d_1,d_2}{}^*(\CS(E_\cT))$.

Applying $\fpb^{d_1,d_2}_!$ to $\on{Kosz}_{d_1,d_2}(\Eis_!(E_\cT))^\bullet$, we obtain a complex
\begin{multline} \label{Kosz GL(2) downstairs}
...\to H^d(X,(E_2\otimes E_1^{-1})^{(d)})\otimes \Eis_!^{d_1+d,d_2-d}(E_\cT)\to...\to \\
\to H^2(X,(E_2\otimes E_1^{-1})^{(2)})\otimes \Eis_!^{d_1+2,d_2-2}(E_\cT)\to \\
\to H^1(X,(E_2\otimes E_1^{-1}))\otimes \Eis_!^{d_1+1,d_2-1}(E_\cT)\to
\Eis_!^{d_1,d_2}(E_\cT),
\end{multline}
which is in turn quasi-isomorphic to  $\Eisb^{d_1,d_2}(E_\cT)$.

\medskip

Let us recall that $H^d(X,(E_2\otimes E_1^{-1})^{(d)})\simeq \Lambda^d
\left(H^1(X,E_2\otimes E_1^{-1})\right)=:\Lambda^d(\sW)$. Thus, we obtain
that the terms of the complex \eqref{Kosz GL(2) downstairs} coincide with those 
of the complex $\sK_{d_1,d_2}(\Eis_!(E_\cT))^\bullet$ appearing 
in \thmref{main GL(2)}(b).

In particular, for each pair $(d'_1,d'_2)$ we obtain a map
$$\partial_1:\sW\otimes \Eis_!^{d'_1,d'_2}(E_\cT) \to \Eis_!^{d'_1-1,d'_2+1}(E_\cT).$$ 
Therefore, to finish the proof of the theorem, it suffices to show that for each $d$
the differential in \eqref{Kosz GL(2) downstairs}
\begin{multline}  
H^d(X,(E_2\otimes E_1^{-1})^{(d)})\otimes \Eis_!^{d_1+d,d_2-d}(E_\cT)\to \\
H^{d-1}(X,(E_2\otimes E_1^{-1})^{(d-1)})\otimes \Eis_!^{d_1+d-1,d_2-+1}(E_\cT)
\end{multline}
equals
\begin{multline*}
\Lambda^d(\sW)\otimes \Eis_!^{d_1+d,d_2-d}(E_\cT)\hookrightarrow
\Lambda^{d-1}(\sW)\otimes \sW\otimes \Eis_!^{d_1+d,d_2-d}(E_\cT)
\overset{\on{id}\otimes \partial_1}\longrightarrow \\ 
\longrightarrow \Lambda^{d-1}(\sW) \otimes \Eis_!^{d_1+d-1,d_2+1}(E_\cT).
\end{multline*}

Note that we have a commutative diagram of stacks
$$
\CD
X^{(d-1)}\times X\times \Bun_B^{d_1+d,d_2-d} 
@>{\on{id}\times \imath_1^{d_1+d-1,d_2-d+1}}>> 
X^{(d-1)}\times \BunBb^{d_1+d-1,d_2-d+1}  \\
@V{\on{sym}_{d-1,1}}VV    @V{\imathb_{d-1}^{d_1,d_2}}VV   \\
X^{(d)}\times \Bun_B^{d_1+d,d_2-d} 
@>{\imath_d^{d_1,d_2}}>> \BunBb^{d_1,d_2},
\endCD
$$
where $\on{sym}_{d-1,1}$ is the natural map $X^{(d-1)}\times X\to X^{(d)}$.

Hence, the map
$$(\imath_d^{d_1,d_2})_!\left(\IC_{X^{(d)}}\boxtimes \IC_{\Bun_B^{d_1+d,d_2-d}}\right)\to
(\imath_{d-1}^{d_1,d_2})_!\left(\IC_{X^{(d-1)}}\boxtimes \IC_{\Bun_B^{d_1+d-1,d_2-d+1}}\right),$$
appearing in \eqref{Kosz GL(2) upstairs non-twisted} equals the composition
\begin{multline*}
(\imath_d^{d_1,d_2})_!\left(\IC_{X^{(d)}}\boxtimes \IC_{\Bun_B^{d_1+d,d_2-d}}\right)\to
(\imath_d^{d_1,d_2})_!\left(\on{sym}_{d-1,1}{}_!(\IC_{X^{(d-1)}\times X})\boxtimes
\IC_{\Bun_B^{d_1+d,d_2-d}}\right) \simeq \\
\simeq (\imathb_{d-1}^{d_1,d_2})_!
\left(\IC_{X^{(d-1)}}\boxtimes 
(\imath_{1}^{d_1+d-1,d_2-d+1})_!(\IC_{X\times \Bun_B^{d_1+d,d_2-d}})\right)\to \\
\to (\imathb_{d-1}^{d_1,d_2})_! \left(\IC_{X^{(d-1)}}\boxtimes 
(\imath_{0}^{d_1+d-1,d_2-d+1})_!(\IC_{\Bun_B^{d_1+d-1,d_2-d+1}})\right)\simeq \\
\simeq 
(\imath_{d-1}^{d_1,d_2})_!\left(\IC_{X^{(d-1)}}\boxtimes \IC_{\Bun_B^{d_1+d-1,d_2-d+1}}\right),
\end{multline*}
implying the desired equality after taking the direct image with respect to
$\fpb^{d_1,d_2}$.

\section{Deforming local systems}   \label{deformations}

The goal of this section is to describe explicitly the topological algebra
$\wh\CO_{\DefbE}$.

\ssec{}   \label{intr Upsilon}

Let $\cn_{X}$ be the constant sheaf of Lie algebras over $X$ with
fiber $\cn$. Let $\cn_{X,E_\cT}$ denote its twist by means of $E_\cT$
with respect to the adjoint action of $\cT$ on $\cn$. It is naturally
graded by elements of $\cLambda^{pos}$, where the latter
is a sub-semigroup of $\cLambda$ equal to the positive span
of simple co-roots.
We will consider the standard complex 
$\bC_\bullet(\cn_{X,E_\cT})$ as a sheaf of co-commutative
DG co-algebras on $X$, also endowed with a grading by means of
$\cLambda^{pos}$.

\medskip

Recall that to $\cla\in \cLambda^{pos}-0$ we have attached the corresponding
partially symmetrized power of $X$, denoted $X^\cla$. 
We are going to associate to each such $\cla$ a certain perverse
sheaf $\Omega(\cn_{X,E_\cT})^{-\cla}$ on $X^\cla$. (In the non-twisted case, i.e.,
for $E_\cT$ trivial, we will denote it simply by $\Omega(\cn_{X})^{-\cla}$.)
This will be based on the following general construction.

\medskip

Let $A$ be a sheaf of co-commutative DG co-algebras on $X$,
and let $n$ be a non-negative integer. By the procedure of
\cite{CHA}, Sect. 3.4, to $A$ one can associate a sheaf, denoted,
$A_{X^{(n)}}$ on $X^{(n)}$,
whose fiber at $D=\Sigma\, m_k\cdot x_k \, |\, x_{k_1}\neq x_{k_2}$
is $\underset{k}\bigotimes\, A_{x_{k}}$. 

Suppose now that $A$ is $\cLambda^{pos}$-graded, and let
$\cla$ be an element of $\cLambda^{pos}$. We have a natural
map $X^\cla\to X^{(|\cla|)}$. Then the *-pull-back of $A_{X^{(n)}}$
to $X^\cla$ admits a sub-sheaf, denoted $A_{X^\cla}$, whose fiber
at $D=\cla_k\cdot x_k \, |\, x_{k_1}\neq x_{k_2}$ is the subspace
$$\underset{k}\bigotimes\, \left(A_{x_{k}}\right)^{\cla_k}\subset
\underset{k}\bigotimes\, A_{x_{k}}.$$

We apply this to $A=\bC_\bullet(\cn_{X,E_\cT})$ and obtain a complex
of sheaves that we denote by $\Upsilon(\cn_{X,E_\cT})^{\cla}$.
Let us describe it even more explicitly:

\medskip

Recall that $X^\cla$ admits a stratification, numbered by partitions 
$$\fP^\lambda: \cla=\underset{k}\Sigma\, m_k\cdot \cla_k\, |\, m_k\in \BZ^{>0},
\cla_k\in \cLambda^{pos}-0, \cla_{k}\neq \cla_{k'}$$
The corresponding stratum in $X^\cla$ is isomorphic to
$\left(\underset{k}\Pi\, X^{(n_k)}\right){}_{disj}$, where the subscript "$disj$"
denotes the complement to the diagonal divisor in the above product.
Its dimension is $|\fP(\cla)|=\underset{k}\Sigma\, m_k$.
Let $j^{\fP(\cla)}$ denote the corresponding locally closed embedding.

Let $\on{Cous}(\cn_{X,E_\cT})^{\fP(\cla)}$ denote the complex of (both, sheaves,
and (shifted) perverse) sheaves
$$j^{\fP(\cla)}_!\left(\underset{k}\boxtimes\,
\left(\Lambda^\bullet(\cn_{X,E_\cT})^{\cla_k}\right)^{(m_k)}\right).$$
The direct sum
$$\on{Cous}(\cn_{X,E_\cT})^\cla:=
\underset{\fP(\cla)}\bigoplus\, \on{Cous}(\cn_{X,E_\cT})^{\fP(\cla)},$$
viewed as a graded perverse sheaf,
admits a natural differential and the resulting total complex is, by definition,
quasi-isomorphic to $\Upsilon(\cn_{X,E_\cT})^{\cla}$.

\begin{prop} \label{Upsilon}
$\Upsilon(\cn_{X,E_\cT})^{\cla}$ is acyclic off cohomolgical degree $0$
in the perverse t-structure.
\end{prop}

\ssec{}   \label{gr of Omega}

In order to prove \propref{Upsilon} we need to introduce some notation
that will be also useful in the sequel. Let $\cla_1,\cla_2$ be two elements
of $\cLambda^{pos}$. Note we have a natural addition map
$$X^{\cla_1}\times X^{\cla_2}\to X^{\cla_1+\cla_2}.$$
This map is finite, and it induces an exact functor
$$\Perv(X^{\cla_1})\times \Perv(X^{\cla_2})\to \Perv(X^{\cla_1+\cla_2})$$
that we will denote by $\CT_1,\CT_2\mapsto \CT_1\star \CT_2$.

\medskip

The natural increasing filtration on $\bC_\bullet(\cn_{X,E_\cT})$ defines
a filtration on $\Upsilon(\cn_{X,E_\cT})^{\cla}$. The associated graded
can be obtained by the same procedure, where instead of $\cn_{X,E_\cT}$
we use the abelian Lie algebra structure on the same sheaf. Hence,
$$\on{gr}^j\left(\Upsilon(\cn_{X,E_\cT})^{\cla}\right)=
\underset{\cla=\Sigma\, n_k\cdot \check\alpha_k,\, \check\alpha_k\in \check\Delta^+,\,
\Sigma\, n_k=j}\bigoplus\, 
\underset{k}\star\, \left(\Lambda^{(n_k)}(E_\cT^{\check\alpha_k})\right)[j],$$
where in the above formula for a local system $\sF$ we denote by $\Lambda^{(n)}(\sF)$
its external exterior power.

This description of the associated graded readily implies that $\Upsilon(\cn_{X,E_\cT})^{\cla}$
is a perverse sheaf. In addition, we obtain the following:

\begin{cor}  \label{cohomology of ups}
Assume that $E_\cT$ is regular. Then for every $\cla\in \cLambda^{pos}$
the cohomology
$$H^\bullet\left(X^\cla,\Upsilon(\cn_{X,E_\cT})^{\cla}\right)$$
is concentrated in degree $0$.
\end{cor}

\begin{proof}

It is enough to prove that each 
$H^\bullet\left(X^\cla,\on{gr}^j(\Upsilon(\cn_{X,E_\cT})^{\cla})\right)$ 
is concentrated in degree $0$. The latter follows from the fact that
$$H^\bullet\left(X^{(n)}, \Lambda^{(n)}(E_\cT^{\check\alpha})[n]\right)\simeq
\Sym^{n}\left(H^\bullet(X,E_\cT^{\check\alpha})[1]\right)\simeq
\Sym^n\left(H^1(X,E_\cT^{\check\alpha})\right).$$

\end{proof}

\ssec{}   \label{discuss omega}

Let $\Omega(\cn_{X,E_\cT})^{-\cla}$ denote the Verdier dual of 
$\Upsilon(\cn_{X,E_\cT})^{\cla}$. From \corref{cohomology of ups}
it follows that $H^\bullet\left(X^\cla,\Upsilon(\cn_{X,E_\cT})^{\cla}\right)$
is also concentrated in degree $0$. 

From the construction of $\Upsilon(\cn_{X,E_\cT})^{\cla}$ it follows
that we have natural maps
\begin{equation} \label{comult upsilon}
\Upsilon(\cn_{X,E_\cT})^{\cla_1+\cla_2}\to
\Upsilon(\cn_{X,E_\cT})^{\cla_1}\star \Upsilon(\cn_{X,E_\cT})^{\cla_2},
\end{equation}
that are co-associative and co-commutative in the natural sense. Hence,
we obtain the maps
$$\Omega(\cn_{X,E_\cT})^{-\cla_1}\star \Omega(\cn_{X,E_\cT})^{-\cla_2}\to
\Omega(\cn_{X,E_\cT})^{-\cla_1-\cla_2}$$
that are associative and commutative.

In addition the perverse sheaves $\Omega(\cn_{X,E_\cT})^{-\cla}$
possess the following factorization
property. For $\cla=\cla_1+\cla_2$ as above, let
$\left(X^{\cla_1}\times X^{\cla_2}\right)_{disj}\subset X^{\cla_1}\times X^{\cla_2}$
be the open subset corresponding to the condition that the two divisors have
a disjoint support. 
We have a natural isomorphism:
\begin{equation} \label{factorization of omega}
\Omega(\cn_{X,E_\cT})^{-\cla}|_{\left(X^{\cla_1}\times X^{\cla_2}\right)_{disj}}\simeq
\left(\Omega(\cn_{X,E_\cT})^{-\cla_1}\boxtimes \Omega(\cn_{X,E_\cT})^{-\cla_2}\right)
|_{\left(X^{\cla_1}\times X^{\cla_2}\right)_{disj}},
\end{equation}
and similarly for $\Upsilon(\cn_{X,E_\cT})^{\cla}$.

\medskip

Set
$$R_{E_\cT}^{-\cla}:=H(X^\cla, \Omega(\cn_{X,E_\cT})^{-\cla}) \text{ and }
R_{E_\cT}:=\underset{\cla\in \cLambda^{pos}}\oplus\, R_{E_\cT}^{-\cla},\,\,
\wh{R}_{E_\cT}:=\underset{\cla\in \cLambda^{pos}}\Pi\, R_{E_\cT}^{-\cla}.$$

We obtain that $R_{E_\cT}$ is a commutative $-\cLambda^{pos}$-graded algebra,
and $\wh{R}_{E_\cT}$ is isomorphic to the completion of $R_{E_\cT}$ at the
natural augmentation ideal.

Note that the computation in \secref{gr of Omega} implies that $R_{E_\cT}$
admits a natural decreasing filtration and
$$\on{gr}(R_{E_\cT})\simeq
\Sym\left(\underset{\check\alpha\in \check\Delta^+}\bigoplus\,
H^1(X,E_\cT^{\check\alpha})^*\right)\simeq \Sym\left(H^1(X,\fn_{X,E_\cT})^*\right).$$
It is easy to see that the above filtration is given by powers of the augmentation ideal.

\medskip

We have:

\begin{thm} \label{deformation base}
We have a canonical isomorphism of topological algebras
$$\wh\CO_{\DefbE}\simeq \wh{R}_{E_\cT}.$$
\end{thm}
The proof will be given in \secref{proof of deformation base}.

\medskip

Let $\CO_{\DefbE}$ denote the algebra equal to the direct sum
of homogeneous components of $\wh\CO_{\DefbE}$ with respect
to the natural $\cT$-action. The above theorem implies:
\begin{cor} \label{graded deformation base}
$\CO_{\DefbE}\simeq R_{E_\cT}$.
\end{cor}



Let us end this section with the following remark.
In the context of D-modules, we can look at the {\it scheme},
and not just the formal scheme, of $\cB$-local systems $E_\cB$, such that
the induced $\cT$-local system is identified with $E_\cT$.

It is easy to see that the regularity assumption on $E_\cT$ implies
that this scheme indeed exists; moreover, it is isomorphic to
$\Spec(\CO_{\DefbE})$. In other words, the algebra of functions
on this scheme is isomorphic to $R_{E_\cT}$.

\section{Structure of the extension by zero}    \label{structure of open}
 
\ssec{}   \label{statement of extension}

For $\cla\in \cLambda^{pos}$ recall the maps
$$\imathb_\cla:X^\cla\times \BunBb\to \BunBb \text{ and }
\imath_\cla:X^\cla\times \Bun_B\to \BunBb.$$

In what follows we will use the following notation. For $\CT\in D^b(X^\cla)$
and $\CS\in D^b(\BunBb)$, we will denote by $\CT\star \CS\in D^b(\BunBb)$
the object $(\imathb_\cla)_!(\CT\boxtimes \CS)$. This operation is clearly
associative with respect to $\star: D^b(X^{\cla_1})\times D^b(X^{\cla_2})\to
D^b(X^{\cla_1+\cla_2})$.

The main result of this section is the following:

\begin{thm}   \label{extension by zero}  \hfill

\smallskip

\noindent{\em (1)}
The $0$-th perverse cohomology of $\imath_\cla^!\left(\jmath_!(\IC_{\Bun^\cmu_B})\right)$
is canonically isomorphic to the product
$\Omega(\cn_X)^{-\cla}\boxtimes \IC_{\Bun^{\cmu+\cla}_B}$. 
In particular, by adjunction we obtain a map
\begin{equation} \label{action upstairs}
\Omega(\cn_X)^{-\cla} \star \jmath_!(\IC_{\Bun^{\cmu+\cla}_B})\to 
\jmath_!(\IC_{\Bun^{\cmu}_B}),
\end{equation}

\smallskip

\noindent{\em (2)}
For two elements 
$\cla_1,\cla_2\in \cLambda^{pos}$ the diagram
$$
\CD
\Omega(\cn_X)^{-\cla_1}\star \Omega(\cn_X)^{-\cla_2}\star  
\jmath_!(\IC_{\Bun^{\cmu+\cla_1+\cla_2}_B}) @>>>
\Omega(\cn_X)^{-\cla_1-\cla_2}\star \jmath_!(\IC_{\Bun^{\cmu+\cla_1+\cla_2}_B}) \\
@VVV      @VVV  \\
\Omega(\cn_X)^{-\cla_1} \star \jmath_!(\IC_{\Bun^{\cmu+\cla_1}_B}) @>>>
\jmath_!(\IC_{\Bun^{\cmu}_B}).
\endCD
$$
is commutative.
\end{thm}

We shall now explain how this theorem implies point (1) of \thmref{main},
i.e., that $\Eis_!(E_\cT)$ carries an action of $\CO_{\DefbE}$. Taking
into account \corref{graded deformation base}, we need to
exhibit the maps
\begin{equation}  \label{action maps}
H(X^\cla, \Omega(\cn_{X,E_\cT})^{-\cla})\otimes \Eis^{\cmu+\cla}_!(E_\cT)\to 
\Eis^\cmu_!(E_\cT), \,\, \cla\in \cLambda^{pos}
\end{equation}
that are associative in the natural sense.

\medskip

Let us tensor both sides of \eqref{action upstairs} with $(\fqb^\cmu)^*(\CS(E_\cT))$. 
We obtain a map
$$\imath_\cla{}_!\biggl(\Omega(\cn_{X,E_\cT})^{-\cla}\boxtimes
\left(\IC_{\Bun^{\cmu+\cla}_B}\otimes (\fq^{\cmu+\cla})^*(\CS(E_\cT))\right)\biggr)\to
\jmath_!\left(\IC_{\Bun^{\cmu}_B}\otimes (\fq^{\cmu})^*(\CS(E_\cT))\right),$$
and taking the direct image with respect to $\fpb^\cmu:\BunBb^\cmu\to \Bun_G$,
we obtain the map of \eqref{action maps}, as required. The associativity
of the action follows from the commutativity of the diagram in \thmref{extension by zero}.

\medskip

Let us consider the pro-object in $\on{Perv}(\Bun_G)$ defined as
$$\wh\Eis_!(E_\cT):=\Eis_!(E_\cT)\underset{\CO_{\DefbE}}\otimes \wh\CO_{\DefbE},$$
and let us describe it in intrinsic terms.

For $\cla\in \cLambda$, let $\jmath_!^{\tau(\cla)}(\IC_{\Bun^\cmu_B})$ denote the
quotient of $\jmath_!(\IC_{\Bun^\cmu_B})$ by the image of the maps
$$\Omega(\cn_X)^{\cla'}\star \jmath_!(\IC_{\Bun^{\cmu+\cla'}_B})\to
\jmath_!(\IC_{\Bun^{\cmu}_B}),$$
given by \eqref{action upstairs} for $\cla'>\cla$. Set
$$\Eis^\cmu_!(E_\cT)^{\tau(\cla)}:=
\fpb^\cmu_!\left(\jmath_!^{\tau(\cla)}(\IC_{\Bun^\cmu_B})\otimes
(\fqb^\cmu)^*(\CS(E_\cT))\right).$$
As in \thmref{Eisenstein is perverse}, each $\Eis^\cmu_!(E_\cT)^{\tau(\cla)}$ is perverse.
Moreover, for every open substack $U\subset \Bun_G$ of finite type
there are only finitely many $\cmu$, for which the restriction
$\Eis_!^\cmu(E_\cT)^{\tau(\cla)}|_U$ is non-zero. Hence, the direct
sum $\Eis^\cmu_!(E_\cT)^{\tau(\cla)}:=\underset{\cmu}\oplus\,
\Eis^\cmu_!(E_\cT)^{\tau(\cla)}$ makes sense as an object of $\on{Perv}(\Bun_G)$.
We have: 
$$\wh\Eis_!(E_\cT)\simeq \underset{\underset{\cla}{\longleftarrow}}{"\on{lim}"}\,
\Eis^\cmu_!(E_\cT))^{\tau(\cla)}.$$

\ssec{}     \label{intr Zastava}

The rest of the present section is be devoted to the proof of \thmref{extension by zero}.

For a local system $E_\cT$ as above, let $\fU(\cn_{X,E_\cT})^\cla$ denote the 
following {\it sheaf} on $X^\cla$ for $\cla\in \cLambda^{pos}$.  Its fiber at a 
point $\Sigma\, \cla_k\cdot x_k$ with $x_k$'s distinct is the tensor product
$$\underset{k}\bigotimes\, U(\cn_{E_{\cT,x}})^{\cla_k},$$
where the superscript $\cla_k$ refers to the corresponding weight 
component in $U(\cn)$ and the subscript $E_{\cT,x}$ to the twist
by the fiber of $E_\cT$ at $x$. These fibers glue to a sheaf by means
of the co-multiplication map on $U(\cn)$.
When $E_\cT$ is trivial, we will denote the corresponding sheaf simply by
$\fU(\cn_{X})^\cla$.

One of the ingredients in the proof \thmref{extension by zero} is the following result, which
essentially follows from \cite{BFGM}:

\begin{prop}  \label{IC calculation}
There exists a canonical isomorphism in $D^b(X^\cla\times \Bun_B^{\cmu+\cla})$:
$$\imath_\cla^!(\IC_{\BunBb^\cmu})\simeq \fU(\cn_X)^\cla\boxtimes \IC_{\Bun_B^{\cmu+\cla}}.$$
\end{prop}

We should remark that the proof of \propref{IC calculation}, that we give,
uses one piece of unpublished work (see below) that has to do with
the identification of the sheaf $\fU(\cn_X)^\cla$. However, for the proof
of \thmref{extension by zero} we will need only to know the image of
$\fU(\cn_X)^\cla$ in the Grothendieck group, a computation which is
fully carried out in \cite{BFGM}.

\begin{proof}

Let us recall the construction of the Zastava spaces, $\ol{\CZ}^{\cla}$ for 
$\cla\in \cLambda^{pos}$. Let $B^-$ (resp., $N^-$) be the negative 
Borel subgroup of $G$ (resp., its unipotent radical).

By definition, $\ol{\CZ}^{\cla}$ is the open subscheme of 
$\BunBb{}^{-\cla}\underset{\Bun_G}\times \Bun_{N^-}$, corresponding
to the condition that the reduction to $N^-$ and the generalized reduction to $B$
on the given $G$-bundle are transversal at the generic point of the curve. 
The stack $\ol{\CZ}^{\cla}$ is naturally fibered over $X^\cla$ by means of a projection
denoted $\pi^\cla$, and with a section of this projection, denoted $\fs^\cla$
(see \cite{BFGM}, Sect. 2).

A basic feature of the Zastava spaces is the following factorization property
with respect to the projection $\pi^\cla$ (see \cite{BFGM}, Proposition 2.4):
\begin{equation} \label{factorization of Zastava}
\ol\CZ{}^{\cla_1+\cla_2}\underset{X^{\cla_1+\cla_2}}\times 
\left(X^{\cla_1}\times X^{\cla_2}\right)_{disj}\simeq
\left(\ol\CZ{}^{\cla_1}\times \ol\CZ{}^{\cla_2}\right)\underset{X^{\cla_1}\times X^{\cla_2}}
\times \left(X^{\cla_1}\times X^{\cla_2}\right)_{disj}.
\end{equation}

Let $\CZ{}^\cla\subset \ol{\CZ}^\cla$ denote the open subscheme,
corresponding to the condition that the $B$-structure is non-degenerate;
let $\jmath^{\CZ^\cla}$ denote its open embedding.
Both this subscheme and the section $\fs^\cla$ are compatible with the
isomorphisms \eqref{factorization of Zastava} above.

It was shown in {\it loc. cit.} that $\ol{\CZ}^{\cla}$ is locally in the smooth topology isomorphic
to $\BunBb$, in such a way that $\Bun_B\subset \BunBb$ corresponds to
$\CZ{}^\cla$, and the locally closed subvariety 
$X^\cla\times \Bun_B\subset \BunBb$ corresponds to a subscheme $\fs^\cla(X^\cla)$
that projects isomorphically onto $X^\cla$.

\medskip

We claim that to prove \propref{IC calculation}, it is sufficient to establish the isomorphism
\begin{equation} \label{IC on Z}
\fs_\cla^!(\ol{\CZ}^\cla)\simeq \fU(\cn_X)^\cla.
\end{equation}
Indeed, the local isomorphism between the triples
$$(X^\cla\times \Bun_B)\overset{\imath_\cla}\hookrightarrow 
\BunBb\overset{\jmath}\hookleftarrow \Bun_B \text{ and }
X^\cla\overset{\fs^\cla}\hookrightarrow \ol\CZ{}^\cla \overset{\jmath^{\CZ^\cla}}
\hookleftarrow \CZ^\cla$$
implies that $\imath_\cla^!(\IC_{\BunBb^\cmu})$ has the desired shape,
except for a possible twist by local systems along the $\Bun^{\cmu+\cla}_B$ 
multiple. The fact that no such twist occurs can be seen using the action of
the Hecke stack, as in \cite{BG}, Sects. 5.2 and 6.2.

\medskip

In \cite{BFGM}, Proposition 5.2,
it was also shown that there is a canonical isomorphism
$$\fs_\cla^!(\IC_{\ol{\CZ}^\cla})\simeq \pi^\cla_!(\IC_{\ol{\CZ}^\cla}),$$
and that the expression appearing in the above formula is a {\it sheaf}
isomorphic to the top (=$2|\cla|$) cohomology in the {\it usual t-structure} of
$\pi^\cla_!(\BC_{\CZ{}^\cla})$. Thus, we have to show
that the latter is isomorphic as a sheaf to $\fU(\cn_X)^\cla$.

For that we note that $\CZ{}^\cla$ is naturally a subscheme in
the Beilinson-Drinfeld Grassmannian $\Gr_{G,X^\cla}$, equal to the intersection
of the corresponding semi-infinite orbits. The identification between 
$\fU(\cn_X)^\cla$ and the top cohomology of $\pi^\cla_!(\BC_{\CZ{}^\cla})$
at the level of fibers follows from the realization of $U(\cn)$ is 
the top cohomology of the intersection of semi-infinite orbits, see 
\cite{BFGM}, Theorem 5.9.

In order to see that this identification glues to an isomorphism of
sheaves one needs to express the co-product on $U(\fn)$ in terms of
$\pi^\cla_!(\BC_{\CZ{}^\cla})$. This relationship has been recently 
established in \cite{Kam}

\end{proof}

To state a corollary of the above proposition that will be used in the
proof of \thmref{extension by zero}, let us consider the following 
version of the Grothendieck group of perverse sheaves on $\BunBb^\cmu$.
We start with the usual Grothendieck group of the Artinian category of
perverse sheaves on $\BunBb$ of {\it finite length}, and we complete
it with respect to the topology, where the system of neighbourhoods
of zero is given by classes of perverse sheaves, supported on 
closures of $\imath_\cla(X^\cla\times \Bun_B)$, $\cla\in \cLambda^{pos}$. 
In particular, the class of each $(\imath_\cla)_!(\CT)$, 
$\CT\in D^b(X^\cla\times \Bun^{\cmu+\cla}_B)$
is a well-defined element of this group.

\begin{cor}   \label{Groth group comp}
There is an equality $[\jmath_!(\IC_{\Bun^\cmu_B})]=\underset{\cla\in \cLambda^{pos}}\Sigma\,
[\Omega(\cn_X)^{-\cla}\star \IC_{\BunBb^{\cmu+\cla}}]$.
\end{cor}

\begin{proof}

From \propref{IC calculation} it follows that for $\CT\in D^b(X^\cla)$,
$$[\CT\star \IC_{\BunBb^{\cmu+\cla}}]=\underset{\cla'\in \cLambda^{pos}}\Sigma\,
[\CT\star \fU(\cn_X)^{\cla'}\star \jmath_!(\IC_{\Bun^{\cmu+\cla+\cla'}_B})].$$

Moreover, by \secref{abstract Koszul}, for $0\neq \cla\in \cLambda^{pos}$,
$$\underset{\cla_1+\cla_2=\cla}{\underset{\cla_1,\cla_2\in \cLambda^{pos}}\Sigma}\,
[\fU(\cn_X)^{\cla_1}\star \Omega(\cn_X)^{-\cla_2}]=0.$$

This implies the assertion of the corollary.
\end{proof}

\begin{cor}  \label{Groth for h 0}
In the Grothendieck group of perverse sheaves 
on $X^\cla\times \Bun_B^{\cmu+\cla}$ we have the following equality:
$$[h^0\left(\imath_\cla^!\left(\jmath_!(\IC_{\Bun_B^\cmu})\right)\right)]=
[\Omega(\cn_X)^{-\cla} \boxtimes \IC_{\Bun^{\cmu+\cla}_B}].$$
\end{cor}

\begin{proof}

Let $\BunBb^{\cmu,\leq\cla}$ by the open substack of $\BunBb^\cmu$,
obtained by removing the closed substack equal to the
union $\imath_{\cla'}(X^{\cla'}\times \Bun_B^{\cmu+\cla'})$ for 
$\cla'-\cla\in \cLambda^{pos}-0$.

Arguing as in the proof of \propref{IC calculation}, we can replace
the original question about $\BunBb^{\cmu,\leq \cla}$ for one about the
Zastava space $\ol\CZ{}^\cla$. Using the factorization property 
\eqref{factorization of Zastava}, and arguing by induction on $|\cla|$, 
can assume that the desired equality holds in the Grothendieck group
of perverse sheaves over the open substack 
$(X^\cla-\Delta(X))\times \Bun_B^{\cmu+\cla}$, where 
$\Delta(X)\subset X^\cla$ denotes the main diagonal.

\medskip

Taking into account \corref{Groth group comp}, we have to show that there
does not exist $\cla'\in \cLambda^{pos}$ with $0\neq \cla'\neq\cla$ and a perverse sheaf $\CT$ on $X^{\cla'}$, 
appearing as subquotient of $\Omega(\cn_X)^{-\cla'}$, and a non-trivial extension of 
perverse sheaves on $\BunBb^{\cmu,\leq\cla}$:
$$0\to (\CT\star \IC_{\BunBb^{\cmu+\cla'}})|_{\BunBb^{\cmu,\leq\cla}}\to \CT'\to
(\Delta_!(\BC_X)[1]\boxtimes \IC_{\Bun^{\cmu+\cla}_B})\to 0.$$

Obviously, an extension as above does not exist unless $\cla-\cla'\in \Lambda^{pos}$.
In the latter case, it would be given by a morphism from $\Delta_!(\BC_X[1])
\boxtimes \IC_{\Bun^{\cmu+\cla}_B}$ to
$$h^1\biggl(\imath_\cla^!\left(\imathb_{\cla'}{}_!(\CT\boxtimes 
\IC_{\BunBb^{\cmu+\cla'}})\right)\biggr).$$
By \propref{IC calculation}, the latter expression is isomorphic to
$$h^1\biggl((\CT\star \fU(\cn_X)^{\cla-\cla'})\boxtimes \IC_{\Bun^{\cmu+\cla}_B}\biggr).$$

Note that $\fU(\cn_X)^{\cla-\cla'}$ is concentrated in the perverse cohomological
degrees $\geq 1$. The assertion of the corollary follows now from the fact that
$\Delta^!(\CT\star \fU(\cn_X)^{\cla-\cla'})$ is concentrated in the cohomological
degrees $\geq 2$.

\end{proof}



\ssec{}

Our present goal is to construct a map
$$\Omega(\cn_X)^{-\cla} \star \jmath_!(\IC_{\Bun^{\cmu+\cla}_B}) \to
\jmath_!(\IC_{\Bun^{\cmu}_B}),$$
or, equivalently, a map
\begin{equation} \label{need to construct}
\Omega(\cn_X)^{-\cla} \boxtimes \IC_{\Bun^{\cmu+\cla}_B}\to
h^0\left(\imath_\cla^!\left(\jmath_!(\IC_{\Bun_B^\cmu})\right)\right).
\end{equation}

Let $\overset{\circ}X{}^\cla$ be the open subset of $X^\cla$,
corresponding to the {\it full} partition. I.e., it corresponds
to coloured divisors of the form $\Sigma\, \cla_k\cdot x_k$
with $x_k$'s pairwise distinct and each $\cla_k$ being a 
simple coroot. Let $j^\cla$ denote the open embedding
$\overset{\circ}X{}^\cla\hookrightarrow X^\cla$. 

First, we claim that here exists a map (in fact, an isomorphism)
\begin{equation} \label{map over open}
j^\cla{}^*\left(\Omega(\cn_X)^{-\cla}\right)\boxtimes \IC_{\Bun^{\cmu+\cla}_B}\simeq
(j^\cla\times \on{id})^*\biggl(h^0\left(\imath_\cla^!\circ \jmath_!(\IC_{\Bun_B^\cmu})\right)\biggr)
\end{equation}
over $\overset{\circ}X{}^\cla\times \Bun_B^{\cmu+\cla}$. In other words, we claim
that the isomorphism stated in \thmref{extension by zero} holds over the open
substack $\overset{\circ}X{}^\cla\times \Bun_B^{\cmu+\cla}$.

\begin{proof}

As in the proof of \propref{IC calculation}, the assertion reduces to one
about the Zastava space. Namely, we have to show that the restriction of
$$h^0\left((\fs^\cla)^!\circ \jmath^{\CZ^\cla}_!(\IC_{\CZ^\cla})\right)$$
to $\overset{\circ}X{}^\cla$ is isomorphic to the restriction of $\Omega(\cn_X)^{-\cla}$.

Write $\cla=\underset{i\in I}\Sigma\, n_i\cdot \check\alpha_i$. Then
$$\overset{\circ}X{}^\cla\simeq \left(\underset{i\in I}\Pi\, X^{(n_i)}\right)_{disj},$$
and the restriction of $\Omega(\cn_X)^{-\cla}$ to it isomorphic to
$$\underset{i\in I}\boxtimes\, \Lambda^{(n_i)}((\cn^*)^{-\check\alpha_i}_X[1]),$$
where $(\cn^*)^{-\check\alpha_i}_X$ denotes the constant sheaf on $X$ with
fiber $(\cn^*)^{-\check\alpha_i}$, and $\Lambda^{(n_i)}\left(\cdot\right)$ the
external exterior power of a local system. 

\medskip

By the same argument as in \cite{BFGM}, Proposition 5.2, 
$(\fs^\cla)^!\circ \jmath^{\CZ^\cla}_!(\IC_{\CZ^\cla})\simeq 
(\pi^\cla\circ \jmath^{\CZ^\cla})_!(\IC_{\CZ^\cla})$. So, we have to 
calculate the $0$-th perverse cohomology
of the direct image with compact supports of the constant sheaf on $\CZ^\cla$,
cohomologically shifted by $[2|\cla|]$.

However, by the factorization property given by \eqref{factorization of Zastava},
$$\CZ^\cla\underset{X^\cla}\times \overset{\circ}X{}^\cla\simeq
\underset{i\in I}\Pi\, \left(\CZ^{\check\alpha_i}\right)^{\times n_i}/\Sigma_{n_i},$$
where $\Sigma_{n_i}$ is the corresponding symmetric group. Moreover, each
$\CZ^{\check\alpha_i}$ is isomorphic to the product $X\times \BG_m$,
and 
$$h^0\left((\pi^{\check\alpha_i} \circ \jmath^{\CZ^{\check\alpha_i}})_!
(\IC_{\CZ^{\check\alpha_i}})\right)\simeq \BC_X[1]\otimes H^1(\BG_m,\BC).$$

This makes the required isomorphism manifest once we identify each of the lines
$(\cn^*)^{-\check\alpha_i}$ with $H^1(\BG_m,\BC)$. There exists a natural identification
like this, when we realize $U(\cn)^\cla$ as the top cohomology of the corresponding
intersection in the affine Grassmannian.

\end{proof}

Thus, our task is to show that the isomorphism \eqref{map over open} extends
to a map (and an isomorphism) over the entire $X^\cla\times \Bun_B^{\cmu+\cla}$.

\begin{lem}  \label{injectivity for omega}
The canonical map
$$\Omega(\cn_X)^{-\cla}\to j^\cla_*\circ j^\cla{}^*\left(\Omega(\cn_X)^{-\cla})\right)$$
is injective.
\end{lem}

\begin{proof}
Using \eqref{factorization of omega} we can apply induction on $|\cla|$. The
assertion clearly holds for $|\cla|=1$, i.e., when $\cla$ is a simple root. Thus,
we can assume that $|\cla|\geq 2$ and that the injectivity assertion
holds over $X^\cla-\Delta(X)$. 

We have to show that $\Delta^!(\Omega(\cn_X)^{-\cla})$
has no cohomologies in degrees $\leq 0$. This amounts to the fact that the complex
$\bC^\bullet(\cn)^{-\cla}$ has no cohomologies in degrees $\leq 1$. But this is evidently
so: the kernel of the co-bracket $\cn^*\to \Lambda^2(\cn^*)$ is spanned by the duals of
the simple roots.
\end{proof}

The main geometric ingredient in the construction of the map 
\eqref{need to construct} is the following:

\begin{prop}   \label{injectivity for !}
The map of perverse sheaves on $X^\cla\times \Bun_B^{\cmu+\cla}$
$$h^0\left(\imath_\cla^!\circ \jmath_!(\IC_{\Bun_B^\cmu})\right)\to
(j^\cla\times \on{id})_*\circ 
(j^\cla\times \on{id})^*\biggl(h^0\left(\imath_\cla^!\circ \jmath_!(\IC_{\Bun_B^\cmu})\right)\biggr)$$
is injective.
\end{prop}

Let us assume this proposition, and proceed with the construction of the map
\eqref{need to construct}. \footnote{Added in Jan. 2008: a simpler proof
of this proposition has been found that does not rely on \cite{FFKM} and makes
\secref{FFKM constr} redundant. Namely, one can argue by induction on $|\cla|$
analyzing $\Delta^!\circ \pi_!^\cla \circ \jmath_!(\IC_{\CZ^\cla})$.}

\ssec{}

By \propref{injectivity for !}, we obtain that if the map \eqref{map over open} extends to 
a map as in \eqref{need to construct}, then it does so uniquely. 
Furthermore, by \lemref{injectivity for omega}, the latter map is automatically injective. 
We claim that in this case, it is surjective as well, implying the isomorphism statement of 
\thmref{extension by zero}. Indeed, the equivalence of injectivity and surjectivity
properties of the map in question follows immediately from 
\corref{Groth for h 0}. 

The commutativity of the diagram of 
\thmref{extension by zero} also follows from \lemref{injectivity for omega}.

\medskip

Thus, let us assume by induction that the map \eqref{map over open} has been shown to
extend to a map \eqref{need to construct} for all parameters $\cla'$ with $|\cla'|<|\cla|$. 
Since the the question of extension is local, we can pass to the Zastava space $\ol\CZ{}^\cla$
as in the proof of \propref{IC calculation}. 

Using the factorization property \eqref{factorization of Zastava},
and by the induction hypothesis, we can assume that the map
\eqref{need to construct} has been extended over 
$(X^\cla-\Delta(X))\times \Bun_B^{\cmu+\cla}$. Let us distinguish two cases.

\medskip

\noindent Case 1, when $\cla$ is not a root of $\cg$. In this case, by
\corref{Groth for h 0} and \secref{gr of Omega}, we obtain that 
neither $\Omega(\cn_X)^{-\cla} \boxtimes \IC_{\Bun^{\cmu+\cla}_B}$
nor $h^0\left(\imath_\cla^!\left(\jmath_!(\IC_{\Bun_B^\cmu})\right)\right)$ has
sub-quotients supported on $\Delta(X)\times \Bun^{\cmu+\cla}_B$. Hence, both
sides of \eqref{need to construct} are the minimal extensions of their respective
restrictions to $(X^\cla-\Delta^!(X))\times \Bun_B^{\cmu+\cla}$, and the extension
assertion follows by the functoriality of the minimal extension operation. 

\bigskip

\noindent Case 2, when $\cla$ is a root of $\cg$. Consider the sum of the images of
$\Omega(\cn_X)^{-\cla} \boxtimes \IC_{\Bun^{\cmu+\cla}_B}$ and 
$h^0\left(\imath_\cla^!\left(\jmath_!(\IC_{\Bun_B^\cmu})\right)\right)$ in
$$j^\cla_*\circ j^\cla{}^*\left(\Omega(\cn_X)^{-\cla}\right)\boxtimes \IC_{\Bun^{\cmu+\cla}_B}\simeq
(j^\cla\times \on{id})_*\circ
(j^\cla\times \on{id})^*\biggl(h^0\left(\imath_\cla^!\circ \jmath_!(\IC_{\Bun_B^\cmu})\right)\biggr).$$
Let us denote this perverse sheaf by $\CT$. We claim that 
$\Omega(\cn_X)^{-\cla} \boxtimes \IC_{\Bun^{\cmu+\cla}_B}$ maps isomorphically to $\CT$.
Indeed, if it did not, we would have a non-trivial extension
$$0\to \Omega(\cn_X)^{-\cla} \boxtimes \IC_{\Bun^{\cmu+\cla}_B}\to \CT\to
\Delta_!(\BC_X[1])\boxtimes \IC_{\Bun^{\cmu+\cla}_B}\to 0.$$
However, we have:
\begin{lem}
$\on{Ext}^1_{X^\cla}(\BC_{\Delta(X)}[1],\Omega(\cn_X)^{-\cla})=0$
for $\cla$ being a root of $\cg$.
\end{lem}

\begin{proof}
By adjunction, it suffices to show that $h^1\left(\Delta^!(\Omega(\cn_X)^{-\cla})\right)=0$.
This amounts to the calculation of the component of weight $-\cla$ in
$H^2(\cn,\BC)$. It is known that the weights that appear in $H^2(\cn,\BC)$ 
are of the form $-\cla=-\check\alpha_i-s_i(\check\alpha_j)$ for pairs of simple roots
$\alpha_i$ and $\alpha_j$. However, a weight $\cla$ of the above form is
never a root.
\end{proof}

Thus, we obtain that 
$h^0\left(\imath_\cla^!\left(\jmath_!(\IC_{\Bun_B^\cmu})\right)\right)$ is
a sub-perverse sheaf of $\Omega(\cn_X)^{-\cla} \boxtimes \IC_{\Bun^{\cmu+\cla}_B}$.
The fact that this inclusion is an isomorphism follows from 
\corref{Groth for h 0}.

\section{The FFKM construction and proof of \propref{injectivity for !}} 
\label{FFKM constr}

\ssec{}

The goal of this section is to prove \propref{injectivity for !}, which is, in a way,
the main technical point of this paper. We will first make a reduction
to the case when $\cla$ is a root of $\cg$, and in the latter case we will
deduce the required assertion from a construction of \cite{FFKM}.

\ssec{}

We argue by induction on $|\cla|$, and we claim that we can assume 
that the morphism 
\begin{equation} \label{map in question}
h^0\left(\imath_\cla^!\circ \jmath_!(\IC_{\Bun_B^\cmu})\right)\to
(j^\cla\times \on{id})_*\circ (j^\cla\times \on{id})^*
\biggl(h^0\left(\imath_\cla^!\circ \jmath_!(\IC_{\Bun_B^\cmu})\right)\biggr)
\end{equation}
is injective over the open substack $(X^\cla-\Delta(X))\times \Bun_B^{\cmu+\cla}$.

Indeed, the injectivity statement is local, so we can replace $\BunBb^{\cmu,\leq \cla}$
by the Zastava space $\ol\CZ{}^\cla$, and apply the factorization property 
\eqref{factorization of Zastava}.

This reduces the assertion of the proposition to the following one:
\begin{prop}  \label{injectivity for !, reform}
Assume that $|\cla|>1$. Then the perverse sheaf 
$h^0\left(\imath_\cla^!\circ \jmath_!(\IC_{\Bun_B^\cmu})\right)$
does not have sub-objects supported on $\Delta(X)\times \Bun_B^{\cmu+\cla}$.
\end{prop}

The rest of this section is essentially devoted to the proof of this proposition.
Let us assume first that $\cla$ is not a root of $\cg$. Then by
\corref{Groth for h 0} and \secref{gr of Omega}, the Jordan-H\"older series of the 
perverse sheaf $h^0\left(\imath_\cla^!\circ \jmath_!(\IC_{\Bun_B^\cmu})\right)$ does
not contain terms that are supported on the closed substack
$\Delta(X)\times \Bun_B^{\cmu+\cla}$,
implying, in particular, that it does not have such sub-objects.

\ssec{}    \label{meat}

Hence, it remains to analyze the case when $\cla$ is a root $\check\alpha$, but not a 
simple root. We will argue by contradiction, assuming that 
$h^0\left(\imath_{\check\alpha}^!\circ \jmath_!(\IC_{\Bun_B^\cmu})\right)$ admits 
$\Delta_!(\BC_X)[1]\boxtimes \IC_{\Bun_B^{\cmu+\check\alpha}}$ as a sub-object.

\medskip

Consider the quotient 
\begin{equation} \label{bad quotient}
h^0\left(\imath_{\check\alpha}^!\circ \jmath_!(\IC_{\Bun_B^\cmu})\right)/
\left(\Delta_!(\BC_X)[1]\boxtimes \IC_{\Bun_B^{\cmu+\check\alpha}}\right).
\end{equation}
By \corref{Groth for h 0}, it is isomorphic to the intermediate extension of its own 
restriction to the open substack $(X^\cla-\Delta(X))\times \Bun_B^{\cmu+\cla}$. 
By the induction hypothesis
and the above factorization argument, we can assume that \thmref{extension by zero} 
holds over $(X^\cla-\Delta(X))\times \Bun_B^{\cmu+\cla}$. Hence, the above
quotient perverse sheaf is isomorphic to
$$\biggl(\on{ker}\left(\Omega(\cn_X)^{-\check\alpha}\to 
\Delta_!((\cn^*)^{-\check\alpha}_X[1])\right)\biggr)\boxtimes \IC_{\Bun_B^{\cmu+\check\alpha}},$$
since the first multiple in the above formula equals the intermediate extension
of the restriction of $\Omega(\cn_X)^{-\check\alpha}$ to $X^\cla-\Delta(X)$.

\medskip

Let $\check\beta$ and $\check\gamma$ be two roots such that
$\check\alpha=[\check\beta,\check\gamma]$. By \secref{gr of Omega}, the perverse
sheaf of \eqref{bad quotient} admits a further quotient, isomorphic to 
\begin{equation} \label{two root quotient}
\CF'_1:=\biggl(\left((\cn^*)^{-\check\beta}_X[1]\right)\star 
\left((\cn^*)^{-\check\gamma}_X[1]\right)\biggr)
\boxtimes \IC_{\Bun_B^{\cmu+\check\alpha}}.
\end{equation}
Let $\CF'$ denote the corresponding quotient of 
$\jmath_!(\IC_{\Bun_B^\cmu})|_{\BunBb^{\cmu,\leq \check\alpha}}$.
This perverse sheaf has a 3-step filtration with $\CF'_1$ being the above
perverse sheaf on $X^{\check\alpha}\times \Bun_B^{\cmu+\check\alpha}$,
thought of as a closed substack of $\BunBb^{\cmu,\leq \check\alpha}$,
and $$\CF'_2:=
\on{ker}(\CF'\to \IC_{\BunBb^\cmu}|_{\BunBb^{\cmu,\leq \check\alpha}}),$$
so that
$\CF'_3/\CF'_2\simeq \IC_{\BunBb^\cmu}|_{\BunBb^{\cmu,\leq \check\alpha}}$
and 
$$\CF'_2/\CF'_1\simeq \left(\on{ker}(\jmath_!(\IC_{\Bun_B^\cmu})\to 
\IC_{\BunBb^\cmu})|_{\BunBb^{\cmu,\leq \check\alpha}}\right)
/\left((\imath_{\check\alpha})_!
h^0\left(\imath_{\check\alpha}^!\circ \jmath_!(\IC_{\Bun_B^\cmu})\right)\right).$$

\medskip

By \corref{Groth group comp} and \propref{IC calculation}, the only terms in the 
Jordan-H\"older series of $\jmath_!(\IC_{\BunBb^\cmu})|_{\IC_{\BunBb^{\cmu,\leq \check\alpha}}}$
that can have a non-trivial $\on{Ext}^1$ to the perverse sheaf \eqref{two root quotient}
are 
\begin{equation} \label{two perv sheaves}
\left((\cn^*)^{-\check\beta}_X[1]\right)\star \IC_{\BunBb^{\cmu+\check\gamma}}
|_{\BunBb^{\cmu,\leq \check\alpha}} \text{ and } 
\left((\cn^*)^{-\check\gamma}_X[1]\right)\star \IC_{\BunBb^{\cmu+\check\beta}}
|_{\BunBb^{\cmu,\leq \check\alpha}}.
\end{equation}

Moreover, no other terms in the Jordan-H\"older series of 
$\jmath_!(\IC_{\BunBb^\cmu})|_{\IC_{\BunBb^{\cmu,\leq \check\alpha}}}$
apart from $\IC_{\BunBb^\cmu}|_{\BunBb^{\cmu,\leq \check\alpha}}$ admit
a non-trivial $\on{Ext}^1$ to either of the perverse sheaves appearing in 
\eqref{two perv sheaves}.

Thus, we obtain that the perverse sheaf $\CF'$ on $\BunBb^{\cmu+\check\beta}$
admits a quotient that we shall denote by $\CF$, endowed with a 3-step filtration
$$0=\CF_0\subset \CF_1\subset \CF_2\subset \CF_3=\CF,$$
such that 
$$\CF_1\simeq \CF'_1\simeq 
\left(\BC^{\check\beta}_X[1]\star \BC^{\check\gamma}_X[1]\right)\star \IC_{\Bun_B^{\cmu+\check\alpha}},$$
$$\CF_2/\CF_1\simeq \BC^{\check\beta}_X[1]\star \IC_{\BunBb^{\cmu+\check\gamma}}
|_{\BunBb^{\cmu,\leq \check\alpha}} \bigoplus
\BC^{\check\gamma}_X[1]\star \IC_{\BunBb^{\cmu+\check\beta}}
|_{\BunBb^{\cmu,\leq \check\alpha}}$$
and
$\CF_3/\CF_2\simeq \IC_{\BunBb^\cmu}|_{\BunBb^{\cmu,\leq \check\alpha}}$.

The corresponding elements in
$$\on{Ext}^1_{\BunBb^{\cmu,\leq \check\alpha}}\left(\IC_{\BunBb^\cmu},
\BC^{\check\beta}_X[1]\star \IC_{\BunBb^{\cmu+\check\gamma}}\right), \text{ and }
\on{Ext}^1_{\BunBb^{\cmu,\leq \check\alpha}}\left(\IC_{\BunBb^\cmu},
\left(\BC^{\check\gamma}_X[1]\right)\star \IC_{\BunBb^{\cmu+\check\beta}}\right)$$
are non-zero, and correspond, therefore, to nonzero multiples of the maps
$$\BC_X\to \fU(\cn_X)^{\check\beta} \text{ and }
\BC_X\to \fU(\cn_X)^{\check\gamma},$$
given by $\cn^{\check\beta}\to U(\cn)^{\check\beta}$ and
$\cn^{\check\gamma}\to U(\cn)^{\check\gamma}$, respectively.

The extension
$$0\to \CF_1\to \CF_2\to \CF_2/\CF_1\to 0$$
is non-trivial, since otherwise we would obtain that there exist
a non-trivial element in
$$\on{Ext}^1_{\BunBb^{\cmu,\leq \check\alpha}}\left(\IC_{\BunBb^\cmu},
\left(\BC^{\check\beta}_X[1]\star 
\BC^{\check\gamma}_X[1]\right)\star \IC_{\Bun_B^{\cmu+\check\alpha}}\right),$$
which is impossible by \propref{IC calculation}.

We are going to show now, using results of \cite{FFKM} that
a perverse sheaf $\CF$ with a filtration having such properties does not
exist.

\ssec{}

Recall the {\it sheaf} $\fU(\cn_X)^\cla$ on $X^\cla$. The multiplication operation
on $U(\cn)$ defines a map
\begin{equation} \label{mult U}
\fU(\cn_X)^{\cla_1}\star \fU(\cn_X)^{\cla_2}\to \fU(\cn_X)^{\cla_1+\cla_2}.
\end{equation}

The result of \cite{FFKM}, Sect. 2 that we will use can be summarized as follows:

\begin{thm} \label{FFKM thm}  \hfill

\smallskip

\noindent{\em (1)}
For every $\cla\in \cLambda^{pos}$ there exists an isomorphism
$\fU(\cn_X)^\cla\boxtimes \IC_{\Bun_B^{\cmu+\cla}}\simeq
\imath_\cla^!(\IC_{\BunBb^\cmu})$.

\medskip

\noindent{\em (2)}
The above isomorphism extends to a unique morphism
$\fU(\cn_X)^\cla\boxtimes \IC_{\BunBb^{\cmu+\cla}}\to 
\imathb_\cla^!(\IC_{\BunBb^\cmu})$,
or, equivalently (by adjunction), to a morphism
$$\fU(\cn_X)^\cla\star \IC_{\BunBb^{\cmu+\cla}} \to \IC_{\BunBb^\cmu}.$$

\medskip

\noindent{\em (3)}
For $\cla=\cla_1+\cla_2$ the diagram 
$$
\CD
\fU(\cn_X)^{\cla_1}\star \fU(\cn_X)^{\cla_2}\star
\IC_{\BunBb^{\cmu+\cla}} @>>> \fU(\cn_X)^\cla\star \IC_{\BunBb^{\cmu+\cla}}  \\
@VVV    @VVV   \\
\fU(\cn_X)^{\cla_1}\star \IC_{\BunBb^{\cmu+\cla_1}} @>>> \IC_{\BunBb^\cmu},
\endCD
$$
commutes.
\end{thm}

Let us add several remarks. First, comparing point (1) of
\thmref{FFKM thm} and that of \propref{IC calculation}, we obtain that
there are {\it a priori} two isomorphisms
$$\fU(\cn_X)^\cla\boxtimes \IC_{\Bun_B^{\cmu+\cla}}\rightrightarrows
\imath_\cla^!(\IC_{\BunBb^\cmu}).$$
At this stage it is not clear why these two maps coincide. 

The existence of the map stated in point (2) of the theorem is proved
in \cite{FFKM} by purity considerations. The statement about 
uniqueness of this extension (which is omitted in {\it loc. cit.})
follows by analyzing cohomological degrees of various subquotients.

\ssec{}

We shall analyze the commutative diagram of \thmref{FFKM thm} in
the following particular case. 

Let us first take $\cla_1=\check\beta$ and $\cla_2=\check\gamma$,
and identify $\cn^{\check\beta}_X$ and $\cn^{\check\gamma}_X$
with $\BC_X$ (up to a scalar). We obtain extensions
$$0\to \BC_X[1]^{\check\gamma}\star 
\IC_{\BunBb^{\cmu+\check\gamma}}|_{\BunBb^{\cmu,\leq \check\alpha}}
\to \wt\CF_{2,\beta}\to \IC_{\BunBb^{\cmu}}
|_{\BunBb^{\cmu,\leq \check\alpha}}\to 0$$
and
\begin{equation}  \label{second extension}
0\to \BC_X[1]^{\check\beta}\star 
\BC_X[1]^{\check\gamma}\star \IC_{\Bun_B^{\cmu+\check\alpha}}
\to \wt\CF_{1,\beta}\to \BC_X[1]^{\check\gamma}\star 
\IC_{\BunBb^{\cmu+\check\gamma}}|_{\BunBb^{\cmu,\leq \check\alpha}}\to 0,
\end{equation}
whose cup product is an element in
\begin{equation} \label{Ext 2}
\on{Ext}^2_{\BunBb^{\cmu,\leq \check\alpha}}\left(\IC_{\BunBb^{\cmu}},
\BC_X[1]^{\check\beta}\star 
\BC_X[1]^{\check\gamma}\star \IC_{\Bun_B^{\cmu+\check\alpha}}\right),
\end{equation}
corresponding to the map
$$\BC_{X\times X}\to \cn^{\check\beta}_X\star \cn^{\check\gamma}_X\to
\fU(\cn_X)^{\check\alpha}.$$

Interchanging the roles of $\check\beta$ and $\check\gamma$ we obtain
another element in the above $\on{Ext}^2$ group. The difference of these
two elements is non-zero for any choice of non-zero scalars as above, since
the roots $\check\beta$ and $\check\gamma$ do not commute.

\medskip

Going back to the
perverse sheaf $\CF$ of \secref{meat} we obtain that the corresponding
element in 
$$\on{Ext}^1_{\BunBb^{\cmu,\leq \check\alpha}}\left(
\left((\cn^*)^{-\check\beta}_X[1]\right)\star \IC_{\BunBb^{\cmu+\check\gamma}},\CF_1\right)$$
is non-zero, and hence equals, up to a scalar to the one of 
\eqref{second extension}. A similar assertion holds for $\check\beta$ replaced
by $\check\gamma$.

Let us now interpret what the existence of a sheaf $\CF$ with a 3-step extension 
having the properties specified above would mean:

It implies that that the difference of the two resulting elements in
\eqref{Ext 2} is zero, contradicting the assertion made earlier.

\section{The Koszul complex}   \label{Koszul section}

\ssec{}

In the previous sections we have endowed $\Eis_!(E_\cT)$ with an action of
$\CO_{\DefbE}$. The goal of this section is to show that
$$\BC\overset{L}{\underset{\CO_{\DefbE}}\otimes}\Eis_!(E_\cT)\simeq
\Eisb(E_\cT),$$
thereby proving point (2) of \thmref{main}.

In order to do this we shall first construct a certain Koszul complex, by
means of which one can compute $\BC\overset{L}{\underset{\CO_{\DefbE}}\otimes}\CM$
for any $\CO_{\DefbE}$-module $\CM$.

\ssec{}     

Recall again the {\it sheaf} $\fU(\cn_{X,E_\cT})^\cla\in D^b(X^\cla)$, introduced
in \secref{intr Zastava}. We are now going to represent it 
by an explicit complex of perverse sheaves.  

\medskip

For a positive integer $m$ consider the direct sum 
$$\fU(\cn_{X,E_\cT})^{m,\cla}:=
\underset{\underset{\cla_j\neq 0, \Sigma \cla_j=\cla}{\cla_1,...,\cla_m\in \cLambda^{pos}}}
\bigoplus\, \Upsilon(\cn_{X,E_\cT})^{\cla_1}\star...\star \Upsilon(\cn_{X,E_\cT})^{\cla_m},$$
where $\Upsilon(\cn_{X,E_\cT})^{\cmu}$'s are as in \secref{intr Upsilon}.
Note that we have natural maps
\begin{equation}  \label{add on U dot}
\Upsilon(\cn_{X,E_\cT})^{\cla'}\star \fU(\cn_{X,E_\cT})^{m,\cla''}\to
\fU(\cn_{X,E_\cT})^{m+1,\cla'+\cla''}.
\end{equation}

We define a differential $\fd_\fU^{m,\cla}:\fU(\cn_{X,E_\cT})^{m,\cla}\to \fU(\cn_{X,E_\cT})^{m+1,\cla}$
as follows. By induction, assume that
$$\fd_{\fU}^{m-1,\cla-\cla_1}:\Upsilon(\cn_{X,E_\cT})^{\cla_2}\star...\star \Upsilon(\cn_{X,E_\cT})^{\cla_m}\to
\fU(\cn_{X,E_\cT})^{m,\cla-\cla_1}$$
has been defined. We let $\fd_\fU^{m,\cla}$ be the sum of
$$\on{id}_{\Upsilon(\cn_{X,E_\cT})^{\cla_1}}\star \fd_{\fU}^{m-1,\cla-\cla_1}$$
and 
\begin{multline*}
\Upsilon(\cn_{X,E_\cT})^{\cla_1} \star \left(\underset{j=2,...,m}\star \, 
\Upsilon(\cn_{X,E_\cT})^{\cla_j}\right)\to  \\
\to \underset{\underset{\cla'_1,\cla''_1\neq 0}{\cla'_1+\cla''_1=\cla_1}}\bigoplus\, 
\Upsilon(\cn_{X,E_\cT})^{\cla'_1}\star \Upsilon(\cn_{X,E_\cT})^{\cla''_1}
\star \left(\underset{j=2,...,m}\star \, \Upsilon(\cn_{X,E_\cT})^{\cla_j}\right),
\end{multline*}
coming from \eqref{comult upsilon}. It is straightforward to check that 
$\fd_\fU^{m,\cla}\circ \fd_\fU^{m-1,\cla}=0$, so we can form a complex
of perverse sheaves that we will denote by $\fU(\cn_{X,E_\cT})^{\bullet,\cla}$. 

\begin{lem}  \label{dot and without}
We have an isomorphism in $D^b(X^\cla)$:
$$\fU(\cn_{X,E_\cT})^{\bullet,\cla}\simeq \fU(\cn_{X,E_\cT})^{\cla}.$$
\end{lem}

\begin{proof}

Let us compute the fiber of $\fU(\cn_{X,E_\cT})^{\bullet,\cla}$ at a point
$\Sigma\, \cla_k\cdot x_k$ of $X^\cla$. It is easy to see that this fiber
will be isomorphic, as a complex, to the product of the corresponding
complexes for the $x_k$'s; so we can assume that our point is $\cla\cdot x$.

In this case the resulting complex is quasi-isomorphic to the complex
associated to the bi-complex, whose $m$-th column is
$$\underset{\underset{\cla_j\neq 0, \Sigma \cla_j=\cla}{\cla_1,...,\cla_m\in \cLambda^{pos}}}
\bigoplus\, \left(\Lambda^\bullet(\cn_{E_{\cT,x}})\right)^{\cla_1}\otimes...
\otimes \left(\Lambda^\bullet(\cn_{E_{\cT,x}})\right)^{\cla_m}.$$
The vertical differential in this bi-complex comes from the Lie algebra structure 
on $\cn$, and the horizontal one is defined as in the case of 
$\fU(\cn_{X,E_\cT})^{m,\cla}$ via the co-multiplication on 
$\Lambda^\bullet(\cn_{E_{\cT,x}})$.

The $0$-th term of this complex is isomorphic to
$$\underset{\underset{\cla_j\neq 0,\Sigma \cla_j=\cla}{m,\cla_1,...,\cla_m\in \cLambda^{pos}}}
\bigoplus\, \cn^{\cla_1}_{E_{\cT,x}}\otimes...\otimes \cn^{\cla_m}_{E_{\cT,x}},$$
and it maps to $U(\cn_{E_{\cT,x}})^\cla$ via the product operation
in this algebra.

This map is easily seen to induce a quasi-isomorphism from the above
complex to $U(\cn_{E_{\cT,x}})^\cla$. In addition, it is straightforward
to check that the above fiber-wise calculation identifies 
$\fU(\cn_{X,E_\cT})^{\cla}$ with the cohomology of $\fU(\cn_{X,E_\cT})^{\bullet,\cla}$.

\end{proof}

The complexes $\fU(\cn_{X,E_\cT})^{\bullet,\cla}$ possess the following structures.
For $\cla=\cla_1+\cla_2$ there is a multiplication map
\begin{equation} \label{mult dot U}
\fU(\cn_{X,E_\cT})^{\bullet,\cla_1}\star \fU(\cn_{X,E_\cT})^{\bullet,\cla_2}\to
\fU(\cn_{X,E_\cT})^{\bullet,\cla},
\end{equation}
inducing the map \eqref{mult U}. 
In addition, they have a factorization property similar
to that of \eqref{factorization of omega}:
\begin{equation} \label{factorization of U}
\fU(\cn_{X,E_\cT})^{\bullet,\cla}|_{\left(X^{\cla_1}\times X^{\cla_2}\right)_{disj}}\simeq
\left(\fU(\cn_{X,E_\cT})^{\bullet,\cla_1}\boxtimes \fU(\cn_{X,E_\cT})^{\bullet,\cla_2}\right)
|_{\left(X^{\cla_1}\times X^{\cla_2}\right)_{disj}},
\end{equation}
also compatible with the corresponding isomorphism for $\fU(\cn_{X,E_\cT})^{\cla}$.

\ssec{}    \label{abstract Koszul}

Consider now the direct sum
$$\on{Kosz}(E_\cT)^{\bullet,\cla,*}:=\underset{\cla_1+\cla_2=\cla}\bigoplus\, 
\fU(\cn_{X,E_\cT})^{\bullet,\cla_1}\star \Upsilon(\cn_{X,E_\cT})^{\cla_2}.$$

We can view it as a "module" over the $\fU(\cn_{X,E_\cT})^{\bullet,\cla}$'s via
\eqref{mult dot U}:
\begin{equation} \label{Koszul module}
\fU(\cn_{X,E_\cT})^{\bullet,\cmu}\star \on{Kosz}(E_\cT)^{\bullet,\cla,*}\to
\on{Kosz}(E_\cT)^{\bullet,\cla+\cmu,*}.
\end{equation}
It is also a "co-module" over the $\Upsilon(\cn_{X,E_\cT})^{\cla}$'s
via the maps \eqref{comult upsilon}:
\begin{equation} \label{Koszul comodule}
\on{Kosz}(E_\cT)^{\bullet,\cla,*}\to \underset{\cla=\cla_1+\cla_2}\bigoplus\,
\on{Kosz}(E_\cT)^{\bullet,\cla_1,*}\star \Upsilon(\cn_{X,E_\cT})^{\cla_2}.
\end{equation}

Now, by the very construction of the complexes $\fU(\cn_{X,E_\cT})^{\bullet,\cla}$, 
we can endow $\on{Kosz}(E_\cT)^{\bullet,\cla,*}$ with a differential, which makes
it into an acyclic complex for all $\cla\neq 0$. This differential is compatible
with the maps \eqref{Koszul module} and \eqref{Koszul comodule}.

\medskip

Let $\fU(\cn_{X,E_\cT})^{\bullet,-\cla,*}$ be the complex obtained from
$\fU(\cn_{X,E_\cT})^{\bullet,\cla}$ by applying Verdier duality term-wise.
Similarly, set 
$$\on{Kosz}(E_\cT)^{\bullet,-\cla}:=\underset{\cla_1+\cla_2=\cla}\bigoplus\, 
\fU(\cn_{X,E_\cT})^{\bullet,-\cla_1,*}\star \Omega(\cn_{X,E_\cT})^{-\cla_2}.$$
This is a complex possessing the structures dual to those of \eqref{Koszul module}
and \eqref{Koszul comodule}.

\medskip

Let us assume now that $E_\cT$ is regular. Note that by \corref{cohomology of ups}, 
the terms of $\fU(\cn_{X,E_\cT})^{m,\cla}$ are such that 
$H^j(X^\cla, \fU(\cn_{X,E_\cT})^{m,\cla})=0$ unless $j=0$. 
Applying the functor $H^0(X^\cla,?)$ term-wise to $\fU(\cn_{X,E_\cT})^{\bullet,\cla}$,
we obtain a $\BZ^{\geq 0}$-graded 
vector space that we shall denote by $\fu_{E_\cT}^{\bullet,\cla}$. The direct
sum 
$$\fu_{E_\cT}^{\bullet}:=\underset{\cla\in \cLambda^{pos}} \bigoplus\,
\fu_{E_\cT}^{\bullet,\cla}$$
is a $\BZ$-graded associative algebra. Dually, we set
$$\fu_{E_\cT}^{\bullet,-\cla,*}=H(X^\cla, \fU(\cn_{X,E_\cT})^{\bullet,-\cla,*}) \text{ and }
\fu_{E_\cT}^{\bullet,*}:=\underset{\cla\in \cLambda^{pos}} \bigoplus\,
\fu_{E_\cT}^{\bullet,-\cla,*},$$
the latter being a $\BZ^{\leq 0}$-graded co-associative co-algebra.

In addition, for $\cla\in \cLambda^{pos}$, we define the complex
$\sK(E_\cT)^{\bullet,-\cla}$ by applying $H(X^\cla,?)$ term-wise to
$\on{Kosz}(E_\cT)^{\bullet,-\cla}$. In other words,
$$\sK(E_\cT)^{\bullet,-\cla}=\underset{\cla_1+\cla_2=\cla}\bigoplus\,
\fu_{E_\cT}^{\bullet,-\cla_1,*}\otimes R_{E_\cT}^{-\cla_2}.$$

Set also $\sK(E_\cT)^\bullet=\underset{\cla\in \cLambda}\bigoplus\,
\sK(E_\cT)^{\bullet,-\cla}$. This is a DG module over $R_{E_\cT}$.
By the acyclicity of $\on{Kosz}(E_\cT)^{\bullet,-\cla}$, 
this DG module is quasi-isomorphic to $\BC$. Moreover, by construction,
when we disregard the differential, it is free and isomorphic 
to $\fu^{\bullet,*}_{E_\cT}\otimes R_{E_\cT}$.

\medskip

Let $\CM$ is a module over $R_{E_\cT}\simeq \CO_{\DefbE}$. We obtain that 
$\BC\overset{L}{\underset{\CO_{\DefbE}}\otimes}\CM$ can be computed by means of the
complex
\begin{equation} \label{Tor comp}
\fu^{\bullet,*}_{E_\cT}\otimes \CM\simeq \sK(E_\cT)^\bullet\underset{R_{E_\cT}}\otimes \CM,
\end{equation}
where the differential is obtained as a composition
$$\fu^{\bullet,*}_{E_\cT}\otimes \CM \to \sK(E_\cT)^\bullet\otimes \CM\simeq
\fu^{\bullet,*}_{E_\cT}\otimes R_{E_\cT}\otimes \CM\to
\fu^{\bullet,*}_{E_\cT}\otimes \CM,$$
where the first arrow is given by the differential on $\sK(E_\cT)^\bullet$, and the last 
arrow is given by the action of $R_{E_\cT}$ on $\CM$.

\ssec{}   \label{Koszul upstairs}

Consider now the direct sum
$$\on{Kosz}^{\bullet}_{\BunBb^\cmu}:=\underset{\cla'\in \cLambda^{pos}}\bigoplus\,
\fU(\cn_{X})^{\bullet,-\cla',*}\star \jmath_!(\IC_{\Bun^{\cmu+\cla'}})$$
as a graded perverse sheaf on $\BunBb^\cmu$. 

It acquires a differential by means of
\begin{multline*}
\underset{\cla'\in \cLambda^{pos}}\bigoplus\,
\fU(\cn_{X})^{\bullet,-\cla',*}\star \jmath_!(\IC_{\Bun^{\cmu+\cla'}}) \to
\underset{\cla'=\cla'_1+\cla'_2}\bigoplus\,
\fU(\cn_{X})^{\bullet,-\cla'_1,*}\star \Omega(\cn_{X})^{-\cla'_2}\star \jmath_!(\IC_{\Bun^{\cmu+\cla'}})\to \\
\to \underset{\cla'_1\in \cLambda^{pos}}\bigoplus\,
\fU(\cn_{X})^{\bullet,-\cla'_1,*}\star \jmath_!(\IC_{\Bun^{\cmu+\cla'_1}}),
\end{multline*}
where the first arrow is given by the differential on $\on{Kosz}(E_\cT)^{\bullet,-\cla'}$
(for $E_\cT$ being trivial) and the second arrow is given by \eqref{action upstairs}.

\begin{thm}  \label{exactness of Koszul upstairs}
The map $\jmath_!(\IC_{\Bun^\cmu_B})\to \IC_{\BunBb^\cmu}$ defines
a quasi-isomorphism 
$$\on{Kosz}^\bullet_{\BunBb^\cmu}\to \IC_{\BunBb^\cmu}.$$
\end{thm}

Before proving this theorem, let us show how it implies the 
assertion of point (2) of \thmref{main}. 

\medskip

We have to show that
\begin{equation}  \label{Eis tor}
\BC\overset{L}{\underset{\CO_{\DefbE}}\otimes}\Eis_!(E_\cT)\simeq
\Eisb(E_\cT).
\end{equation}

Let us tensor the complex $\on{Kosz}^\bullet_{\BunBb^\cmu}$ by
$(\fqb^\cmu)^*(\CS(E_\cT))$. We obtain a complex
$$\on{Kosz}_\cmu(\Eis_!(E_\cT))^\bullet:=
\underset{\cla'\in \cLambda^{pos}}\bigoplus\,
\fU(\cn_{X,E_\cT})^{\bullet,-\cla',*}\star 
\jmath_!\left(\IC_{\Bun^{\cmu+\cla'}}\otimes (\fq^{\cmu+\cla'})^*(\CS(E_\cT))\right).$$
 
Applying the functor $\fpb^\cmu_!$ to it term-wise, we obtain a complex
$$\sK_\cmu(\Eis_!(E_\cT))^\bullet:=
\underset{\cla'\in \cLambda^{pos}}\bigoplus\, \fu_{E_\cT}^{\bullet,-\cla',*}
\otimes \Eis_!^{\cmu+\cla'}(E_\cT).$$
The differential on $\sK_\cmu(\Eis_!(E_\cT))^\bullet$ is given by the formula for
the differential on \eqref{Tor comp}, using the action maps \eqref{action maps}.
Hence, by \secref{abstract Koszul}, the direct sum
\begin{equation} \label{dir sum Kosz}
\sK(\Eis_!(E_\cT))^\bullet:=\underset{\cmu}\bigoplus\, \sK_\cmu(\Eis_!(E_\cT))^\bullet\simeq
\fu_{E_\cT}^{\bullet,*}\otimes \Eis_!(E_\cT)
\end{equation}
is quasi-isomorphic to
the LHS of \eqref{Eis tor}. 

Now, by \thmref{exactness of Koszul upstairs}, the complex $\on{Kosz}_\cmu(\Eis_!(E_\cT))^\bullet$
is quasi-isomorphic to the perverse sheaf
$\IC_{\BunBb^\cmu}\otimes (\fqb^\cmu)^*(\CS(E_\cT))$.
We obtain that $\sK_\cmu(\Eis_!(E_\cT))^\bullet$ is quasi-isomorphic 
to $\Eisb_\cmu(E_\cT)$. 

Therefore, $\sK(\Eis_!(E_\cT))^\bullet$ is quasi-isomorphic 
to the RHS of \eqref{Eis tor}, as required.
 
\ssec{Proof of \thmref{exactness of Koszul upstairs}} 
 
We proceed by induction on $|\cla|$ by showing that the map in question is
a quasi-isomorphism over the open substack $\BunBb^{\cmu,\leq\cla}$.
The base of the induction, i.e., $\cla=0$ evidently holds. Thus, we assume
that the quasi-isomorphism in question is valid for all $\cla'$ with
$|\cla'|<|\cla|$.

Thus, it is sufficient to show that the map
\begin{equation} \label{Koszul on stratum}
\imath_\cla^*\left(\on{Kosz}^\bullet_{\BunBb^\cmu}\right)\to
\imath_\cla^*\left(\IC_{\BunBb^\cmu}\right)
\end{equation}
is a quasi-isomorphism.

Note now that by \lemref{dot and without}, the LHS of
\eqref{Koszul on stratum} is quasi-isomorphic to 
\begin{equation} \label{stalk on stratum}
\fU(\cn_X)^{-\cla,*}\boxtimes \IC_{\Bun_B^{\cmu+\cla}}.
\end{equation}
By \propref{IC calculation}, the RHS of \eqref{Koszul on stratum} 
is also isomorphic to the expression in \eqref{stalk on stratum}.
Thus, to prove the theorem, we need to show that the
resulting endomorphism of \eqref{stalk on stratum} is
an isomorphism. \footnote{We do not claim at this
stage that this is the identity automorphism.}

\medskip

Consider the canonical filtration on $\fU(\cn_X)^{\cla,*}$, corresponding
to the perverse t-structure. The associated graded is isomorphic to
$$\underset{\cla=m_1\cdot \check\beta_1+...+m_k\cdot \check\beta_k}
\bigoplus\,
((\cn^*)^{-\check\beta_1}_X)^{(n_1)}\star...\star 
((\cn^*)^{-\check\beta_k}_X)^{(n_k)}[2(n_1+...+n_k)],$$
where $\check\beta_1,...,\check\beta_k$ are {\it not necessarily simple}
roots of $\cg$. Each of the summands, except the one corresponding to $k=1$
(in which case $\cla$ is itself a root), is a cohomologically shifted
perverse sheaf, which is the intermediate restriction of its extension to
$X^\cla-\Delta(X)$.

However, we claim that by the induction hypothesis we can assume that 
\eqref{Koszul on stratum} is an isomorphism over
$(X^\cla-\Delta(X))\times \Bun_B^{\cmu+\cla}$. Indeed, 
since the assertion is local, we can replace $\BunBb^{\cmu,\leq\cla}$
by the Zastava space, and then apply the factorization principle,
\eqref{factorization of Zastava}. 

\medskip

Hence, the map \eqref{Koszul on stratum} induces an isomorphism
on the associated graded pieces of \eqref{stalk on stratum}, except,
possibly, on $(\cn^*)^{-\cla}_X[2]\boxtimes \IC_{\Bun_B^{\cmu+\cla}}$,
where $(\cn^*)^{-\cla}_X[2]$ is the last quotient of $\fU(\cn_X)^{\cla,*}$
(and which can only occur if $\cla$ is a root). 

\medskip

Suppose, by contradiction, that this map was not an isomorphism, i.e., equal to zero. 
We would obtain that 
$$\on{Cone}\left(\on{Kosz}^\bullet_{\BunBb^\cmu}\to \IC_{\BunBb^\cmu}\right)\simeq
\on{Cone}\left((\imath_\cla)_!((\cn^*)^{-\cla}_X[2])\to (\imath_\cla)_!((\cn^*)^{-\cla}_X[2])\right),$$
and therefore
$$h^0\biggl(\on{Cone}\left(\on{Kosz}^\bullet_{\BunBb^\cmu}\to \IC_{\BunBb^\cmu}\right)\biggr)
|_{\BunBb^{\cmu,\leq \cla}}\neq 0.$$
But this is a contradiction, since $\IC_{\BunBb^\cmu}$ is an irreducible perverse sheaf,
and the complex $\on{Kosz}^\bullet_{\BunBb^\cmu}$ is concentrated in non-positive
perverse cohomolgical degrees.

\ssec{A comparison of two isomorphisms}   \label{comparison}

Let us denote by $\on{Kosz}^{\bullet,*}_{\BunBb^\cmu}$ the complex
obtained from $\on{Kosz}^\bullet_{\BunBb^\cmu}$ by applying term-wise
Verdier duality. Its terms are given by
$$\underset{\cla\in \cLambda^{pos}}\bigoplus\,
\fU(\cn_{X})^{\bullet,\cla}\star \jmath_*(\IC_{\Bun^{\cmu+\cla}}),$$
and the differential by that on $\on{Kosz}(E_\cT)^{\bullet,\cla,*}$ and maps
dual to those of \eqref{action upstairs}. 

As this complex is quasi-isomorphic to $\IC_{\BunBb^\cmu}$, it 
can be used to calculate $\imath^!_\cla(\IC_{\BunBb^\cmu})$. From
\lemref{dot and without}, we obtain a quasi-isomorphism
\begin{equation} \label{Koszul calc}
\imath^!_\cla(\IC_{\BunBb^\cmu})\simeq \fU(\cn_X)^{\cla}\boxtimes \IC_{\Bun_B^{\cmu+\cla}}.
\end{equation}

Note that \propref{IC calculation} and \thmref{FFKM thm}(1) give two more
isomorphisms between the same objects.

\begin{conj}
The isomorphisms of \propref{IC calculation} and \thmref{FFKM thm}(1)
coincide.
\end{conj}

In the remainder of this section we will prove that the
isomorphisms of \eqref{Koszul calc} and \thmref{FFKM thm}(1) coincide.

\ssec{}   
The starting point is the following observation:

\begin{lem}  \label{U rigid}
Suppose that we have a system of automorphisms $\phi^\cla$ of the sheaves
$\fU(\cn_X)^{\cla}$ with the following two properties:

\begin{itemize}

\item $\phi^\cla=\on{id}$ for $\cla$ being a simple root $\check\alpha_i$.

\item The system $\{\phi^\cla\}$ is compatible with the maps \eqref{mult U}.

\end{itemize}

Then $\phi^\cla=\on{id}$ for all $\cla$.
\end{lem}

The lemma follows from the fact that the simple root spaces generate
$U(\cn)$. We apply it in our situation for $\phi^\cla$ being the discrepancy
of the maps \eqref{Koszul calc} and \thmref{FFKM thm}(1). The fact that
$\phi^{\check\alpha_i}=\on{id}$ follows from the construction of both
maps. Thus, it remains to check the compatibility with the product
operation \eqref{mult U}.

\medskip

First, let us observe that the map 
\begin{equation} \label{expl adjunction}
\fU(\cn_X)^\cla\star \jmath_!(\IC_{\Bun_B^{\cmu+\cla}})\to \IC_{\BunBb^\cmu},
\end{equation}
corresponding by adjunction to \eqref{Koszul calc}, is represented by
$$\fU(\cn_X)^{\bullet,\cla}\star \jmath_!(\IC_{\Bun_B^{\cmu+\cla}})\to
\fU(\cn_X)^{\bullet,\cla}\star \IC_{\BunBb^{\cmu+\cla}}\hookrightarrow 
\on{Kosz}^{\bullet,*}_{\BunBb^\cmu},$$
where the image of the last arrow belongs to the kernel of the differential,
because $\IC_{\BunBb^{\cmu'}}\subset \jmath_*(\IC_{\Bun_B^{\cmu'}})$
lies in the kernel of the maps dual to those of \eqref{action upstairs}.
By construction, the map \eqref{expl adjunction} extends to a map
\begin{equation} \label{IC map}
\fU(\cn_X)^{\cla}\star \IC_{\BunBb^{\cmu+\cla}}\simeq
\fU(\cn_X)^{\bullet,\cla}\star \IC_{\BunBb^{\cmu+\cla}}\to \on{Kosz}^{\bullet,*}_{\BunBb^\cmu}\simeq
\IC_{\BunBb^\cmu},
\end{equation}
whose existence is asserted in \thmref{FFKM thm}(2).

Using \thmref{FFKM thm}(3), it remains to establish the commutativity
of the following diagram in $D^b(\BunBb^\cmu)$, whose arrows are induced by \eqref{IC map}:
$$
\CD
\fU(\cn_X)^{\cla_1}\star \fU(\cn_X)^{\cla_2}\star \IC_{\BunBb^{\cmu+\cla_1+\cla_2}}  @>>>
\fU(\cn_X)^\cla\star \IC_{\BunBb^{\cmu+\cla}}  \\
@VVV   @VVV   \\
\fU(\cn_X)^{\cla_1}\star \IC_{\BunBb^{\cmu+\cla_1}}  @>>> \IC_{\BunBb^\cmu}.
\endCD
$$

Consider the map
\begin{equation} \label{mult Koszul up}
\fU(\cn_X)^{\bullet,\cla}\star \on{Kosz}^{\bullet,*}_{\BunBb^{\cmu+\cla}}\to
\on{Kosz}^{\bullet,*}_{\BunBb^{\cmu}},
\end{equation}
induced by \eqref{mult dot U}. The compatibility of the latter with the
differential on the $\fU(\cn_X)^{\bullet,\cla}$'s implies that
\eqref{mult Koszul up} is a map of complexes.

The commutativity of the above diagram follows now from the next statement:

\begin{lem}    \label{expl Fink action}
The diagram
$$
\CD
\fU(\cn_X)^{\cla}\star \IC_{\BunBb^{\cmu+\cla}} @>>> \IC_{\BunBb^\cmu} \\
@V{\sim}VV   @V{\sim}VV \\
\fU(\cn_X)^{\bullet,\cla}\star \on{Kosz}^{\bullet,*}_{\BunBb^{\cmu+\cla}} @>>>
\on{Kosz}^{\bullet,*}_{\BunBb^{\cmu}}
\endCD
$$
is commutative in the derived category.
\end{lem}

\begin{proof}
This follows from the fact that the composition
$$\fU(\cn_X)^{\bullet,\cla}\star \IC_{\BunBb^{\cmu+\cla}}\to 
\fU(\cn_X)^{\bullet,\cla}\star \on{Kosz}^{\bullet,*}_{\BunBb^{\cmu+\cla}} 
\overset{\text{\eqref{mult Koszul up}}}\longrightarrow
\on{Kosz}^{\bullet,*}_{\BunBb^{\cmu}}$$
equals the map of \eqref{IC map}.
\end{proof}
 
\section{The Hecke property in the case of $GL_2$}   \label{Hecke prop for GL(2)}

In this section we will prove \thmref{Hecke} for $GL_2$ by a direct calculation
for $V\in \Rep(\cG)$ being the standard $2$-dimensional
representation of $\cG=GL_2$. We retain the notation of \secref{case of GL(2)}. 

\ssec{}

Let $E_{\cB,\wh{univ}}$ denote the tautological $\wh\CO_{\DefbE}$-family of 
2-dimensional local systems on $X$. Let $E_{\cB,\wh{univ},x}$ be its fiber at $x$; this
is a locally free $\wh\CO_{\DefbE}$-module of rank $2$. Let 
$E_{\cB,univ,x}$ be the corresponding $\CO_{\DefbE}\simeq \Sym(\sW)$-module, 
i.e., the direct sum of homogeneous components of $E_{\cB,\wh{univ},x}$ 
with respect to the natural $\cT$-action. We have a short exact sequence
$$0\to \Sym(\sW)\otimes E_{1,x}\to E_{\cB,univ,x}\to \Sym(\sW)\otimes E_{2,x}\to 0.$$

We need to prove that for any $V\in \Rep(\cG)$ being the standard 
$2$-dimensional representation of $GL_2$ there exists a canonical 
isomorphism of $\cLambda$-graded perverse sheaves on $\Bun_G$:
\begin{equation} \label{local Hecke for GL(2)}
\on{H}^{\on{basic}}_x(\Eis_!(E_\cT))\simeq E_{\cB,univ,x}
\underset{\Sym(\sW)}\otimes \Eis_!(E_\cT),
\end{equation}
where $\on{H}^{\on{basic}}_x$ is the Hecke functor corresponding to the
standard $2$-dimensional representation of $GL_2$.

\ssec{}

Let us first describe the fiber $E_{\cB,univ,x}$ explicitly.

\medskip

Set $\sW':=H^1(X-x,E_2\otimes E_1^{-1})$.
The cokernel $\sW'/\sW$ is canonically isomorphic to $E_{2,x}\otimes E_{1,x}^{-1}$. 
Consider the map
\begin{equation} \label{quotient}
\on{Sym}(\sW)\otimes \sW\to \left(\on{Sym}(\sW)\otimes \sW'\right)\bigoplus \on{Sym}(\sW),
\end{equation}
where the first component corresponds to the embedding of $\sW$ into $\sW'$, and the
second component is given by the multiplication map.

\begin{lem} \label{descr of extenstion}
The $\on{Sym}(\sW)$-module $E_{\cB,univ,x}$
is canonically isomorphic cokernel of the map of \eqref{quotient},
tensored by $E_{1,x}$.
\end{lem}

\begin{cor}
Let $\CF$ be an object of some abelian category endowed with an action of
$\on{Sym}(\sW)$. Then
$$\CF\overset{L}{\underset{\on{Sym}(\sW)}\otimes} E_{\cB,univ,x}$$ is 
canonically  quasi-isomorphic
to the complex 
$$\sW\otimes \CF\to \left(\sW'\otimes \CF\right)\bigoplus \CF,$$ tensored by
$E_{1,x}$.
\end{cor}

Thus, we obtain that the existence of the isomorphism \eqref{local Hecke for GL(2)}
is equivalent to the following assertion:

\begin{prop}   \label{verify Hecke for GL(2)}
The object $\on{H}^{\on{basic}}_x(\Eis_!^{d_1,d_2}(E_\cT))\otimes E^{-1}_{1,x}$ 
is canonically quasi-isomorphic
to the complex
$$H^1(X,E_2\otimes E_1^{-1})\otimes \Eis_!^{d_1,d_2-1}(E_\cT)\to
H^1(X-x,E_2\otimes E_1^{-1})\otimes \Eis_!^{d_1,d_2-1}(E_\cT)
\bigoplus \Eis_!^{d_1-1,d_2}(E_\cT).$$
\end{prop}

\ssec{Proof of \propref{verify Hecke for GL(2)}}   \label{expl for sl(2)}

Let $\H^{\on{basic}}_x$ denote the stack classifying triples
$$(\CM,\CM',\beta:\CM\hookrightarrow \CM'),$$ where $(\CM,\CM')\in \Bun_G=\Bun_2$,
and $\beta$ is an embedding such that the quotient $\CM'/\CM$ has length $1$
and is supported at $x$. Let $\hl,\hr$ denote the two projections from $\H^{\on{basic}}_x$
to $\Bun_G$ that remember $\CM$ and $\CM'$, respectively.

Consider the Cartesian product
$$\H_x'{}^{\on{basic}}:=\H^{\on{basic}}_x\underset{\Bun_G}\times \BunBb^{d_1,d_2},$$
where $\H^{\on{basic}}_x$ maps to $\Bun_G$ via $\hr$. This stack classifies quintuples
$$(\CL'\overset{\kappa'}\hookrightarrow \CM'; \,\,\CM \overset{\beta}\hookrightarrow \CM'),$$
where $(\CM,\CM',\beta)$ are as above, $\CL'$ is a line bundle on $X$, and $\kappa'$
is an embedding of $\CL'$ into $\CM'$ as a coherent sub-sheaf. 

Let us denote by $\hr'$ the natural projection $\H_x'{}^{\on{basic}}\to \BunBb^{d_1,d_2}$
that remembers the data of $\CL'\overset{\kappa'}\hookrightarrow \CM'$. Let
$\hl'$ denote the map $\H_x'{}^{\on{basic}}\to \BunBb^{d_1,d_2+1}$ that sends a quintuple
as above to $(\CL,\kappa,\CM)$, where $\CL:=\CL'(-x)$ and $\kappa$ is the
(unique and well-defined) embedding of $\CL$ into $\CM$, such that $\beta\circ \kappa$
equals
$$\CL\hookrightarrow \CL'\overset{\kappa'}\to \CM'.$$

\medskip

By construction and base change,
$$\on{H}_x(\Eis_!^{d_1,d_2}(E_\cT))\simeq
\fpb^{d_1-1,d_2}_!\circ \hr'{}_! \circ \hl'{}^*\left(
(\imath_0^{d_1,d_2})_!(\IC_{\Bun_B^{d_1,d_2}})\otimes (\fqb^{d_1,d_2})^*(\CS(E_\cT))\right),$$
which, in turn, is isomorphic to
$$\fpb^{d_1-1,d_2}_!\left(\hr'{}_! \circ \hl'{}^*
\left((\imath_0^{d_1,d_2})_!(\IC_{\Bun_B^{d_1,d_2}})\right)
\otimes (\fqb^{d_1-1,d_2})^*(\CS(E_\cT))\right)\otimes E_{1,x}[1].$$

Thus, to prove \propref{verify Hecke for GL(2)}, it is sufficient to show that
\begin{multline} \label{sl(2) Hecke map}
\hr'{}_! \circ \hl'{}^* \left((\imath_0^{d_1,d_2})_!(\IC_{\Bun_B^{d_1,d_2}})\right)
\simeq \\
\simeq \on{Co-Ker}\Biggl(
(\imathb_1^{d_1,d_2-1})_!\left(\IC_X\boxtimes 
(\mathi_0^{d_1,d_2-1})_!(\IC_{\Bun^{d_1,d_2-1}_B})\right)\to \\
\to (\imathb_1^{d_1,d_2-1})_!\left(j_x{}_*(\IC_{X-x})\boxtimes 
(\mathi_0^{d_1,d_2-1})_!(\IC_{\Bun^{d_1-1,d_2}_B})\right)\bigoplus
(\imath_0^{d_1-1,d_2})_!\left(\IC_{\Bun_B^{d_1-1,d_2}}\right)\Biggr),
\end{multline}
where $j_x$ denotes the open embedding $X-x\hookrightarrow X$.

\medskip

To establish the required isomorphism, note that both the LHS and the RHS
are extensions by $0$ from the open substack 
$$\imath_0^{d_1-1,d_2}(\Bun_B^{d_1-1,d_2})\cup 
\imath_1^{d_1,d_2-1}(\Bun_B^{d_1,d_2-1})\subset \BunBb^{d_1-1,d_2}.$$

Over this open subset, $\imath_1^{d_1,d_2-1}(X\times \Bun_B^{d_1,d_2-1})$ is a
smooth divisor, which itself contains the divisor,
corresponding to the point $x\in X$. The map $\hl'$ is an isomorphism away from
$x\times \Bun_B^{d_1,d_2-1}$, and over this codimension-$2$ closed substack
it is a fibration with typical fiber $\BP^1$.

\medskip

Therefore, our situation admits the following local model. Let ${\mathsf f}:\wt{\BA}{}^2\to \BA^2$ be the 
blow-up of the affine plane at the origin. Let $\imath^1$ be the embedding of a fixed line
$\BA^1\hookrightarrow \BA^2$, and let $\imath^0$ be the embedding of its
complement; we will denote by $\wt{\imath}{}^0$ the embedding of
the complement of the proper transform of $\BA^1$ into $\wt{\BA}{}^2$. 
Finally, let $j$ denote the embedding $\BA^1-0\hookrightarrow \BA^1$. We have:

\begin{lem}
$${\mathsf f}_!\left(\wt{\imath}{}^0_!(\IC_{\wt{\BA}{}^2-\BA^1})\right)\simeq
\on{Co-Ker}\Biggl(\imath^1_!(\IC_{\BA^1})\to \imath^1_!(j_*(\IC_{\BA^1-0}))\bigoplus 
\imath^0_!(\IC_{\BA^2-\BA^1})\Biggr).$$
\end{lem}

The proof is a straightforward verification. 

\section{The Hecke property}    \label{Hecke section}

\ssec{}

Let $E_{\cB,\wh{univ}}$ the canonical $\cB$-local system over $X$
over the formal scheme $\DefbE$. For a point $x\in X$ and
$V\in \Rep(\cG)$, let $V_{E_{\cB,\wh{univ},x}}$ be the fiber at $x$
of the local system associated with $E_{\cB,\wh{univ}}$ and $V$.
This is a locally free $\wh\CO_{\DefbE}$-module of rank equal
to $\dim(V)$.

As in the case of $GL_2$, the first step is to describe 
$V_{E_{\cB,\wh{univ},x}}$ explicitly as an $\wh\CO_{\DefbE}$-module,
in terms of the isomorphism of \thmref{deformation base}.

\medskip

Let $\ceta\in \cLambda$ be a large enough weight, so that
$\ceta-\cnu\in \cLambda^{pos}$ whenever $V(\cnu)\neq 0$.
(If $G$ is not of the adjoint type, we will assume that $Z(G)$
acts on $V$ by a single character.)
For $\cla\in \cLambda^{pos}$ we consider the following
complex on $X^\cla$:

We consider the stratification of $X^\cla$, numbered by triples:
$(\cla_1,\cla_2\, |\, \cla_1+\cla_2=\cla, \fP(\cla_1))$, which each
stratum corresponds to the configuration
$$\cla^k_1\cdot x_k, \cla_2\cdot x,\, x_{k}\neq x_{k'}\neq x, \Sigma\, \cla^k_1=\cla_1.$$

On each such stratum we put the locally-constant complex, denoted
$\on{Cous}(\cn^*_{E,\cT},V)^{\cla_1,\cla_2,\fP(\cla_1)}$, whose $!$-stalk
at the above point is
$$\underset{k}\otimes\, \left(\cLambda^\bullet(\cn^*_{E_{\cT,x_k}})\right)^{-\cla^k_1}\bigotimes
\left(\cLambda^\bullet(\cn^*_{E_{\cT,x}})\otimes V_{E_{\cT,x}}\right)^{\ceta-\cla_2},$$
with the standard Chevalley differential. 

Let $j^{\fP(\cla_1),\cla_2}$ denote the embedding of the corresponding stratum 
into $X^\cla$. The direct sum of complexes
$$\underset{\cla_1,\cla_2,\fP(\cla_1)}\bigoplus\,
j^{\fP(\cla_1),\cla_2}_*\left(\on{Cous}(\cn^*_{E,\cT},V)^{\cla_1,\cla_2,\fP(\cla_1)}\right)$$
acquires a natural differential. We shall denote the complex, associated to the resulting
bi-complex by $\Omega(\cn_{X,E_{\cT}},V_{E_{\cT,x}})^{\ceta-\cla}$. When $E_\cT$
is trivial (i.e., when there is no twisting), we shall denote this complex simply by
$\Omega(\cn_X,V_{x})^{\ceta-\cla}$.

\begin{prop}  \label{Upsilon V}
The complex $\Omega(\cn_{X,E_{\cT}},V_{E_{\cT,x}})^{\ceta-\cla}$
is a perverse sheaf. If $E_\cT$ is regular, then the cohomology 
$H(X^\cla,\Omega(\cn_{X,E_{\cT}},V_{E_{\cT,x}})^{\ceta-\cla})$ is concentrated
in degree $0$.
\end{prop}

\begin{proof}

When we regard $V$ as a $\cB$-module, it carries a canonical
filtration, parametrized by the partially order set $\cLambda$
(with the order relation $\cla'\geq \cla''$ if $\cla'-\cla''\in \cLambda^{pos}$),
such that $\on{gr}^{\cnu}(V)\simeq V(\cnu)$ (here and in the sequel $V(\cnu)$
denotes the $\cnu$ weight space of $V$).

This filtration induces a filtration on 
$\Omega(\cn_{X,E_{\cT}},V_{E_{\cT,x}})^{\ceta-\cla}$,
such that
$$\on{gr}^{\cnu}\left(\Omega(\cn_{X,E_{\cT}},V_{E_{\cT,x}})^{\ceta-\cla}\right)\simeq
\Omega(\cn_{X,E_{\cT}})^{\ceta-\cla-\cnu}\otimes V(\cnu),$$
where $\Omega(\cn_{X,E_{\cT}})^{\ceta-\cla-\cnu}$, which is by definition a
perverse sheaf on $X^{\cla-(\ceta-\cnu)}$, is viewed as a perverse sheaf on
$X^\cla$ via the map 
$$X^{\cla-(\ceta-\cnu)}\to X^{\cla-(\ceta-\cnu)}\times X^{\ceta-\cnu}\to X^\cla,$$
where the first arrow corresponds to the point 
$(\ceta-\cnu)\cdot x\in X^{\ceta-\cnu}$.

This proves both points of the proposition in view of \propref{Upsilon}
and \corref{cohomology of ups}.

\end{proof}

\ssec{}

Let us make the following observation. Let us replace the element $\ceta$
by another element $\ceta'$; with no restriction of generality we can assume
that $\ceta'-\ceta\in \cLambda^{pos}$. Let $\cnu=\ceta-\cla=\ceta'-\cla'$ with
$\cla,\cla'\in \cLambda^{pos}$. Then we have the perverse sheaf
$\Omega(\cn_{X,E_{\cT}},V_{E_{\cT,x}})^{\ceta-\cla}$ on $X^\cla$ and
the perverse sheaf $\Omega(\cn_{X,E_{\cT}},V_{E_{\cT,x}})^{\ceta'-\cla'}$
on $X^{\cla'}$. However, it is easy to see that the latter is canonically
isomorphic to the direct image of the former under the closed embedding 
$$X^\cla\hookrightarrow X^{\cla'},$$
corresponding to adding the coloured divisor $(\cla'-\cla)\cdot x$.

For $\cnu\in \cLambda$, let ${}_{\infty\cdot x}X^\cnu$ denote 
the ind-scheme $\underset{\cla\in \cLambda^{pos}}{\underset{\longrightarrow}{lim}}\,
X^\cla$, which we think of as classifying divisors of the form
$\cla'\cdot x-\Sigma\, \cla_k\cdot x_k$, where $\cla_k\in \cLambda^{pos}$ for 
$x_k\neq x$, and $\cla'\in \cLambda$ arbitrary, but so that 
$\cla'-\Sigma\, \cla_k=\cnu$. 

We obtain that for $\cnu\in \cLambda$ we have a well-defined perverse
sheaf $\Omega(\cn_{X,E_{\cT}},V_{E_{\cT,x}})^{\cnu}$ on ${}_{\infty\cdot x}X^\cnu$,
which equals $\Omega(\cn_{X,E_{\cT}},V_{E_{\cT,x}})^{\ceta-\cla}$ on $X^\cla$
for $\ceta$ and $\cla$ large enough with $\ceta-\cla=\cnu$. 
\footnote{This construction has the advantage that it makes sense whether or not 
$Z(G)$ acts on $V$ by a single character.} Letting $V$ be the trivial representation,
we recover $\Omega(\cn_{X,E_{\cT}})^{-\cla}$ as a perverse sheaf on
$X^\cla\subset {}_{\infty\cdot x}X^{-\cla}$.

\medskip

For each $\cla\in \cLambda^{pos}$
We have natural addition maps
$${}_{\infty\cdot x}X^{\cnu_1}\times {}_{\infty\cdot x}X^{\cnu_2}\to {}_{\infty\cdot x}X^{\cnu_1+\cnu_2},$$
and the corresponding functors
$$\star:D^b({}_{\infty\cdot x}X^{\cnu_1})\times 
D^b({}_{\infty\cdot x}X^{\cnu_2})\to D^b({}_{\infty\cdot x}X^{\cnu_1+\cnu_2}),.$$

By construction, there exists a canonical map
\begin{equation} \label{tensor Omega V}
\Omega(\cn_{X,E_{\cT}},V^1_{E_{\cT,x}})^{\cnu_1}\star 
\Omega(\cn_{X,E_{\cT}},V^2_{E_{\cT,x}})^{\cnu_2}\to
\Omega(\cn_{X,E_{\cT}},(V^1\otimes V^2)_{E_{\cT,x}})^{\cnu_1+\cnu_2},
\end{equation}
which is associative in the natural sense. Letting $V_2$ be the trivial
representation, we obtain the map
\begin{equation} \label{mult Omega V}
\Omega(\cn_{X,E_{\cT}},V_{E_{\cT,x}})^{\cnu}\star 
\Omega(\cn_{X,E_{\cT}})^{-\cla}\to 
\Omega(\cn_{X,E_{\cT}},V_{E_{\cT,x}})^{\cnu-\cla}.
\end{equation}

In particular, assuming that $E_\cT$ is regular, set
$$R(V_x)_{E_\cT}^\cnu:=
H\left({}_{\infty\cdot x}X^\cnu,\Omega(\cn_{X,E_{\cT}},V_{E_{\cT,x}})^{\cnu}\right) \text{ and }
R(V_x)_{E_\cT}:=\underset{\cnu}\oplus\, R(V_x)_{E_\cT}^\cnu.$$
We obtain that $R(V_x)_{E_\cT}$ is a ($\cLambda$-graded) module over
the ($\cLambda$-graded) commutative algebra $R_{E_\cT}$.

This module is finitely generated and projective. Indeed, the filtration,
introduced in the proof of \propref{Upsilon V}, induces a filtration on the
above module, with the associated graded being the free module
on the vector space $V$.

Set also
$$\wh{R}(V_x)_{E_\cT}=\underset{\cnu}\Pi\, R(V_x)_{E_\cT}^\cnu\simeq
R(V_x)_{E_\cT}\underset{R_{E_\cT}}\otimes \wh{R}_{E_\cT}.$$

\begin{lem}  \label{fibers over deformation}
Under the isomorphism $\wh\CO_{\DefbE}\simeq \wh{R}_{E_\cT}$,
the module $V_{E_{\cB,\wh{univ},x}}$ corresponds to 
$\wh{R}(V_x)_{E_\cT}$.
\end{lem}
The proof will be given in \secref{proof of deformation base}.

\medskip

Let $V_{E_{\cB,univ,x}}$ be the $\cLambda$-graded version of
$V_{E_{\cB,\wh{univ},x}}$. We obtain:

\begin{cor}   \label{graded fibers over deformation}
Under the isomorphism $\CO_{\DefbE}\simeq R_{E_\cT}$, the
$\CO_{\DefbE}$-module $V_{E_{\cB,univ,x}}$ corresponds to
to the $R_{E_\cT}$-module $R(V_x)_{E_\cT}$.
\end{cor}

\ssec{}

Let ${}_{\infty\cdot x}\BunBb^\cmu$ denote the ind-version of $\BunBb^\cmu$,
where the maps $\kappa^\lambda$ are allowed to have poles of arbitrary
order at $x$, see \cite{BG}, Sect. 4.1.1. The stack $\BunBb^\cmu$ is a 
closed substack of 
${}_{\infty\cdot x}\BunBb^\cmu$; hence perverse sheaves (or objects of the
derived category) on the former can be thought of as corresponding objects
on the latter. We will denote by $_{\infty\cdot x}\BunBb$ the union
of the $_{\infty\cdot x}\BunBb^\cmu$'s over $\cmu\in \cLambda$. 

For $\cnu\in \cLambda$ we have a natural map
$$_{\infty\cdot x}\imathb_\cnu:\,{}_{\infty\cdot x}X^\cnu\times 
{}_{\infty\cdot x}\BunBb^{\cmu-\cnu}\to {}_{\infty\cdot x}\BunBb^\cmu,$$
defined in the same way as $\imathb_\cla$.
Let $_{\infty\cdot x}\imath_\cnu$ denote the restriction of $_{\infty\cdot x}\imathb_\cnu$
to the locally closed substack
$$_{\infty\cdot x}X^\cnu\times \Bun_B^{\cmu-\cnu}\subset{}
_{\infty\cdot x}X^\cnu\times {}_{\infty\cdot x}\BunBb^{\cmu-\cnu}.$$
The images of the maps $_{\infty\cdot x}\imath_\cnu$ for
$\cnu\in \cLambda$ define a stratification of $_{\infty\cdot x}\BunBb^\cmu$.

Using the maps $_{\infty\cdot x}\imathb_\cnu$ we define the functors
\begin{equation} \label{star conv}
\star:D^b({}_{\infty\cdot x}X^\cnu)\times D^b({}_{\infty\cdot x}\BunBb^{\cmu+\cnu})
\to D^b({}_{\infty\cdot x}\BunBb^{\cmu}).
\end{equation}

\medskip

Let $\H_x$ be the Hecke stack for $G$, and let $\H'_x$ be its version over
$_{\infty\cdot x}\BunBb$ (see \cite{BG}, 4.1.2), so that we have a 
commutative diagram with both squares Cartesian:
$$
\CD
_{\infty\cdot x}\BunBb @<{\hl{}'}<<  \H'_x  @>{\hr{}'}>> _{\infty\cdot x}\BunBb \\
@V{\fpb}VV   @V{\fpb'}VV    @V{\fpb}VV  \\
\Bun_G    @<{\hl}<< \H_x  @>{\hr}>> \Bun_G.
\endCD
$$

Thus, for each $V\in \Rep(\cG)$ we can associate the perverse sheaves 
$\CV$ and $\CV'$ on $\H_x$ and $\H'_x$, respectively, and the Hecke functors
$$\on{H}^V_x: D^b(\Bun_G)\to D^b(\Bun_G), \text{ given by }
\CT\mapsto \hl_!(\CV\otimes \hr{}^*(\CT))$$ and 
$$\on{H}'{}^V_x: D^b({}_{\infty\cdot x}\BunBb)\to D^b({}_{\infty\cdot x}\BunBb), \text{given by }
\CT'\mapsto \hl{}'_!(\CV'\otimes \hr{}'{}^*(\CT')).$$

Note that the convolution functors \eqref{star conv} introduced
above and the Hecke functors $\on{H}'{}^V_x$ naturally commute.

\ssec{}   \label{structure of Hecke}

The main geometric ingredient, from which we shall deduce \thmref{Hecke} is 
a description of the object 
\begin{equation} \label{Hecke of ext by zero}
\on{H}'{}^V_x\left(\jmath_!(\IC_{\Bun^\cmu_B})\right).
\end{equation}
Before stating the corresponding theorem, let us make several observations.

First, let us recall from \cite{BG}, Theorem 3.3.2,
that there is a canonical isomorphism:
\begin{equation} \label{IC Hecke}
\on{H}'{}^V_x(\IC_{\BunBb^\cmu})\simeq \underset{\cnu}\oplus\, 
V(\cnu)\otimes \delta^\cnu_{x}\star \IC_{\BunBb^{\cmu-\cnu}},
\end{equation}
where $\delta^\cnu_{x}$ denotes the sky-scraper at the point
$\cnu\cdot x\in {}_{\infty\cdot x}X^\cnu$. Here the weight spaces $V(\cnu)$
are realized as cohomologies of the corresponding spherical 
sheaves on the affine Grassmannian along semi-infinite orbits
(see \cite{BG}, proof of Theorem 3.3.2).

Combining this with \corref{Groth group comp}, we obtain that 
the object \eqref{Hecke of ext by zero} is a perverse sheaf.

\medskip

Moreover, the fiber product
$$\H'_x \underset{\Bun_G,\hr}\times \Bun_B^\cmu$$ is naturally
stratified, according to the order of zero/pole at $x$ of the
corresponding generalized $B$-structure under the projection $\hl'$,
see \cite{BG}, Sect. 3.3.4. Hence,
we obtain a filtration on \eqref{Hecke of ext by zero},
parametrized by $\cLambda$ with
$$\on{gr}^\cnu\left(\on{H}'{}^V_x\left(\jmath_!(\IC_{\Bun^\cmu_B})\right)\right)
\simeq V(\cnu)\otimes \delta^\cnu_{x}\star 
\jmath_!(\IC_{\Bun^{\cmu+\cnu}_B}).$$

We will prove:

\begin{thm} \label{Hecke of extension by zero}  \hfill 

\smallskip

\noindent{\em (A)}
We have an isomorphism of perverse sheaves:
\begin{equation} \label{need to construct V}
\Omega(\cn_X,V_x)^{\cnu}\boxtimes \IC_{\Bun^{\cmu-\cnu}_B}\simeq
h^0\biggl({}_{\infty\cdot x}\imath_\cnu^!
\left(\on{H}'{}^V_x\left(\jmath_!(\IC_{\Bun^\cmu_B})\right)\right)\biggr).
\end{equation}
In particular, by adjunction we obtain a canonical map of perverse sheaves 
\begin{equation} \label{action upstairs V}
\Omega(\cn_X,V_x)^{\cnu}\star \jmath_!(\IC_{\Bun^{\cmu-\cnu}_B})\to
\on{H}'{}^V_x\left(\jmath_!(\IC_{\Bun^\cmu_B})\right).
\end{equation}

\smallskip

\noindent{\em (B)}
For $\cla\in \cLambda^{pos}$ the following diagram is commutative:
$$
\CD
\Omega(\cn_X,V_x)^{\cnu}\star \Omega(\cn_X)^{-\cla}\star 
\jmath_!(\IC_{\Bun^{\cmu-\cnu+\cla}_B}) @>{\text{\eqref{action upstairs}}}>>
\Omega(\cn_X,V_x)^{\cnu}\star \jmath_!(\IC_{\Bun^{\cmu-\cnu}_B}) \\
@V{\text{\eqref{mult Omega V}}}VV     @V{\text{\eqref{action upstairs V}}}VV   \\
\Omega(\cn_X,V_x)^{\cnu-\cla}\star \jmath_!(\IC_{\Bun^{\cmu-\cnu+\cla}_B}) 
@>{\text{\eqref{action upstairs V}}}>>
\on{H}'{}^V_x\left(\jmath_!(\IC_{\Bun^\cmu_B})\right).
\endCD
$$

\smallskip

\noindent{\em (C)} The map \eqref{action upstairs V} is compatible with
filtrations, and resulting map on the associated graded level makes the
following diagram commutative:
$$
\CD
\on{gr}^{\cnu'}\left(\Omega(\cn_X,V_x)^{\cnu}\right)\star 
\jmath_!(\IC_{\Bun^{\cmu-\cnu}_B}) @>>> \on{gr}^{\cnu'}
\left(\on{H}'{}^V_x\left(\jmath_!(\IC_{\Bun^\cmu_B})\right)\right) \\
@V{\sim}VV     @V{\sim}VV  \\
V(\cnu')\otimes \delta^{\cnu'}_{x}\star 
\Omega(\cn_X)^{\cnu-\cnu'} \star \jmath_!(\IC_{\Bun^{\cmu-\cnu}_B})
@>{\text{\eqref{action upstairs}}}>> V(\cnu')\otimes \delta^{\cnu'}_{x}\star 
\jmath_!(\IC_{\Bun_B^{\cmu-\cnu'}}).
\endCD
$$

\end{thm}

The proof will be given in \secref{proof of Hecke ext}.

\ssec{}

Let us show how \thmref{Hecke of extension by zero} implies \thmref{Hecke}.
Tensoring both sides of \eqref{action upstairs V} by $(\fqb^\cmu)^*(\CS(E_\cT))$,
we obtain a map
\begin{multline*}
({}_{\infty\cdot x}\imathb_\cnu)_!\biggl(
\Omega(\cn_{X,E_{\cT}},V_{E_{\cT,x}})^{\cnu}\boxtimes 
\jmath_!\left(\IC_{\Bun_B^{\cmu-\cnu}}\otimes (\fq^{\cmu-\cnu})^*(\CS(E_\cT))\right)\biggr)\to \\
\to \on{H}'{}^V_x\left(\jmath_!\left(\IC_{\Bun^\cmu_B}\otimes (\fq^{\cmu})^*(\CS(E_\cT))\right)\right).
\end{multline*}
Taking the direct image of both sides of the above expression with respect
to the natural morphism $\fpb^\cmu:{}_{\infty\cdot x}\BunBb^\cmu\to \Bun_G$,
we obtain a map
$$R(V_x)_{E_\cT}^\cnu\otimes \Eis_!^{\cmu-\cnu}(E_\cT)\simeq
H\left({}_{\infty\cdot x}X^\cnu,\Omega(\cn_{X,E_{\cT}},V_{E_{\cT,x}})^{\cnu}\right)
\otimes \Eis_!^{\cmu-\cnu}(E_\cT)\to \on{H}_x^V(\Eis^\cmu_!(E_\cT)).$$
Summing up over all $\cnu$ and $\cmu$,
we obtain a map
\begin{equation} \label{map before ten prod}
R(V_x)_{E_\cT}\otimes
\Eis_!(E_\cT)\to \on{H}_x^V(\Eis_!(E_\cT)).
\end{equation}

Point (B) of \thmref{Hecke of extension by zero} implies that the latter map
respects the action of the algebra $R(V_x)_{E_\cT}$, i.e., we obtain a map
\begin{equation} \label{desired Hecke map for Eis}
R(V_x)_{E_\cT}\underset{R_{E_\cT}}\otimes 
\Eis_!(E_\cT)\to \on{H}_x^V(\Eis_!(E_\cT)).
\end{equation}

We claim that the above map is an isomorphism, implying \thmref{Hecke}.

\medskip

Indeed, the filtration on $\H'_x \underset{\Bun_G,\hr}\times \Bun_B^\cmu$
defines a filtration on each $\on{H}_x^V(\Eis^\cmu_!(E_\cT))$ with
$$\on{gr}^\cnu\left(\on{H}_x^V(\Eis^\cmu_!(E_\cT))\right)\simeq V(\cnu)\otimes
\Eis^{\cmu-\cnu}_!(E_\cT).$$ By point (C) of \thmref{Hecke of extension by zero} 
we obtain that the map \eqref{map before ten prod} respects the filtrations,
and the corresponding
map on the associated graded level can be identified with
$$V\otimes R_{E_\cT} \otimes
\Eis_!(E_\cT) \to V\otimes \Eis_!(E_\cT),$$
given by the $R_{E_\cT}$-action on $\Eis_!(E_\cT)$.

Hence, the map \eqref{desired Hecke map for Eis} also respects the
filtrations, and on the associated graded level induces  
the identity isomorphism of $V\otimes \Eis_!(E_\cT)$, implying that
\eqref{desired Hecke map for Eis} is itself an isomorphism.

\ssec{}

Recall from \thmref{FFKM thm} that for each $\cla\in \cLambda^{pos}$ we have 
a canonical map in $D^b(\BunBb^\cmu)$:
\begin{equation} \label{Fink action}
\fU(\fn_X)^\cla\star \IC_{\BunBb^{\cmu+\cla}}\to \IC_{\BunBb^{\cmu}},
\end{equation}
and the Verdier dual map
\begin{equation} \label{dual Fink}
\IC_{\BunBb^{\cmu}}\to \fU(\fn_X)^{*,-\cla}\star \IC_{\BunBb^{\cmu+\cla}}.
\end{equation}

In this subsection we will study how these maps are compatible with
the isomorphisms of \eqref{IC Hecke}. Namely, we will establish the following:

\begin{thm}  \label{compat Hecke}
For $V\in \Rep(\cG)$ and $\cla$ as above the following diagram is 
commutative
$$
\CD
\on{H}'{}^V_x(\IC_{\BunBb^{\cmu}})  @>{\text{\eqref{dual Fink}}}>>
\on{H}'{}^V_x\left(\fU(\fn_X)^{*,-\cla}\star \IC_{\BunBb^{\cmu+\cla}}\right) \\
& & @V{\sim}VV  \\
@V{\text{\eqref{IC Hecke}}}VV  \fU(\fn_X)^{*,-\cla}\star \on{H}'{}^V_x(\IC_{\BunBb^{\cmu+\cla}}) \\
& &   @V{\text{\eqref{IC Hecke}}}VV  \\
\underset{\cnu}\oplus\, V(\cnu)\otimes \delta_x^\cnu\star
\IC_{\BunBb^{\cmu-\cnu}} @>>>
\underset{\cnu'}\oplus\, V(\cnu')\otimes \delta_x^{\cnu'}\star 
\fU(\fn_X)^{*,-\cla}\star \IC_{\BunBb^{\cmu-\cnu'+\cla}},
\endCD
$$
where the lower horizontal arrow is comprised of the maps
$$V(\cnu)\otimes \delta_x^\cnu\star
\IC_{\BunBb^{\cmu-\cnu}}   \overset{\text{\eqref{dual Fink}}}
\longrightarrow V(\cnu)\otimes \delta_x^{\cnu}\star 
\fU(\fn_X)^{*,-\cla}\star \IC_{\BunBb^{\cmu-\cnu+\cla}}$$
and the map
$$V(\cnu)\otimes \delta_x^\cnu\star
\IC_{\BunBb^{\cmu-\cnu}} \to
V(\cnu+\cla)\otimes \fU(\fn_X)^{*,-\cla}\star \delta_x^{\cnu+\cla}\star
\IC_{\BunBb^{\cmu-\cnu}},$$
given by the identification of the !-stalk of $\fU(\fn_X)^{*,-\cla}$
at $\cla\cdot x\in X^\cla$ with $(U(\cn)^\cla)^*$ and the map
$$V(\cnu)\to V(\cnu+\cla)\otimes (U(\cn)^\cla)^*,$$
given by the action of $U(\cn)$ on $V$.
\end{thm}

The rest of this section will be devoted to the proof of this theorem.
We will need some preliminaries. Consider the graded perverse sheaf
on $_{\infty\cdot x}X^\cnu$
$$\on{Kosz}(V_x)^{\bullet,\cnu}:=\underset{\cla\in \cLambda^{pos}}\bigoplus\,
\fU(\cn_X)^{\bullet,-\cla,*}\star \Omega(\cn_X,V_x)^{\cnu+\cla}.$$

The maps \eqref{mult Omega V} define on $\on{Kosz}(V_x)^{\bullet,\cnu}$
a differential, by the same formula as in the case of 
$\on{Kosz}_{\BunBb^\cmu}^{\bullet,\cnu}$. The canonical projection
$$\Omega(\cn_X,V_x)^\cnu\to \on{gr}^\cnu(\Omega(\cn_X,V_x)^\cnu)\simeq
\delta^\cnu_x\otimes V(\cnu)$$ defines a map
\begin{equation} \label{Kosz res V}
\on{Kosz}(V_x)^{\bullet,\cnu}\to \delta^\cnu_x\otimes V(\cnu).
\end{equation}

\begin{lem}  \label{Kosz res V lem} \hfill

\smallskip

\noindent{\em (1)} The map \eqref{Kosz res V} is a quasi-isomorphism.

\smallskip

\noindent{\em (2)} For $\cla\in \cLambda^{pos}$, the diagram
$$
\CD
\on{Kosz}(V_x)^{\bullet,\cnu} @>>> \fU(\cn_X)^{\bullet,-\cla,*}\star
\on{Kosz}(V_x)^{\bullet,\cnu+\cla} \\
@VVV    @VVV    \\
\delta^\cnu_x\otimes V(\cnu) @>>>  \fU(\cn_X)^{-\cla,*}\star 
\delta^{\cnu+\cla}_x\otimes V(\cnu+\cla)
\endCD
$$
commutes in the derived category, where the top horizontal arrow
is induced by \eqref{Koszul module}, and the bottom horizontal arrow
is as \thmref{compat Hecke}.
\end{lem}

Recall the complex
$$\on{Kosz}_{\BunBb^\cmu}^{\bullet}:=
\underset{\cla\in \cLambda^{pos}}\bigoplus\,
\fU(\cn_X)^{\bullet,-\cla,*}\star \jmath_!(\IC_{\Bun_B^{\cmu+\cla}}).$$
According to \thmref{exactness of Koszul upstairs}, it is
quasi-isomorphic to $\IC_{\BunBb^\cmu}$. Moreover, by
\lemref{expl Fink action}, the map \eqref{dual Fink} is
represented by the map
\begin{equation} \label{explicit dual Fink}
\on{Kosz}_{\BunBb^\cmu}^{\bullet}\to 
\fU(\cn_X)^{\bullet,-\cla,*}\star \on{Kosz}_{\BunBb^{\cmu+\cla}}^{\bullet},
\end{equation}
induced by the map dual to \eqref{mult dot U}.

\medskip

Next, we shall reinterpret the isomorphisms of \thmref{Hecke of extension by zero}
and \eqref{IC Hecke}.
Consider the graded perverse sheaf on $_{\infty\cdot x}\BunBb^\cmu$
$$\on{Kosz}(V_x)_{\BunBb^\cmu}^{\bullet}:=
\underset{\cnu\in \cLambda,\cla\in \cLambda^{pos}}\bigoplus\,
\Omega(\cn_X,V_x)^\cnu\star \fU(\cn_X)^{\bullet,-\cla,*}
\star \jmath_!(\IC_{\Bun^{\cmu-\cnu+\cla}_B}).$$
It also acquires a natural differential differential. The projection
on the direct summand corresponding to $\cla=0$ and 
the map \eqref{action upstairs V} define a map of complexes:
\begin{equation} \label{Hecke for !}
\on{Kosz}(V_x)_{\BunBb^\cmu}^{\bullet}\to 
\on{H}'{}^V_x\left(\jmath_!(\IC_{\Bun_B^{\cmu}})\right).  
\end{equation}

\begin{prop}  \label{expl Hecke for !}
The map \eqref{Hecke for !} is a quasi-isomorphism.
\end{prop}

\begin{proof}

Considering the filtrations on both sides as in 
\thmref{Hecke of extension by zero}(C) and passing
to the associated graded we obtain a map, whose $\cnu'$-component
is
$$V(\cnu') \otimes \delta_x^{\cnu'}\star
\left(\underset{\cla,\cnu}\oplus\, \Omega(\cn_X)^{\cnu-\cnu'}\star
\fU(\cn_X)^{\bullet,-\cla,*}\star \jmath_!(\IC_{\Bun^{\cmu-\cnu+\cla}_B})\right)\to
V(\cnu')\otimes \delta_x^{\cnu'}\star
\jmath_!(\IC_{\Bun^{\cmu-\cnu'}}).$$
Regrouping the terms, the RHS of the above expression is the direct sum
over $\ceta$ of 
$$V(\cnu') \otimes \delta_x^{\cnu'}\star\left(\underset{\cla}\oplus\, 
\Omega(\cn_X)^{\ceta+\cla}\star 
\fU(\cn_X)^{\bullet,-\cla,*}\star \jmath_!(\IC_{\Bun^{\cmu-\ceta-\cnu'}_B})\right).$$

However, it follows from the acyclicity of $\on{Kosz}(E_\cT)^{\bullet,-\ceta}$
that the latter expression is acyclic unless $\ceta=0$, and quasi-isomorphic to
$V(\cnu') \otimes \delta_x^{\cnu'}\star\jmath_!(\IC_{\Bun^{\cmu-\cnu'}_B})$ in the latter 
case, as required.

\end{proof}

Let us apply the functor $\on{H}'{}^V_x$ term-wise to the complex 
$\on{Kosz}_{\BunBb^\cmu}^{\bullet}$. By \propref{expl Hecke for !}, the result
is quasi-isomorphic to the complex
\begin{equation}  \label{Hecke to Kosz}
\underset{\cla}\bigoplus\, \fU(\cn_X)^{\bullet,-\cla,*}\star 
\on{Kosz}(V_x)_{\BunBb^{\cmu+\cla}}^{\bullet}
\end{equation}
with a natural differential. Using \eqref{Kosz res V} we obtain a map
from the expression in \eqref{Hecke to Kosz} to
$$\underset{\cnu'}\bigoplus\, V(\cnu')\otimes \delta_x^{\cnu'}\star 
\on{Kosz}_{\BunBb^{\cmu-\cnu'}}^{\bullet}.$$

The next assertion follows from \thmref{Hecke of extension by zero}:

\begin{lem}  \label{expl IC Hecke}
The following diagram commutes in the derived category:
$$
\CD
\underset{\cla}\bigoplus\, \fU(\cn_X)^{\bullet,-\cla,*}\star 
\on{Kosz}(V_x)_{\BunBb^{\cmu+\cla}}^{\bullet}   @>>>  
\underset{\cnu'}\bigoplus\, V(\cnu')\otimes \delta_x^{\cnu'}\star 
\on{Kosz}_{\BunBb^{\cmu-\cnu'}}^{\bullet}   \\
@V{\sim}VV     @V{\sim}VV   \\
\on{H}'{}^V_x(\IC_{\BunBb^\cmu})   @>>>  
\underset{\cnu'}\bigoplus\, V(\cnu')\otimes \delta_x^{\cnu'}\star \IC_{\BunBb^{\cmu-\cnu'}}.
\endCD
$$
\end{lem}

Now, using \lemref{Kosz res V lem}, and \lemref{expl IC Hecke}, we can prove
\thmref{compat Hecke}. Indeed, the diagram, whose commutativity we have to
prove, is equivalent to the following one, which is manifestly commutative:
$$
\CD
\on{H}'{}^V_x\left(\on{Kosz}_{\BunBb^{\cmu}}^{\bullet}\right) 
@>{\text{\eqref{explicit dual Fink}}}>>
\on{H}'{}^V_x\left(\fU(\cn_X)^{\bullet,-\cla,*}\star 
\on{Kosz}_{\BunBb^{\cmu+\cla}}^{\bullet}\right) \\
& & @V{\sim}VV  \\
@V{\sim}VV    \fU(\cn_X)^{\bullet,-\cla,*}\star 
\on{H}'{}^V_x\left(\on{Kosz}_{\BunBb^{\cmu+\cla}}^{\bullet}\right) \\
& & @VVV  \\
\underset{\cla'}\bigoplus\, \fU(\cn_X)^{\bullet,-\cla',*}\star 
\on{Kosz}(V_x)_{\BunBb^{\cmu+\cla'}}^{\bullet} @>>>
\underset{\cla''}\bigoplus\, \fU(\cn_X)^{\bullet,-\cla,*}\star \fU(\cn_X)^{\bullet,-\cla'',*}\star
\on{Kosz}(V_x)_{\BunBb^{\cmu+\cla+\cla''}}^{\bullet},
\endCD
$$
where the bottom arrow is given by
$$\fU(\cn_X)^{\bullet,-\cla',*}\to \fU(\cn_X)^{\bullet,-\cla,*}\star \fU(\cn_X)^{\bullet,-\cla'',*}$$
for $\cla'=\cla+\cla''$.


\section{Proof of \thmref{Hecke of extension by zero}}    \label{proof of Hecke ext}

\ssec{}

The strategy of the proof of \thmref{Hecke of extension by zero} will be
largely parallel to that of \thmref{extension by zero}. Without restricting
the generality, we can assume that the representation $V$ is irreducible,
i.e., $V=V^\ceta$ for some highest weight $\ceta\in \cLambda^+$.
Then the support of $\on{H}'{}^V_x\left(\jmath_!(\IC_{\Bun_B^\cmu})\right)$
is contained in the closed substack of $_{\infty\cdot x}\BunBb^\cmu$
equal to the image of 
$$\BunBb^{\cmu-\ceta}\simeq (\ceta\cdot x)\times \BunBb^{\cmu-\ceta}
\overset{_{\infty\cdot x}\imathb_{\ceta}}\hookrightarrow {}
_{\infty\cdot x}\BunBb^\cmu.$$

The perverse sheaf $\Omega(\cn_X,V^\ceta_x)^\cnu$ is non-zero only
for $\cnu\leq \ceta$ and is supported on the subscheme
$$X^{\ceta-\cnu}\times (\ceta\cdot x)\subset{}_{\infty\cdot x}X^\cnu.$$

From now on, we will work on $\BunBb^{\cmu-\ceta}$ and
$X^{\ceta-\cnu}$ rather than on $_{\infty\cdot x}\BunBb^\cmu$ and
$_{\infty\cdot x}X^\cnu$.

\medskip

Consider the open subset $\overset{\circ}{(X-x)}{}^{\ceta-\cnu}\subset
X^{\ceta-\cnu}$. Let us denote its open embedding by $j_x^{\ceta-\cnu}$.
We have:
$$j_x^{\ceta-\cnu}{}^*\left(\Omega(\cn_X,V^\ceta_x)^\cnu\right)\simeq 
j_x^{\ceta-\cnu}{}^*\left(\Omega(\cn_X)^{\cnu-\ceta}\right) \otimes
V^\ceta(\ceta).$$
First, we claim that the isomorphism stated in \eqref{need to construct V}
holds over 
$$\overset{\circ}{(X-x)}{}^{\ceta-\cnu}\times \Bun_B^{\cmu-\cnu}\subset
X^{\ceta-\cnu}\times \Bun_B^{\cmu-\cnu},$$
i.e.,
\begin{equation} \label{map over open V}
j_x^{\ceta-\cnu}{}^*\left(\Omega(\cn_X,V^\ceta_x)^\cnu\right)\boxtimes
\IC_{\Bun_B^{\cmu-\cnu}}\simeq
(j_x^{\ceta-\cnu}\times \on{id})^*\biggl(
h^0\biggl({}_{\infty\cdot x}\imath_\cnu^!
\left(\on{H}'{}^{V^\ceta}_x\left(\jmath_!(\IC_{\Bun^\cmu_B})\right)\right)\biggr)\biggr).
\end{equation}
This follows from the definition of the functor $\on{H}'{}^V_x$
and \eqref{map over open}.

\medskip

We have the following assertion, parallel to \lemref{injectivity for omega}:

\begin{lem} \label{injectivity for omega V}
The canonical map
$$\Omega(\cn_X,V^\ceta_x)^\cnu\to
j_x^{\ceta-\cnu}{}_*\circ j_x^{\ceta-\cnu}{}^* 
\left(\Omega(\cn_X,V^\ceta_x)^\cnu\right)$$
is injective.
\end{lem}

\begin{proof}
Let $i_x^{\ceta-\cnu}$ denote the embedding of the point 
$(\ceta-\cnu)\cdot x$ into $X^{\ceta-\cnu}$. By induction on $|\ceta-\cnu|$ 
we have to show that
$$i_x^{\ceta-\cnu}{}^!\left(\Omega(\cn_X,V^\ceta_x)^\cnu\right)$$
has no cohomologies in degrees $\leq 0$, unless $\cnu=\ceta$.

By the definition of $\Omega(\cn_X,V^\ceta_x)^\cnu$, the above
!-stalk is quasi-isomorphic to the weight $\cnu$-component 
in the cohomology $H^\bullet(\cn,V^\ceta)$. Hence,
the $0$-th cohomology, which corresponds to the highest
weight in $V^\ceta$, is of weight $\ceta$, and not $\cnu$.
\end{proof}

Another assertion that we shall need is parallel to 
\propref{injectivity for !}:

\begin{prop}  \label{injectivity for ! V}
The canonical map
$$h^0\biggl({}_{\infty\cdot x}\imath_\cnu^!
\left(\on{H}'{}^{V^\ceta}_x\left(\jmath_!(\IC_{\Bun^\cmu_B})\right)\right)\biggr)
\to (j^\cla\times \on{id})_*\circ (j^\cla\times \on{id})^*
\biggl(h^0\biggl({}_{\infty\cdot x}\imath_\cnu^!
\left(\on{H}'{}^{V^\ceta}_x\left(\jmath_!(\IC_{\Bun^\cmu_B})\right)\right)\biggr)\biggr)$$
is injective.
\end{prop}

Let us assume this proposition and proceed with the proof of
\thmref{Hecke of extension by zero}.

\medskip

As in the proof of \thmref{extension by zero}, from \lemref{injectivity for omega V}
and \propref{injectivity for ! V}, we obtain that if the isomorphism \eqref{map over open V}
extends to a map in one direction
\begin{equation} \label{need to construct V one dir}
\Omega(\cn_X,V^\ceta_x)^\cnu\times \IC_{\Bun_B^{\cmu-\cnu}}\to
h^0\biggl({}_{\infty\cdot x}\imath_\cnu^!
\left(\on{H}'{}^V_x\left(\jmath_!(\IC_{\Bun^\cmu_B})\right)\right)\biggr),
\end{equation}
then it does so uniquely. Moreover, this map will automatically be
an isomorphism and point (B) of \thmref{Hecke of extension by zero}
will hold.

Thus, our present goal will be to establish the required extension property.

\ssec{}

We will distinguish two cases. One is when $\ceta-\cnu$ is a multiple of
a simple coroot, and another when it is not. Let us first treat the latter case.
We will argue by induction on $|\ceta-\cnu|$, so we can assume that
the extension \eqref{need to construct V one dir} exists for all
$\cnu'$ with $\cnu'>\cnu$.

Since the assertion about extension is local, 
as in the proof of \thmref{extension by zero}, we can pass from
$\BunBb^\cmu$ to a suitable version of the Zastava space, 
and using factorization, we can assume that the required
isomorphism holds over the open substack
\begin{equation} \label{remove point}
\left(X^{\ceta-\cnu}-((\ceta-\cnu)\cdot x)\right)\times \Bun_B^{\cmu-\cnu}.
\end{equation}

As in the proof of \thmref{extension by zero}, to establish the required
extension property, it suffices to prove the following:

\begin{lem}
Assume that $|\ceta-\cnu|$ is not a multiple of a simple coroot.
Then 
$$\on{Ext}^1_{X^{\ceta-\cnu}}
\left(i^{\ceta-\cnu}_x{}_!(\BC),\Omega(\cn_X,V^\ceta_x)^\cnu\right)=0.$$
\end{lem}

\begin{proof}
The $\on{Ext}^1$ of the lemma is isomorphic to the weight $\cnu$ component
in $H^1(\cn,V^\ceta)$. However, the 1-st cohomology in question consists
of weights of the form 
$$s_i(\ceta+\check\rho)-\check\rho=\ceta-\check\alpha_i\cdot (\langle \alpha_i,
\ceta\rangle+1),$$
for simple reflections $s_i$. The difference between such a weight and $\ceta$
is a multiple of $\check\alpha_i$.
\end{proof}

\ssec{}

To complete the proof of \thmref{Hecke of extension by zero} it remains to
do three things: (1) prove the extension property for $\cnu=\ceta-n\cdot \check\alpha_i$,
(2) prove \propref{injectivity for ! V}, and (3) establish the compatibility of
point (C) of the theorem. We will do this by reducing to the case of groups
of semi-simple rank one.

\medskip

Thus, let us assume for a moment that $G$ is of semi-simple rank $1$, and show that
\thmref{Hecke of extension by zero} holds in this case. First, it is easy to see that
we can assume that $G=GL_2$, as the statement of the theorem is stable
under isogenies. 

\medskip

Next, we claim that the explicit computation we did in 
\secref{expl for sl(2)} for $V$ being the standard representation $V^{(1,0)}$
amounts to the statement of \thmref{Hecke of extension by zero} in this case. 
Indeed, the perverse sheaf $\Omega(\cn_X,V^\ceta_x)^\cnu$ is non-zero only
for $\cnu$ being $(1,0)$ or $(0,1)$, and is the sky-scraper in the former
case and the sheaf $j_x{}_*(\IC_{X-x})$ in the latter case.

Hence, the isomorphism assertion of \thmref{Hecke of extension by zero}(A)
follows from \eqref{sl(2) Hecke map}. The assertion of 
\thmref{Hecke of extension by zero}(C) is also manifest.

\medskip

Now, we claim that \thmref{Hecke of extension by zero} holds for an arbitrary
representation $V$ of $\cG\simeq GL_2$. Indeed, it is easy to see 
from \eqref{tensor Omega V} (and this is valid for any group $G$) 
that if \thmref{Hecke of extension by zero} 
holds for representations $V^1$ and $V^2$, then it also holds for their
tensor product.

In addition, by \lemref{injectivity for omega V} and \propref{injectivity for ! V},
we obtain that if \thmref{Hecke of extension by zero} holds for some
representation $V$, then it also holds for any of its direct summands.

To prove \thmref{Hecke of extension by zero} for a representation $V$
of $GL_2$ it suffices to observe that $V$ is isomorphic to a direct summand
of a tensor power of the standard representation, up to some power
of the determinant. 

\ssec{}

Next, we shall study how the Hecke functors $\on{H}'{}^V_x$, that act on
the derived category on $_{\infty\cdot x}\BunBb^\cmu$, are related to
similar functors for Levi subgroups of $G$.

Let $P$ be a parabolic in $G$, and let $M$ be the corresponding Levi
quotient; let $B(M)$ denote the Borel subgroup in $M$. If we
choose $P$ so that it contains $B$, then $B(M)$ is the projection of
$B$ to $M$. Let $\BunBMb$ and $_{\infty\cdot x}\BunBMb$
be the corresponding stacks for the group $M$.

Note that we have a natural isomorphism:
$$\Bun_P\underset{\Bun_M}\times \BunBM\simeq \Bun_B,$$
which extends to a locally closed embedding
$$\Bun_P\underset{\Bun_M}\times{}_{\infty\cdot x}\BunBMb\overset{\imath_M}
\hookrightarrow {}_{\infty\cdot x}\BunBb.$$

Recall also that we have a diagram of affine Grassmannians:
$$\Gr_G\overset{\fp_\Gr}\hookleftarrow \Gr_P\overset{\fq_\Gr}\twoheadrightarrow \Gr_M,$$
and the restriction functor $\Rep(\cG)\to \Rep(\cM)$ corresponds to the following
operation on spherical perverse sheaves:
$$\CV\in \on{Perv}(\Gr_G)\mapsto \fq_\Gr{}_!\circ \fq_\Gr^*(\CV)\in \on{Perv}(\Gr_M),$$
up to a cohomological shift.

\medskip

For $U\in \Rep(\cM)$, let $\on{H}'{}^U_{M,x}$ denote the corresponding Hecke functor
acting on the derived category on $_{\infty\cdot x}\BunBMb$. Let us denote by 
$\on{H}'{}^U_{P,x}$ the Hecke functor acting on the derived category on
the stack $\Bun_P\underset{\Bun_M}\times{}_{\infty\cdot x}\BunBMb$. These
functors are compatible via the pull-back functor corresponding to the projection
$$\Bun_P\underset{\Bun_M}\times{}_{\infty\cdot x}\BunBMb\to {}_{\infty\cdot x}\BunBMb.$$

For $\CF\in D^b({}_{\infty\cdot x}\BunBb)$ we have a functorial isomorphism
\begin{equation} \label{Hecke G & P}
\imath_M^*\left((\on{H}'{}^V_{x}(\CF)\right)\simeq \on{H}'{}^{U}_{P,x}
\left(\imath_M^*(\CF)\right),
\end{equation}
where $U=\on{Res}^\cG_\cM(V)$.

Let us observe that when we apply \eqref{Hecke G & P} to
$\CF=\jmath_!(\IC_{\Bun^\cmu_B})$, the resulting isomorphism is
compatible with the filtrations.

\ssec{}

Let us now establish the extension property of the isomorphism \eqref{map over open V}
to a morphism \eqref{need to construct V one dir} when $\ceta-\cnu$ is a multiple of
a simple coroot, say $\check\alpha_i$. We take $P$ to be the minimal parabolic,
corresponding to the chosen vertex of the Dynkin diagram; the corresponding Levi
$M$ is of semi-simple rank one.

Note that the map
$$\Bun_P\underset{\Bun_M}\times \BunBMb^\ceta\to \BunBb^\ceta,$$
induced by $\imath_M$, is an open embedding (which is true for any parabolic),
and its image contains the stratum $X^{\ceta-\cnu}\times \Bun_B^{\cmu-\cnu}$. 

The required assertion follows now from \eqref{Hecke G & P} and the fact that
\thmref{Hecke of extension by zero} holds for $M$.

\ssec{}

Let us now prove \propref{injectivity for ! V}. We argue by induction on $|\ceta-\cnu|$.
The assertion is evidently true for $\cnu=\ceta$, and we assume that it holds for
all $\cnu'>\cnu$. As in the proof of \propref{injectivity for !}, a factorization argument
reduces the assertion of the proposition to the fact that there are no non-zero morphisms
\begin{equation} \label{violating map}
\delta^\cnu_x\star \jmath_!(\IC_{\Bun^{\cmu-\cnu}_B})\to
\on{H}'{}^{V^\ceta}_x\left(\jmath_!(\IC_{\Bun^\cmu_B})\right).
\end{equation}

Let us suppose by contradiction that a morphism like this existed. Let us recall
the filtration on $\on{H}'{}^{V^\ceta}_x\left(\jmath_!(\IC_{\Bun^\cmu_B})\right)$ introduced
in \secref{structure of Hecke}. The only subquotient that admits a map from
$\delta^\cnu_x\star \jmath_!(\IC_{\Bun^{\cmu-\cnu}_B})$ is 
$$\on{gr}^\cnu\biggl(\on{H}'{}^{V^\ceta}_x\left(\jmath_!(\IC_{\Bun^\cmu_B})\right)\biggr)
\simeq V^\ceta(\cnu)\otimes \delta^\cnu_x\star \jmath_!(\IC_{\Bun^{\cmu-\cnu}_B}).$$

Hence, a map in \eqref{violating map} defines an element $\bv\in V^\ceta(\cnu)$. Let $\check\alpha_i$
be a simple root of $\cg$, for which $\bv$ is {\it not} a highest weight vector. (Such $i\in I$ exists,
for otherwise $\bv$ would have been a highest weight vector for $\cg$, which is impossible,
since $\ceta\neq \cnu$.)

\medskip

Let us consider the restriction of the map \eqref{violating map} under $\imath_M$ 
for $M$ being the Levi of the minimal parabolic, corresponding to the vertex $i$.
We obtain a map
\begin{equation} \label{violating rk one}
\delta^\cnu_x\star \jmath_!(\IC_{\Bun^{\cmu-\cnu}_{B(M)}})\to
\on{H}'{}^{U}_{M,x}\left(\jmath_!(\IC_{\Bun^\cmu_{B(M)}})\right),
\end{equation}
whose projection to the subquotient
$$\on{gr}^{\cnu}\biggl(\on{H}'{}^{U}_{M,x}\left(\jmath_!(\IC_{\Bun^\cmu_{B(M)}})\right)\biggr)\simeq
U(\cnu)\otimes \delta^\cnu_x\star \jmath_!(\IC_{\Bun^{\cmu-\cnu}_{B(M)}}),$$
corresponds to the same vector $\bv\in U:=\on{Res}^{\cG}_{\cM}(V^\ceta)$.
 
But this is a contradiction: we know that in the rank one situation the only
non-trivial maps as in \eqref{violating rk one} correspond to 
{\it highest weight vectors} in $U(\cnu)$.

\ssec{}

Finally, let us prove point (C) of \thmref{Hecke of extension by zero}.
By induction, we can assume that the assertion holds for all $\cnu_1<\cnu$.
Let us note that the penultimate terms of the filtration on both sides of
\eqref{action upstairs V} are isomorphic to the intermediate extensions
of their own restrictions to the open substack \eqref{remove point}.

Using the Zastava model and factorization, this implies that the map
\eqref{action upstairs V} is compatible with filtrations, and that the
induced maps on $\on{gr}^{\cnu'}$ for $\cnu'<\cnu$ are as required.
The map on the last term of the filtration
$$V(\cnu)\otimes \delta^{\cnu}_{x}\star 
\jmath_!(\IC_{\Bun^{\cmu-\cnu}_B}) \to
V(\cnu')\otimes \delta^{\cnu}_{x}\star \jmath_!(\IC_{\Bun_B^{\cmu-\cnu}})$$
corresponds, therefore, to an endomorphism $\phi(\cnu)$ of $V(\cnu)$. 

We have to show that $\phi(\cnu)$ equals the identity. For that it is enough to 
show that for every simple root $\check\alpha_i$ of $\cg$ the diagram
\begin{equation} \label{action of i}
\CD
\cn^{\check\alpha_i}\otimes V(\cnu)  @>>>  V(\cnu+\check\alpha_i)  \\
@V{\on{id}\otimes \phi(\cnu)}VV @V{\on{id}}VV  \\
\cn^{\check\alpha_i}\otimes V(\cnu)  @>>> V(\cnu+\check\alpha_i)
\endCD
\end{equation}
commutes.

\medskip

Let $P$ and $M$ be as above. Note that the substack
$$\left(\Bun_P\underset{\Bun_M}\times {}_{\infty\cdot x}\BunBMb\right)
\underset{_{\infty\cdot x}\BunBb}\times 
\left(X^{\ceta-\cnu}\times \Bun_B^{\cmu-\cnu}\right)$$
identifies with 
$$X^{(n_i)}\times  \Bun_B^{\cmu-\cnu}\subset X^{\ceta-\cnu}\times \Bun_B^{\cmu-\cnu},$$
where $n_i=\langle \alpha_i,\ceta-\cnu\rangle$. The $0$-th perverse cohomology
of the *-restriction on $\Omega(\cn_X,V_x)^{\cnu}$ under
$$X^{(n_i)}\hookrightarrow X^{\ceta-\cnu}$$
identifies with the corresponding perverse sheaf 
$\Omega(\cn^{\check\alpha_i}_X,U_x)^{\cnu}$ on $X^{(n_i)}$. It can be also thought
of as the maximal quotient of $\Omega(\cn_X,V_x)^{\cnu}$ supported on the above
subscheme, and it equals the quotient of $\Omega(\cn_X,V_x)^{\cnu}$ corresponding to
the terms of the canonical filtration with $\cnu'\in \ceta-\check\alpha_i\cdot \BZ^{\geq 0}$.

Let us apply the functor $\imath_M^*$ to the morphism \eqref{action upstairs V}. 
We obtain a map
\begin{equation} \label{action M}
\Omega(\cn^{\check\alpha_i}_X,U_x)^{\cnu}\star
\jmath_!(\IC_{\Bun_{B(M)}^{\cmu-\cnu}})\to 
\on{H}'{}^{U}_{M,x}\left(\jmath_!(\IC_{\Bun^\cmu_{B(M)}})\right).
\end{equation}
The restriction of this map to
$$(X^{(n_i)}-n_i\cdot x)\times \Bun^{\cmu-\cnu}_{B(M)}\subset 
X^{(n_i)}\times \Bun^{\cmu-\cnu}_{B(M)}$$
equals the map \eqref{action upstairs V} for the group $M$, by the
induction hypothesis. Hence, by \propref{injectivity for ! V} (applied to $M$),
the map \eqref{action M} equals \eqref{action upstairs V}, at least when
restricted to the direct summand, corresponding to the direct summand of
$U=\on{Res}^{\cG}_{\cM}(V^\ceta)$, denoted $U^{\neq \cnu}$,
complementary to the isotypic component of highest weight $\cnu$.

This implies that the map $\phi:V(\cnu)\simeq U(\cnu)\to U(\cnu)\simeq V(\cnu)$ 
equals the identity map, when restricted to $U^{\neq \cnu}$. This
implies the commutativity of \eqref{action of i}.

\section{Regular Eisenstein series are perverse}   \label{eis is perv}

\ssec{}

The goal of this section is to prove \thmref{Eisenstein is perverse}. Thus, we 
suppose that $E_\cT$ is a $\cT$-local system, which is regular, i.e.,
$\alpha(E_\cT)$ is non-trivial for all $\alpha\in \Delta^+$. 

It is easy to see that we can replace the initial group $G$ by an isogenous one with a 
connected center, and for which $[G,G]$ is still simply connected. In this case, the assumption that
$E_\cT$ is regular implies that it is {\it strongly regular}, i.e., that
$E_\cT^w\neq E_\cT$ for any non-trivial element $w$ of the Weyl group.

In this case, the following strengthening of \thmref{Eisenstein is perverse} will hold:
\begin{thm}  \label{Eisenstein is perverse'}  
Assume that $E_\cT$ is strongly regular. Then:

\smallskip

\noindent{\em(1)}
The complex $\Eisb^\cmu(E_\cT)$ is an irreducible perverse sheaf.

\smallskip

\noindent{\em(2)} 
For $\cmu'\neq \cmu$, the perverse sheaves 
$\Eisb^\cmu(E_\cT)$ and $\Eisb^{\cmu'}(E_\cT)$
are non-isomorphic.

\smallskip

\noindent{\em(3)}
The complex $\Eis_!^\cmu(E_\cT)$ is also a perverse sheaf.

\end{thm}

The main step in the proof is the following:

\begin{prop}  \label{end computation}
$R^0\Hom(\Eisb^\cmu(E_\cT),\Eisb^{\cmu'}(E_\cT))$ is $1$-dimensional
for $\cmu=\cmu'$ and vanishes for $\cmu\neq \cmu'$.
\end{prop}

Let us first explain how \propref{end computation} implies \thmref{Eisenstein is perverse'}.
Indeed, by Kashiwara's conjecture (=Drinfeld's theorem, see \cite{Dr}), the complex
$\Eisb^\cmu(E_\cT)$ is semi-simple and satisfies Hard Lefschetz. Hence,
the first assertion of the corollary implies point (1) of the theorem. Point (2)
of the theorem follows from the second assertion of the corollary.

\medskip

Finally, to see that $\Eis_!^\cmu(E_\cT)$ is perverse, recall that by
\secref{structure of open}, the complex 
$$\jmath_!\left(\IC_{\Bun^\cmu_B}\otimes (\fq^\cmu)^*(\CS(E_\cT))\right)$$ on $\BunBb^\cmu$
admits a filtration by complexes of the form 
$$(\imathb_{\cmu'-\cmu})_!\left(\Omega(\cn_{X,E_\cT})^{\cmu-\cmu'}\boxtimes 
\left(\IC_{\BunBb^{\cmu'}}\otimes (\fq^{\cmu'})^*(\CS(E_\cT))\right)\right).$$
Hence, $\Eis_!^\cmu(E_\cT)$ admits a filtration by complexes of the form 
$$\Eisb^{\cmu'}(E_\cT)\otimes H(X^{\cmu'-\cmu},
\Omega(\cn_{X,E_\cT})^{\cmu-\cmu'}),$$ which are all perverse sheaves by 
point (1) of \thmref{Eisenstein is perverse'} and \corref{cohomology of ups}.

\medskip

\noindent{\it Remark.} Let us assume that the validity of the following (very plausible) 
extension of Kashiwara's conjecture:

\begin{conj}
Let $f:Y'\to Y''$ be a proper map between schemes over an algebraically 
closed field, and let $\CF$ be an irreducible $\ell$-adic perverse sheaf on 
$Y'$, which can be obtained by the 6 operations from a $1$-dimensional
local system on some curve $X$. Then $f_!(\CF)$ is semi-simple and
satisfies Hard Lefschetz.
\end{conj}

Then the above argument deducing \thmref{Eisenstein is perverse} from
\propref{end computation} would be valid in the context of $\ell$-adic sheaves
over an arbitrary ground field.

\ssec{}

From now on our goal is to prove \propref{end computation}. We shall deduce
is from another proposition, which amounts to a computation of the constant term of Eisenstein 
series in the geometric context:

\begin{prop}  \label{CT}
$R^i\Hom(\Eis_!^\cmu(E_\cT),\Eisb^{\cmu'}(E_\cT))=0$ if $i<0$. For $i=0$ it
vanishes if $\cmu\neq \cmu'$ and is $1$-dimensional for 
$\mu=\mu'$. 
\end{prop}

This proposition implies \propref{end computation} as follows:

\medskip

The perverse sheaf $\IC_{\BunBb^\cmu}\otimes \fqb^*(\CS(E_\cT))$ 
admits a filtration in the derived category by complexes of the form
$$(\imath_{\cmu'-\cmu})_!\left(\fU(\cn_{X,E_\cT})^{\bullet,\cmu'-\cmu,*}\boxtimes 
(\IC_{\Bun_B^{\cmu'}}\otimes (\fq^{\cmu'})^*(\CS_E))\right).$$ 

Hence, $\Eisb^\cmu(E_\cT)$ admits a filtration by complexes of the form
$$H\left(X^{\cmu'-\cmu},\fU(\cn_{X,E_\cT})^{\bullet,\cmu'-\cmu,*}\right)
\otimes \Eis_!^{\cmu'}(E_\cT).$$
However, since $\fU(\cn_{X,E_\cT})^{\cmu'-\cmu}$ is a {\it sheaf}, its cohomology
is concentrated in non-negative degrees. Hence, 
$H\left(X^{\cmu'-\cmu},\fU(\cn_{X,E_\cT})^{\bullet,\cmu'-\cmu,*}\right)$ is concentrated
in non-positive degrees. Moreover, the fact that $E_\cT$ is regular implies that
the above cohomology is concentrated in strictly negative degrees unless $\cmu=\cmu'$.

Therefore, $R^i\Hom(\Eis_!^\cmu(E_\cT),\Eisb^{\cmu'}(E_\cT))=0$ 
for $i<0$, as guaranteed by \propref{CT}, implies that 
$$R^0\Hom(\Eisb^\cmu(E_\cT),\Eisb^{\cmu'}(E_\cT))\simeq
R^0\Hom(\Eis_!^\cmu(E_\cT),\Eisb^{\cmu'}(E_\cT)),$$
and we deduce the assertion of \propref{end computation}.

\ssec{Proof of \propref{CT}}

Let  $$\on{CT}^\cmu:D(\Bun_G)\to D(\Bun_T^\cmu)$$ be the constant term functor, the right
adjoint to the functor $\Eis_!^\cmu$ that sends
$$\CF\in D^b(\Bun^\cmu_T)\mapsto \fp^\cmu_!\circ \fq^\cmu{}^*(\CF)[\dim(\Bun_B^\cmu)]\in D^b(\Bun_G).$$
The functor $\on{CT}^\cmu$ is thus given by
$$\CF'\in D^b(\Bun_G)\mapsto \fq^\cmu_*\circ \fp^\cmu{}^!(\CF')[-\dim(\Bun_B^\cmu)]\in D^b(\Bun^\cmu_T),$$
which is well-defined, since $\Bun_B^\cmu$ is of finite type and
the morphism $\fq^\cmu$ is a generalized affine fibration.

\medskip

Thus, we have to compute $R\Hom\left(\CS^\cmu(E_\cT),\on{CT}^\cmu(\Eisb^{\cmu'}(E_\cT))\right)$. 
Consider the diagram
$$
\CD
\Bun^\cmu_B\underset{\Bun_G}\times \BunBb^{\cmu'} @>{\fp'{}^\cmu}>> \BunBb^{\cmu'}  
@>{\fqb^{\cmu'}}>>  \Bun^{\cmu'}_T\\
@V{\fpb'{}^{\cmu'}}VV      @V{\fpb^{\cmu'}}VV  \\
\Bun_B^\cmu  @>{\fp^\cmu}>> \Bun_G  \\
@V{\fq^\cmu}VV  \\
\Bun^\cmu_T.
\endCD
$$

By base change, and using the fact that $\fpb^{\cmu'}$ is proper,
$\on{CT}^\cmu(\Eisb^{\cmu'}(E_\cT))$ is the direct image onto $\Bun^\cmu_T$
from $\Bun^\cmu_B\underset{\Bun_G}\times \BunBb^{\cmu'}$, of the complex
\begin{equation} \label{sheaf on product}
\fp'{}^\cmu{}^!\left(\IC_{\BunBb^{\cmu'}}\otimes (\fqb^{\cmu'})^*(\CS(E_\cT))\right)
[-\dim(\Bun_B^\cmu)].
\end{equation}

The stack $\Bun^\cmu_B\underset{\Bun_G}\times \BunBb^{\cmu'}$ is the union
of locally closed substacks $(\Bun^\cmu_B\underset{\Bun_G}\times \BunBb^{\cmu'})^w$
numbered by elements $w$ of the Weyl group that measure the relative position
of the two flags at the generic point of the curve. We obtain a filtration on
$\on{CT}^\cmu(\Eisb^{\cmu'}(E_\cT))$ as an object of the derived category,
whose subquotients we will denote by $\on{CT}^\cmu(\Eisb^{\cmu'}(E_\cT))^w$.

\begin{prop}   \label{terms of constant term}
The complex $\on{CT}^\cmu(\Eisb^{\cmu'}(E_\cT))^w$ on $\Bun^\cmu_T$ 
is an extension in the derived category of complexes isomorphic to
$\CS_{E_\cT^w}$, where $E_\cT^w$ is the $w$-twist of $E_\cT$.
\end{prop}

This proposition will be proved in the next subsection. Assuming it, 
let us finish the proof of \propref{CT}. First, we obtain that
$$R\Hom\left(\CS^\cmu(E_\cT), \on{CT}^\cmu(\Eisb^{\cmu'}(E_\cT))^w\right)=0$$
for any $1\neq w\in W$. Indeed, since $E_\cT$ is strongly regular, the local 
systems $\CS_{E_\cT^w}$ and $\CS_{E_\cT}$ are non-isomorphic, and hence
$R\Hom$ between them over $\Bun_T$, which is essentially an abelian variety, 
vanishes.

\medskip

Thus, it remains to analyze $\on{CT}^\cmu(\Eisb^{\cmu'}(E_\cT))^1$. We have;
$$(\Bun^\cmu_B\underset{\Bun_G}\times \BunBb^{\cmu'})^1\simeq
X^{\cmu'-\cmu}\times \Bun^{\cmu}_B,$$
and the map $\fp'{}^\cmu$ identifies with $\imath_{\cmu'-\cmu}$.
In terms of this identification, the $!$-restriction of the complex in \eqref{sheaf on product} to
$(\Bun^\cmu_B\underset{\Bun_G}\times \BunBb^{\cmu'})^1$ becomes
$$\fU(\cn_{X,E_\cT})^{\cmu'-\cmu}\boxtimes \IC_{\Bun^{\cmu}_B}[-\dim(\Bun_B^\cmu)].$$ 

Since $\Bun_B^\cmu\to \Bun_T^\cmu$ is a generalized affine fibration,
we obtain:
$$\on{CT}^\cmu(\Eisb^{\cmu'}(E_\cT))^1\simeq
\CS^\cmu(E_\cT)\otimes H\left(X^{\cmu'-\cmu},\fU(\cn_{X,E_\cT})^{\cmu'-\cmu}\right).$$
However, since $\fU(\cn_{X,E_\cT})^{\cmu'-\cmu}$ is a sheaf, the cohomology
$H\left(X^{\cmu'-\cmu},\fU(\cn_{X,E_\cT})^{\cmu'-\cmu}\right)$ is concentrated
in non-negative degrees (and strictly positive if $\cmu\neq \cmu'$ since $E_\cT$
was assumed regular).

This implies that
$$R^i\Hom\left(\CS^\cmu(E_\cT), \on{CT}^\cmu(\Eisb^{\cmu'}(E_\cT))^1\right)=0$$
for $i<0$ and for $i=0$ and $\cmu\neq \cmu'$.

\ssec{Proof of \propref{terms of constant term}}

Let $\Eis_*^{\cmu'}(E_\cT)$ denote
$\fp^\cmu_*\left(\IC_{\BunBb}\otimes (\fq^\cmu)^*(\CS(E_\cT))\right)$. I.e.,
we replace the functor $\fp^\cmu_!$ in the definition of $\Eis_!^{\cmu'}(E_\cT)$
by $\fp^\cmu_*$.

By \corref{Groth group comp} and Verdier duality,
it is enough to calculate $\on{CT}^\cmu(\Eis_*^{\cmu'}(E_\cT))^w$.
This will be parallel to the calculation of the $w$-term in the expression for
the constant term of Eisenstein series in the classical theory of automorphic forms.
Recall that in this theory the sought-for expression is known explicitly, and 
equals (the function corresponding to) $\CS(E^w_\cT)$ multiplied by an appropriate
ratio of L-functions. Unfortunately, in the present geometric context we will not be
able to obtain an explicit formula for $\on{CT}^\cmu(\Eis_*^{\cmu'}(E_\cT))^w$;
we will only show that it is isomorphic to $\CS(E^w_\cT)$ times {\it some} complex
of vector spaces.

\medskip

Let $\Fl$ denote the flag variety $G/B$, and let $\Fl_w$ (resp., $\ol\Fl_w$) denote
the Schubert cell, corresponding to $w$ (resp., its closure). Let us denote the
corresponding locally closed substack in $(\Bun_B\underset{\Bun_G}\times \Bun_B)^w$
(resp., $(\Bun^\cmu_B\underset{\Bun_G}\times \Bun_B^{\cmu'})^w$) by $\CZ_w$
(resp., $\CZ^{\cmu,\cmu'}_w$). For $w=w_0$ this stack is closely related
to the Zastava spaces $\CZ^\cmu$ introduced in \secref{intr Zastava}.

\medskip

By definition, the stack $\CZ_w$ classifies the data of pairs $(\CF_B,\sigma)$,
where $\CF_B$ is a $B$-torsor on $X$, and $\sigma$ is a section of
$\CF_B\overset{B}\times \ol\Fl_w$, such that, over the generic point of the
curve, $\sigma$ hits $\CF_B\overset{B}\times \Fl_w$. Let us call the locus of
$X$ where $\sigma$ does {\it not } hit $\CF_B\overset{B}\times \Fl_w$ the
"locus of degeneration" of $(\CF_B,\sigma)$. We shall now construct a map
$\pi^{\cmu,\cmu',w}:\CZ^{\cmu,\cmu'}_w\to X^{w(\cmu')-\cmu}$, such that
for $(\CF_B,\sigma)$ as above, the support of the divisor
$\pi^{\cmu,\cmu',w}(\CF_B,\sigma)$ is its locus of degeneration.

\medskip

For each irreducible $G$-module $V^\lambda$, let $V^{\lambda,\geq w}$ 
be the canonical $B$-stable subspace, spanned by vectors with weights
$\geq w(\lambda)$. Let $V^{\lambda,> w}\subset V^{\lambda,\geq w}$
be the codimension-$1$ subspace, spanned by vectors with weights
$>w(\lambda)$. Then, a point of $\Fl$, which gives a line $\fl^\lambda$ in each
$V^\lambda$, belongs to $\ol\Fl_w$ if and only if $\fl^\lambda\in V^{\lambda,\geq w}$
for each $\lambda\in \Lambda^+$, and it belongs to $\Fl_w$ if the image of
$\fl^\lambda$ in $V^{\lambda,\geq w}/V^{\lambda,> w}$ is non-zero for every
$\lambda$.

Similarly, a data of $(\CF_B,\sigma)$, where $\sigma$ is a section of
$\CF_B\overset{B}\times \Fl$ defines a sub-bundle for every $\lambda\in \Lambda^+$
$$\CL^\lambda \hookrightarrow V^\lambda_{\CF_B},$$
and such a point belongs to $\CZ_w$ if and only if $\CL^\lambda$ belongs to
$V^{\lambda,\geq w}_{\CF_B}$ and the composed map of line bundles
\begin{equation} \label{defect diagram}
\CL^\lambda\to V^{\lambda,\geq w}_{\CF_B}\to 
\left(V^{\lambda,\geq w}/V^{\lambda,> w}\right)_{\CF_B}
\end{equation}
is non-zero. Moreover, the locus of degeneration of $(\CF_B,\sigma)$ is
where the maps \eqref{defect diagram} have zeroes.

If $(\CF_B,\sigma)$ belongs to the component $\CZ^{\cmu,\cmu'}_w$ of
$\CZ_w$, the degree of the line bundle $\CL^\lambda$ is by definition 
$\langle \lambda,\cmu\rangle$
and that of $\left(V^{\lambda,\geq w}/V^{\lambda,> w}\right)_{\CF_B}$ equals
$\langle \lambda,w(\cmu')\rangle$. The divisors of zeroes of the maps
\eqref{defect diagram} for all $\lambda\in \Lambda^+$ are encoded by a point
of $X^{w(\cmu')-\cmu}$. This point is, by definition, the sought-for 
$\pi^{\cmu,\cmu',w}(\CF_B,\sigma)$.

\medskip

Let us denote by $\fq^\cmu$ (resp., $\fq^{\cmu'}$) the natural map from 
$\CZ^{\cmu,\cmu'}_w$ to $\Bun_T^\cmu$ (resp., $\Bun_T^{\cmu'}$). Note that
$\fq^{\cmu'}$ equals the composition 
$$\CZ^{\cmu,\cmu'}_w\overset{\fq^\cmu\times \pi^{\cmu,\cmu',w}}\to 
\Bun_T^\cmu\times X^{w(\cmu')-\cmu}
\overset{\on{id}\times \on{AJ}}\longrightarrow \Bun_T^\cmu\times \Bun_T^{w(\cmu')-\cmu} \to
\Bun_T^{w(\cmu')}\overset{w^{-1}}\to \Bun^{\cmu'}_T.$$

Applying the Verdier duality and the projection formula, we obtain that in order
to show that $\on{CT}^\cmu(\Eis_*^{\cmu'}(E_\cT))^w$ has the desired form, it
is sufficient to prove the following:

\begin{prop}  \label{constant along Bun T}
The direct image with compact supports of the constant sheaf along the map 
$$\fq^\cmu\times \pi^{\cmu,\cmu',w}:\CZ^{\cmu,\cmu'}_w\to \Bun_T^{\cmu}\times X^{w(\cmu')-\cmu}$$
is an extension of complexes, each of which is a pull-back of a complex
on the second multiple (i.e., $X^{w(\cmu')-\cmu}$).
\end{prop}

\ssec{}

The rest of this section is devoted to the proof of \propref{constant along Bun T}.
Fix a point $x\in X$, and consider the open substack
$$\Bun_T^{w(\cmu)}\times (X-x)^{w(\cmu')-\cmu}\subset
\Bun_T^{w(\cmu)}\times X^{w(\cmu')-\cmu},$$
and let us analyze the restriction to it of the complex
$(\fq^\cmu\times \pi^{\cmu,\cmu',w})_!(\BC_{\CZ^{\cmu,\cmu'}_w})$.

Let us choose a subgroup $N'\subset N$, normalized by $T$,
in such a way that its action on $\Fl_w$ is simply-transitive.
Let $B'=N'\cdot T$ be the corresponding subgroup of $B$.
Consider the sheaves of groups $\on{Maps}(X,B')\subset \on{Maps}(X,B)$.

By conjugating $B(\wh\CO_x)$ inside $B(\wh\CK_x)$ by an element
of $T(\wh\CK_x)$, \footnote{Here $\wh\CO_x$ (resp., $\wh\CK_x$) is the
local ring (resp., field) corresponding to the point $x\in X$.}
we can modify both these sheaves at the point $x$, and
obtain sheaves of groups $\wt{\on{Maps}}(X,B')\subset \wt{\on{Maps}}(X,B)$
in such a way that the forgetful map
$$\Bun^\cmu_B\underset{\left(\wt{\on{Maps}}(X,B)\text{-tors}\right)}\times
\left(\wt{\on{Maps}}(X,B')\text{-tors}\right)\to \Bun^\cmu_B$$
is a smooth generalized affine fibration.

Let us denote by $\wt\CZ^{\cmu,\cmu'}_w$ the Cartesian product
$$(X-x)^{w(\cmu')-\cmu}\underset{X^{w(\cmu')-\cmu}}\times
\CZ^{\cmu,\cmu'}_w\underset{\Bun^\cmu_B}\times 
\left(\Bun^\cmu_B\underset{\left(\wt{\on{Maps}}(X,B)\text{-tors}\right)}\times
\left(\wt{\on{Maps}}(X,B')\text{-tors}\right)\right),$$
and by $\wt\fq^\cmu$ (resp, $\wt\pi^{\cmu,\cmu',w}$) its map to $\Bun_T^\cmu$
(resp., $(X-x)^{w(\cmu')-\cmu}$). It suffices to analyze the
object $(\wt\fq^\cmu\times \wt\pi^{\cmu,\cmu',w})_!(\BC_{\wt\CZ^{\cmu,\cmu'}_w})$.

\medskip

We will show that after a suitable base change, 
in the map 
$$\wt\fq^\cmu\times \wt\pi^{\cmu,\cmu',w}:\wt\CZ^{\cmu,\cmu'}_w\to
\left(\Bun_T^{\cmu}\times (X-x)^{w(\cmu')-\cmu}\right),$$
$\Bun_T^{\cmu}$ splits off a direct factor.

More precisely, we shall exhibit a smooth morphism with connected fibers
(in fact, a principal bundle with respect to a connected group-scheme)
$$\left(\Bun_T^{\cmu}\times (X-x)^{w(\cmu')-\cmu}\right)'\to
\left(\Bun_T^{\cmu}\times (X-x)^{w(\cmu')-\cmu}\right)$$
and a stack $\CW^{\cmu,\cmu'}_w$ over $(X-x)^{w(\cmu')-\cmu}$
such that there exists a Cartesian diagram
$$
\CD
\left(\Bun_T^{\cmu}\times (X-x)^{w(\cmu')-\cmu}\right)'
\underset{\left(\Bun_T^{\cmu}\times (X-x)^{w(\cmu')-\cmu}\right)}\times
\wt\CZ^{\cmu,\cmu'}_w  @>>> \CW^{\cmu,\cmu'}_w \\
@VVV    @VVV  \\
\left(\Bun_T^{\cmu}\times (X-x)^{w(\cmu')-\cmu}\right)' @>>>  
(X-x)^{w(\cmu')-\cmu}.
\endCD
$$

Namely, we take $\left(\Bun_T^{\cmu}\times (X-x)^{w(\cmu')-\cmu}\right)'$
to be the stack, whose fiber over a point 
$$(\CF_T,D)\in \left(\Bun_T^{\cmu}\times (X-x)^{w(\cmu')-\cmu}\right)$$
is the scheme of (sufficiently high) level structures of the $T$-torsor $\CF_T$
over the support of $D$.

The sought-for scheme $\CW^{\cmu,\cmu'}_w$ is constructed as follows.
It classifies triples $(\CF_{N'},D,\sigma')$, where:

\begin{itemize}

\item $\CF_{N'}$ is a principal
$N'$-bundle on $X-x$, 

\item
$D$ is a point of $(X-x)^{w(\cmu')-\cmu}$, and 

\item
$\sigma'$ is a section of $\CF_{N'}\overset{N'}\times \ol\Fl_w$, such that
the resulting divisor (obtained as in \eqref{defect diagram}) equals $D$.

\end{itemize}

In other words, $\CW^{\cmu,\cmu'}_w$ is the scheme, classifying maps
from $X-x$ to the stack $\ol\Fl_w/N'$ of a given degree (cf. \cite{BFG}, Sect. 2.16).

\section{Proof of \thmref{deformation base}}   \label{proof of deformation base}

We will present a short argument, pointed out to us by A.~Beilinson, that
proves the theorem and simultaneously makes the assertion of 
\lemref{fibers over deformation} manifest.

\ssec{}
Let us construct the map
\begin{equation} \label{des map alg}
\wh\CO_{\DefbE}\to \wh{R}_{E_\cT},
\end{equation}
i.e., a map
\begin{equation} \label{des map}
\on{Spf}(\wh{R}_{E_\cT})\to \DefbE.
\end{equation}

We can represent the topological algebra $\wh{R}_{E_\cT}$ as
$\underset{\underset{n\in \BZ^{\geq 0}}{\longleftarrow}}{lim}\, R^n_{E_\cT}$,
where $$R^n_{E_\cT}:=\underset{\cla\in \cLambda^{pos},|\cla|=n}\oplus\, 
\wh{R}^{-\cla}_{E_\cT}.$$

Thus, we need to exhibit a compatible family of $R^n_{E_\cT}$-points of
the functor $\DefbE$. This means that for each $n$ we need to exhibit
a tensor functor
$$\sF_n:\{\cB\mod\}\to \{\text{local systems of }R^n_{E_\cT}\text{-modules over }X\},$$
with identifications 
$$\sF_n(V)\underset{R^n_{E_\cT}}\otimes \BC\simeq V_{E_\cT} \text{ for }
V\in \cB\mod$$ and 
$$\sF_n(V')\simeq R^n_{E_\cT}\otimes V'_{E_\cT} \text{ for }
V'\in \cT\mod\subset \cB\mod.$$

\ssec{}

The sought-for functor $\sF_n$ is constructed as follows. First, for $V\in \cB\mod$
and $\cnu\in \cLambda$, let $\Omega(\cn_{X,E_\cT},V_{E_\cT})$ be the relative
version of $\Omega(\cn_{X,E_\cT},V_{E_{\cT,x}})$, which is a perverse sheaf
on the corresponding ind-scheme $_\infty(X\times X^\cnu)$. Let $R(V)^\cnu_{E_\cT}$
be its direct image under the natural projection 
$$_\infty(X\times X^\cnu)\to X.$$

\medskip

We have the natural maps $\wh{R}^{-\cla}_{E_\cT}\otimes R(V)^\cnu_{E_\cT}\to
R(V)^{\cnu-\cla}_{E_\cT}$. Consider $R(V)_{E_\cT}:=\underset{\cnu}\oplus\, R(V)^\cnu_{E_\cT}$ 
as a $\cLambda$-graded local system on $X$; it is acted on by $R_{E_\cT}$. Set
$$R(V)^n_{E_\cT}:=R(V)_{E_\cT}\underset{R_{E_\cT}}\otimes R^n_{E_\cT};$$
this is an $R^n_{E_\cT}$-local system on $X$ and
$V\mapsto R(V)^n_{E_\cT}$ provides the desired functor $\sF_n$.

\medskip

The functors $\sF_n$ are evidently compatible for different $n$, and hence we obtain
the morphism in \eqref{des map}. Note that by definition, for $V\in \cB\mod$,
the pull-back of
$V_{E_{\cB,\wh{univ},x}}$ under this morphism identifies, by construction, with
$\wh{R}(V_x)_{E_\cT}$, i.e., \lemref{fibers over deformation} holds. Thus, it remains
to show that the above map \eqref{des map} is an isomorphism.

\ssec{}

The regularity assumption on $E_\cT$ implies that the deformation theory of $E_\cB$
is unobstructed. Hence, $\wh\CO_{\DefbE}$ is (non-canonically) isomorphic to the completion of
a polynomial algebra. Likewise, by \secref{discuss omega}, the algebra $R_{E_\cT}$ is regular,
and hence $\wh{R}_{E_\cT}$ is also (non-canonically) isomorphic to the 
completion of a polynomial algebra. 
Therefore, in order to show that the map \eqref{des map} is an isomorphism, it is sufficient
to do so at the level of tangent spaces at the maximal ideal. 

\medskip

From \secref{discuss omega}, it follows that $\wh{R}_{E_\cT}/\fm^2_{\wh{R}_{E_\cT}}$
identifies with 
$$\BC\oplus \underset{\alpha\in \Delta^+}\bigoplus\, H^1(X,\fn^{\check\alpha}_{X,E_\cT})^*
\simeq \BC\oplus  H^1(X,\fn_{X,E_\cT})^*,$$
where $\fm_{\wh{R}_{E_\cT}}$ denotes the maximal ideal in $\wh{R}_{E_\cT}$.
I.e., the tangent space to $\on{Spf}(\wh{R}_{E_\cT})$ at its closed point identifies
canonically with $H^1(X,\fn_{X,E_\cT})$.

In addition, by standard deformation theory, the tangent space to $\DefbE$ at its closed point
also identifies with $H^1(X,\fn_{X,E_\cT})$. Moreover, it is easy to see from the
construction of the map \eqref{des map} that it induces the identity map on
$H^1(X,\fn_{X,E_\cT})$, implying our assertion.

\end{document}